\documentclass[11pt]{article}
\usepackage[english]{babel}

\usepackage{latexsym}
\usepackage{amsmath}
\usepackage{amsfonts}
\usepackage{mathrsfs}
\usepackage{amstext}
\usepackage{amssymb}
\usepackage{amsthm}       
\usepackage{hyperref}

\usepackage{todonotes}

\usepackage{xspace}
\usepackage{xcolor}

\newcommand\cbk{\color{black}}

\newcommand\R{\ensuremath{\mathbb{R}}}

\newcommand\N{\ensuremath{\mathbb{N}}}
\newcommand\eps{\ensuremath{\varepsilon}}

\allowdisplaybreaks

\newcommand{\cE}{\mathcal E}

\newtheorem{theorem}{Theorem}[section]
\newtheorem{prop}[theorem]{Proposition}

\newtheorem{rem}[theorem]{Remark}
\newtheorem{lem}[theorem]{Lemma}

\usepackage[margin=1.2in]{geometry}

\usepackage{graphicx, fullpage}	

\begin{document}

\title{Asymptotic behavior of $u$-capacities and singular perturbations for the Dirichlet-Laplacian}
\author{\ }

\author{
Laura Abatangelo\thanks{Dipartimento di Matematica e Applicazioni, Universit\`a degli Studi di Milano-Bicocca,Via Cozzi 55, 20125 Milano, Italy \texttt{laura.abatangelo@unimib.it}} , 
Virginie Bonnaillie-No\"el\thanks{D\'epartement de math\'ematiques et applications,
\'Ecole normale sup\'erieure, CNRS, PSL University, 45 rue d'Ulm, 75005 Paris, France
\texttt{bonnaillie@math.cnrs.fr}} , 
Corentin L\'ena\thanks{Matematiska institutionen, Stockholms universitet, SE-106 91 Stockholm \texttt{corentin@math.su.se}} , 
and Paolo Musolino\thanks{Dipartimento di Matematica `Tullio Levi-Civita', Universit\`a degli Studi di Padova, via Trieste 63, 35121 Padova, Italy \texttt{musolino@math.unipd.it}}}


\date{November 13, 2019}

\maketitle

\noindent
{\bf Abstract.} In this paper we study the asymptotic behavior of  $u$-capacities of small sets and its application to the analysis of the eigenvalues of the Dirichlet-Laplacian on a bounded {planar} domain with a small hole. More precisely, we consider two (sufficiently regular) bounded open connected sets $\Omega$ and $\omega$ of $\mathbb{R}^2$, containing the origin. First, if $\varepsilon$ is positive and small enough and if $u$ is a function defined on $\Omega$, we compute an asymptotic expansion of the $u$-capacity $\mathrm{Cap}_\Omega(\varepsilon \omega, u)$ as $\varepsilon \to 0$. As a byproduct, we compute an asymptotic expansion for the $N$-th eigenvalues of the Dirichlet-Laplacian in the perforated set $\Omega \setminus (\varepsilon \overline{\omega})$ for $\varepsilon$ close to $0$. Such formula shows explicitly the dependence of the asymptotic expansion on the behavior of the corresponding eigenfunction near $0$ and on the shape $\omega$ of the hole.

\vspace{11pt}


\vspace{11pt}


\section{Introduction}

This paper deals with the asymptotic behavior of the so-called $u$-capacities of small sets and its application to the analysis of the eigenvalues of the Dirichlet-Laplacian  on a bounded domain with a small hole. 

The dependence of the spectrum of the Laplace operator upon regular and singular domain perturbations has been long investigated by several authors with many different techniques. A fundamental tool in the analysis of the eigenvalues of the Dirichlet-Laplacian upon domain perturbation has revealed to be the so-called {\it (condenser) capacity}.

So, if we consider a bounded, connected open set $\Omega$   of $\mathbb{R}^2$, then for every compact subset $K$ of $\Omega$, the {\it (condenser) capacity} of $K$ in $\Omega$ is defined as
\begin{equation} \label{eq:cap}
\mathrm{Cap}_{\Omega}(K)\equiv \mathrm{inf} \bigg\{\int_{\Omega} |\nabla f|^2\, dx \colon f \in H^1_0(\Omega)\ \mathrm{and}\ f-\eta_K \in H^1_0(\Omega \setminus K)\bigg\}\, ,
\end{equation}
where $\eta_K$ is a fixed smooth function such that $\mathrm{supp}\, \eta_K \subseteq \Omega$ and $\eta_K \equiv 1$ in a neighborhood of $K$. The infimum in \eqref{eq:cap} is achieved by a function $V_K \in H^1_0(\Omega)$ such that $V_K-\eta_K \in H^1_0(\Omega\setminus K)$ so that
\[
\mathrm{Cap}_\Omega(K)=\int_{\Omega}|\nabla V_K|^2 \, dx\, ,
\]
where $V_K$ ({\it capacitary potential}) is the unique solution of the Dirichlet problem
\begin{equation}\label{eq:dircapK}
\left\{
\begin{array}{ll}
\Delta V_K=0&\text{ in }\Omega \setminus K\,,\\
V_K=0&\text{on } \partial\Omega\,,\\
V_K=1 &\text{ on }K\,.
\end{array}
\right.
\end{equation}
By saying that $V_K$ solves \eqref{eq:dircapK} we mean that $V_K \in H_0^1(\Omega)$, $V_K - \eta_K \in H_0^1(\Omega \setminus K)$, and
\[
\int_{\Omega \setminus K}\nabla V_K \cdot \nabla \phi\, dx=0 \qquad \forall\phi \in H_0^1(\Omega\setminus K).
\]
Moreover, if $\Omega$ and $K$ are sufficiently regular (for example Lipschitz), one can interpret the boundary conditions of problem \eqref{eq:dircapK} in the trace sense (cf., {\it e.g.}, Costabel \cite{Co88}). 

It is well-known that the spectrum of the Dirichlet-Laplacian on the bounded domain $\Omega$ does not change if we remove a compact subset $K$ of zero capacity (cf., {\it e.g.}, Rauch and Taylor \cite{RaTa75}). {If we denote by
\[
0<\lambda_1(\Omega)<\lambda_2(\Omega)\leq \dots \leq \lambda_N(\Omega)\leq \dots
\]
and
\[
0<\lambda_1(\Omega \setminus K)<\lambda_2(\Omega \setminus K)\leq \dots \leq \lambda_N(\Omega \setminus K)\leq \dots
\]
the sequences of the eigenvalues of the Dirichlet-Laplacian in $\Omega$ and in $\Omega \setminus K$, respectively, then Rauch and Taylor \cite{RaTa75} also proved that} the $N$-th eigenvalue $\lambda_N(\Omega\setminus K)$ of the Dirichlet-Laplacian in $\Omega \setminus K$ is close to $\lambda_N(\Omega)$ if and only if  the capacity $\mathrm{Cap}_\Omega(K)$ of $K$ in $\Omega$ is small. 

The result by Rauch and Taylor \cite{RaTa75} can be seen as a continuity result for the eigenvalues with respect to the capacity. On the other hand, a higher regularity result holds. Indeed, Courtois \cite{Co95} has investigated the behavior of the spectrum of the Dirichlet-Laplacian  in $X\setminus A$, where $X$ is a closed Riemannian manifold and $A$ a ``small'' compact subset. In particular, he has shown that if $K\subseteq \Omega$ is compact and $\mathrm{Cap}_\Omega(K)$ is small then the function
\begin{equation}\label{eq:difflambda}
\lambda_N(\Omega \setminus K)-\lambda_N(\Omega)
\end{equation}
is differentiable with respect to $\mathrm{Cap}_\Omega(K)$. Therefore, one can obtain asymptotic expansions for the difference in \eqref{eq:difflambda} in terms of the capacity $\mathrm{Cap}_\Omega(K)$. 

If, for example, $\lambda_N(\Omega)$ is simple, then in order to obtain more refined asymptotic expansions of the difference $\lambda_N(\Omega \setminus K)-\lambda_N(\Omega)$, one can take into account also the behavior of the corresponding eigenfunction $u_N$. More precisely, one can replace the capacity $\mathrm{Cap}_\Omega(K)$ by the so-called $u_N$-capacity $\mathrm{Cap}_{\Omega}(K,u_N)$.

So, if $u$ is a function in $H^1_0(\Omega)$, we introduce the {\it $u$-capacity} by setting
\begin{equation} \label{eq:ucap}
\mathrm{Cap}_{\Omega}(K,u)\equiv \mathrm{inf} \bigg\{\int_{\Omega} |\nabla f|^2\, dx \colon f \in H^1_0(\Omega)\ \mathrm{and}\ f-u \in H^1_0(\Omega \setminus K)\bigg\}\, .
\end{equation}
The infimum in \eqref{eq:ucap} is achieved by a function $V_{K,u}$ which is the unique solution of  the Dirichlet problem
\begin{equation}\label{eq:dircapKu}
\left\{
\begin{array}{ll}
\Delta V_{K,u}=0&\text{ in }\Omega \setminus K\,,\\
V_{K,u}=0&\text{on } \partial\Omega\,,\\
V_{K,u}=u &\text{ on }K\,,
\end{array}
\right.
\end{equation}
so that
\[
\mathrm{Cap}_\Omega(K,u)=\int_{\Omega}|\nabla V_{K,u}|^2 \, dx\, .
\]
As above, by saying that $V_{K,u}$ solves \eqref{eq:dircapKu} we mean that $V_{K,u} \in H_0^1(\Omega)$, $V_K - u \in H_0^1(\Omega \setminus K)$, and
\[
\int_{\Omega \setminus K}\nabla V_{K,u} \cdot \nabla \phi\, dx=0 \qquad \forall\phi \in H_0^1(\Omega\setminus K).
\]

Definition \eqref{eq:ucap} can be extended to $H^1(\Omega)$ functions, by setting, for any $u \in H^1(\Omega)$, $\mathrm{Cap}_\Omega(K,u) \equiv \mathrm{Cap}_\Omega(K,\eta_K u)$ being $\eta_K$ as in \eqref{eq:cap}.

Such an object can be successfully employed in order to compute asymptotic expansions of (simple) eigenvalues. Indeed, the following result holds (cf.  Courtois \cite[Proof of Theorem 1.2]{Co95} and Abatangelo, Felli, Hillairet, and Lena \cite[Theorem 1.4]{AbFeHiLe}).

\begin{theorem}\label{thm:asy:eig} 
Let $\lambda_N(\Omega)$ be a simple eigenvalue of the Dirichlet-Laplacian in a bounded, connected, and open set $\Omega$. Let $u_N$ be a $L^2(\Omega)$-normalized eigenfunction associated to $\lambda_N(\Omega)$ and let $(K_\varepsilon)_{\varepsilon>0}$ be a family of compact sets contained in $\Omega$ concentrating to a compact set $K$ with $\mathrm{Cap}_\Omega(K) = 0$. Then
\begin{equation}\label{eq:asy:eig}
\lambda_N(\Omega \setminus K_\varepsilon)=\lambda_N(\Omega)+\mathrm{Cap}_\Omega(K_\varepsilon,u_N)+o(\mathrm{Cap}_\Omega(K_\varepsilon,u_N))\, , \qquad \text{as $\varepsilon \to 0$}\, .
\end{equation}
\end{theorem}

The aim of this paper is twofold. On one hand, we wish to investigate the asymptotic behavior of $\mathrm{Cap}_\Omega(K_\varepsilon,u)$ as $\varepsilon \to 0$, when $K_\varepsilon=\varepsilon \overline{\omega}$ (with $\omega$ a sufficiently regular open set) and $u$ a generic function. On the other hand, we want to combine such asymptotic analysis with the formula \eqref{eq:asy:eig} of Theorem \ref{thm:asy:eig} and obtain asymptotic expansions of $\lambda_N(\Omega \setminus (\varepsilon \overline{\omega}))$ where the dependence both on  the structure of the normalized eigenfunction $u_{N}$ around $0$ and on the geometry of $\omega$ is explicit.  We emphasize that in our case, the limit compact $K$ of Theorem \ref{thm:asy:eig} consists of just one point, namely $\{0\}$, and therefore the corresponding capacity is equal to zero.

\subsection{Asymptotic behavior of $u$-capacities}

We will be working in the frame of Schauder classes and thus, in order to introduce the geometric setting of the paper, we now fix 
\[
\alpha\in]0,1[\,, 
\] 
 and we assume that
\begin{equation}\label{e1}
\begin{split}
&\text{$\Omega$ and $\omega$ are open bounded connected subsets of $\mathbb{R}^2$ of}\\
&\text{class $C^{1,\alpha}$ such that $\mathbb{R}^2\setminus\overline{\Omega}$  and $\mathbb{R}^2\setminus\overline{\omega}$ are connected,}\\
&\text{and  such that the origin  $0$ of $\mathbb{R}^2$ belongs both to $\Omega$ and $\omega$.}
\end{split}
\end{equation} 
 For the definition of functions and sets of the Schauder classes $C^{0,\alpha}$ and $C^{1,\alpha}$ we refer for example to Gilbarg and Trudinger~\cite[\S6.2]{GiTr01}.  Condition \eqref{e1} implies that there exists a real number $\varepsilon_0$ such that
\[
\varepsilon_0>0\ \mathrm{and\ }\ \varepsilon\overline{\omega}\subseteq \Omega\ \mathrm{for\ all}\ \varepsilon\in]-\varepsilon_0,\varepsilon_0[\,.
\]  
Then we denote by $\Omega_\varepsilon$ the perforated domain defined by
\[
\Omega_\varepsilon\equiv\Omega\setminus(\varepsilon\overline{\omega})\qquad\quad\forall\varepsilon\in]-\varepsilon_0,\varepsilon_0[\,.
\]
Clearly, $\Omega_\varepsilon$ is an open bounded connected subset of $\mathbb{R}^2$ of class $C^{1,\alpha}$ for all $\varepsilon\in]-\varepsilon_0,\varepsilon_{0}[\setminus\{0\}$. Moreover,  the boundary
$\partial \Omega_\varepsilon$ of $\Omega_\varepsilon$ is the union of the two connected components $\partial \Omega$ and $\partial (\varepsilon \omega)=\varepsilon\partial\omega$, for all
$\varepsilon\in]-\varepsilon_0,\varepsilon_{0}[$. We also note that $\Omega_0=\Omega\setminus\{0\}$.\par
Then we assume that
\begin{equation}\label{eq:assf}
\text{$u \in H^1(\Omega)$ is analytic in a neighborhood of $0$.}
\end{equation}
We are interested in studying the behavior of $\mathrm{Cap}_\Omega(\varepsilon\omega,u)$ as $\varepsilon$ tends to $0$. More precisely, our aim is to obtain as much accurate and constructive as possible expansions for $\mathrm{Cap}_\Omega(\varepsilon\omega,u)$ in terms of the parameter $\varepsilon$. Moreover, besides the dependence on $\varepsilon$, we want to highlight the effect of the geometry of the problem ({\it i.e.},  $\Omega$ and $\omega$) and of the function $u$ on  $\mathrm{Cap}_\Omega(\varepsilon\omega,u)$. 

As we shall see, to reach this goal, one can try to follow different strategies.

\subsubsection{Asymptotic behavior of the capacity and conformal mappings}

A standard method to convert a boundary value problem for the Laplace equation defined in a generic Jordan domain with one hole into a (possibly) easier one is to exploit conformal mapping theory. In this way, one can transform the original problem into a boundary value problem defined in an annular domain of the type $\mathbb{A}(r,1)\equiv \{z \in \mathbb{C}\colon r < |z|< 1\}$ for some $r>0$. Then one can try to find explicit formulas for the solution of the transformed problem in $\mathbb{A}(r,1)$ and finally to exploit those formulas for the representation of the solution of the original problem. Clearly, an approach of this type can be applied also to the computation of the capacity $\mathrm{Cap}_{\Omega}(\varepsilon \omega)$, since it is defined by means of the solution of a Dirichlet problem for the Laplacian with locally constant boundary data.

So we identify $\mathbb{R}^2$ with the complex plane $\mathbb{C}$ and we assume that $\partial \Omega$ and $\partial \omega$ are the image of two simple closed  curves $\zeta^o$ and $\zeta^i$ of class $C^{1,\alpha}$ from the boundary $\partial \mathbb{D}$ of the unit disk $\mathbb{D}$ to $\mathbb{C}$. By the Riemann Mapping Theorem, one deduces that for each $\varepsilon \in ]-\varepsilon_0,\varepsilon_0[\setminus \{0\}$ there exist a unique $r[\varepsilon] \in ]0,1[$ and a unique holomorphic homeomorphism $g[\varepsilon]$ from the set $\mathbb{A}(r[\varepsilon],1)$ onto $\Omega_\varepsilon$ such that the map $g[\varepsilon]$ can be extended to an element of class $C^{1,\alpha}(\overline{\mathbb{A}(r[\varepsilon],1)},\mathbb{C})$ (which we still denote by $g[\varepsilon]$) and such that $g[\varepsilon](1)=\zeta^o(1)$ (cf. Lanza de Cristoforis and Rogosin \cite[Thm.~3.1]{LaRo01}). Moreover, $g[\varepsilon]$ is a homeomorphism of $\overline{\mathbb{A}(r[\varepsilon],1)}$ onto $\overline{\Omega_\varepsilon}$. Now we observe that if we compose the solution of 
\begin{equation}\label{eq:dircapan}
\left\{
\begin{array}{ll}
\Delta \mathcal{V}_{r[\varepsilon]}=0&\text{ in }\mathbb{A}(r[\varepsilon],1)\,,\\
\mathcal{V}_{r[\varepsilon]}=0&\text{on } \partial\mathbb{D}\,,\\
\mathcal{V}_{r[\varepsilon]}=1 &\text{ on }r[\varepsilon] \partial \mathbb{D}\,.
\end{array}
\right.
\end{equation}
with the map $g^{(-1)}[\varepsilon]$ we obtain an harmonic function in $\Omega_\varepsilon$, vanishing on $\zeta^o(\partial \mathbb{D})=\partial \Omega$ and equal to $1$ on $\varepsilon\zeta^i(\partial \mathbb{D})=\varepsilon \partial \omega$.  On the other hand, by a direct computation one verifies that the solution of problem \eqref{eq:dircapan} is delivered by
\[
\mathcal{V}_{r[\varepsilon]}(z)\equiv \frac{\log |z|}{\log r[\varepsilon]} \qquad \forall z \in \overline{\mathbb{A}(r[\varepsilon],1)}\, .
\]
As a consequence,  the capacitary potential $V_{\varepsilon \omega}$ can be represented as
\[
V_{\varepsilon \omega}(z)\equiv \frac{\log |g^{(-1)}[\varepsilon](z)|}{\log r[\varepsilon]} \qquad \forall z \in \overline{\Omega_\varepsilon}\, .
\]
Then one obtains the following formula for the capacity $\mathrm{Cap}_{\Omega}(\varepsilon \omega)$
\[
\mathrm{Cap}_{\Omega}(\varepsilon \omega)=-\frac{2\pi }{\log r[\varepsilon]} \qquad \forall \varepsilon \in ]-\varepsilon_0,\varepsilon_0[\setminus\{0\}\, .
\]
Therefore, if we want to understand the behavior of $\mathrm{Cap}_{\Omega}(\varepsilon \omega)$ as $\varepsilon \to 0$ we need to investigate $r[\varepsilon]$ for $\varepsilon$ close to $0$. On the other hand, by Lanza de Cristoforis \cite{La02,La04}, we know that there exist $\varepsilon_1 \in ]0,\varepsilon_0[$ and a real analytic function $\mathcal{R}$ from $]-\varepsilon_1,\varepsilon_1[$ to $]0,+\infty[$ such that 
\[
r[\varepsilon]=\varepsilon \mathcal{R}[\varepsilon] \qquad \forall \varepsilon \in ]0,\varepsilon_1[\, .
\]
Moreover, $\mathcal{R}[0]>0$ and the term $\mathcal{R}[0]$ depends on the geometry of $\partial \Omega=\zeta^o(\partial \mathbb{D})$ and of $\partial \omega=\zeta^i(\partial \mathbb{D})$. Accordingly, we deduce the formula
\begin{equation}\label{eq:Rancap}
\mathrm{Cap}_{\Omega}(\varepsilon \omega)=-\frac{2\pi }{\log \varepsilon + \log \mathcal{R}[\varepsilon]} \qquad \forall \varepsilon \in ]0,\varepsilon_1[\, .
\end{equation}
Then by \eqref{eq:Rancap} we have that
\[
\mathrm{Cap}_{\Omega}(\varepsilon \omega)=-\frac{1}{\log \varepsilon }\frac{2\pi }{\big(1 + \frac{1}{\log \varepsilon}\log \mathcal{R}[\varepsilon]\big)} \qquad \forall \varepsilon \in ]0,\varepsilon_1[\, .
\]
Hence, there exists a real analytic map $\tilde{\mathcal{R}}$ from a neighborhood of $(0,0)$ in $\mathbb{R}^2$ with values in  $\mathbb{R}$ such that
\[
\mathrm{Cap}_{\Omega}(\varepsilon \omega)=\tilde{\mathcal{R}}\Big[\varepsilon, \frac{1}{\log \varepsilon}\Big] \, ,
\]
for $\varepsilon$ positive and close to $0$.  By the analyticity of $\tilde{\mathcal{R}}$, one immediately deduces that
\begin{equation}\label{eq:serRancap}
\mathrm{Cap}_{\Omega}(\varepsilon \omega)=\sum_{(k,l) \in \mathbb{N}^2} \gamma_{(k,l)}\varepsilon^k \Big(\frac{1}{\log \varepsilon}\Big)^l \, ,
\end{equation}
for $\varepsilon$ positive and small enough, and where the double power series $\sum_{(k,l) \in \mathbb{N}^2} \gamma_{(k,l)}x_1^k x_2^l$ converges absolutely  for $(x_1,x_2)$ in a neighborhood of $(0,0)$. 

Even if one could explicitly deduce from Lanza de Cristoforis \cite{La02,La04} the limiting value $\mathcal{R}[0]$, we emphasize that no attempt has been done so far in order to derive from the real analyticity of $\tilde{\mathcal{R}}$ the exact value of all the coefficients $\gamma_{(k,l)}$ appearing in \eqref{eq:serRancap}. Moreover, if one tries to apply this method for the computation of the $u$-capacity $\mathrm{Cap}_{\Omega}(\varepsilon \omega, u)$, one faces the problem to find an explicit solution of problem \eqref{eq:dircapan} with the third condition replaced by
\[
\mathcal{V}_{r[\varepsilon]}(z)=u(g[\varepsilon]z) \qquad \forall z \in r[\varepsilon] \partial \mathbb{D}\, .
\]
Then clearly such a dependence on $g[\varepsilon]$ and on $u$ of the Dirichlet datum on the hole makes even more involved the computation of the coefficients of the corresponding expansion of the capacity. Therefore, in order to provide an explicit and constructive expansion for $\mathrm{Cap}_{\Omega}(\varepsilon \omega, u)$ we prefer to follow a different strategy, which does not relies on conformal mappings.

\subsubsection{Asymptotic expansion for the capacity}

Boundary value problems in domains with small holes have been largely investigated in the frame of asymptotic analysis. In order to study these problems several asymptotic expansion techniques have been developed: for example,  the method of matching outer and inner asymptotic expansions proposed by Il'in (cf., {\it e.g.}, \cite{Il92}), the compound asymptotic expansion method of Maz'ya, Nazarov, and Plamenevskij~\cite{MaNaPl00i,MaNaPl00ii}, and the asymptotic analysis of Green's kernels in domains with small cavities by mesoscale asymptotic approximations of Maz'ya, Movchan, and Nieves~\cite{MaMoNi13}. 
In Bonnaillie-No\"el and Dambrine \cite{BoDa13} and in Bonnaillie-No\"el, Dambrine, and Lacave \cite{BoDaLa},  the authors have exploited the method of multiscale asymptotic expansions to analyze the two-dimensional Dirichlet-Laplacian in a domain with moderately close small perforations.  The Dirichlet problem in a planar domain with a small hole has received attention also from the numerical point of view. A numerical approach is proposed, {\it e.g.}, in Babu\v{s}ka,  Soane, and Suri \cite{BaSoSu17} and Chesnel and Claeys \cite{ChCl14}. Problems in perforated domains find several applications, as an example, in the frame of shape and topological optimization (cf.~Novotny and Soko\l owsky \cite{NoSo13}) and in inverse problems (cf.~Ammari and Kang \cite{AmKa07} and Ammari, Kang, and Lee \cite{AmKaLe09}).

An asymptotic expansion of the capacity as the hole collapses to a point can be deduced by the analysis of energy integrals in perforated domains that can be found in Maz'ya, Nazarov, and Plamenevskij \cite[\S 8.1]{MaNaPl00i}. In particular, they prove that
\begin{equation}\label{eq:MNPexp}
\mathrm{Cap}_{\Omega}(\varepsilon \omega)=-\frac{2\pi}{\log \varepsilon + 2 \pi \big(H_{(0,0)}+N\big)} +o(\varepsilon^{\delta})\, , \qquad \forall \delta>0\, ,
\end{equation}
for $\varepsilon$ small and positive, where $e^{2\pi N}$ is the {\it logarithmic capacity} (or {\it outer conformal radius}) of $\omega$ and $H_{(0,0)}$ is the value at $x=0$ of the unique harmonic function $h$ in $\Omega$ such that $h(x)=-\log|x|/(2\pi)$  for all $x\in\partial\Omega$. In particular, by combining \eqref{eq:Rancap} and \eqref{eq:MNPexp}, we deduce that
\[
\log \mathcal{R}[0]=2 \pi \big(H_{(0,0)}+N\big)\, .
\]
We also note that expansions for the capacity for the case of several small inclusions can be deduced from the corresponding expansion of the capacitary potential obtained in Maz'ya, Movchan, and Nieves~\cite[\S 3.2.2]{MaMoNi13}. Moreover, one could produce an asymptotic expansion of $\mathrm{Cap}_{\Omega}(\varepsilon \omega)$ in the higher-dimensional case. However, in such a case, the asymptotic behavior would differ from that of \eqref{eq:MNPexp} since  the logarithmic term would not be present in the asymptotic expansion in dimension greater than or equal to three.

Our aim is now two-fold. On the one hand we want to extend the study of the asymptotic behavior of $\mathrm{Cap}_{\Omega}(\varepsilon \omega)$ to the $u$-capacity $\mathrm{Cap}_{\Omega}(\varepsilon \omega, u)$. On the other hand, we want to represent $\mathrm{Cap}_{\Omega}(\varepsilon \omega,u)$ in terms of convergent power series whose coefficients can be explicitly constructed by solving given systems of integral equations on fixed domains (not depending on $\varepsilon$). As we shall see, the computation of higher order terms in the expansion of $\mathrm{Cap}_{\Omega}(\varepsilon \omega,u)$ is necessary if for example $u$ and its derivatives up to a certain order vanish at the origin of $\mathbb{R}^2$.

\subsubsection{The functional analytic approach}

To reach our goal, we adopt the Functional Analytic Approach proposed by Lanza de Cristoforis \cite{La02, La04} for the analysis of singular perturbation problems in perforated domains. This method indeed allows to prove real analyticity properties for the solution of boundary value problems in perforated domains for elliptic equations (see Lanza de Cristoforis \cite{La08} for the Laplace equation) and systems (as the Lam\'e equations in Dalla Riva and Lanza de Cristoforis \cite{DaLa10} and the Stokes system in Dalla Riva \cite{Da13}). Therefore, by this method, one can deduce the possibility to expand the solution or related quantities in convergent power series. Then, to construct these power series, we follow the strategy of \cite{DaMuRo15} and we compute the coefficients in terms of the solutions of recursive systems of integral equations and of quantities related to the data of the problem (such as the unperturbed domain $\Omega$, the inclusion $\omega$, and the function $u$).

We now observe that by assumption \eqref{eq:assf}  on the analyticity of $u$ and by analyticity results for the composition operator (cf.~B\"{o}hme and Tomi~\cite[p.~10]{BoTo73}, 
Henry~\cite[p.~29]{He82}, Valent~\cite[Thm.~5.2, p.~44]{Va88}), possibly shrinking $\varepsilon_0$, there exists a real analytic map $U_\#$ from $]-\varepsilon_0,\varepsilon_0[$ to $C^{1,\alpha}(\partial \omega)$ such that
\begin{equation}\label{eq:Usharp}
u(\varepsilon t)= U_\#[\varepsilon](t) \qquad \forall t \in \partial \omega\, ,\forall \varepsilon \in ]-\varepsilon_0,\varepsilon_0[\, .
\end{equation}
(for the 
definition and properties of analytic maps, we refer to Deimling~\cite[\S 15]{De85}). Then for all $\varepsilon\in]-\varepsilon_{0},\varepsilon_{0}[ \setminus\{0\}$, we denote by $u_\varepsilon$  the unique solution in $C^{1,\alpha}(\overline{\Omega_\varepsilon})$ of the  problem
\begin{equation}\label{eq:direps}
\left\{
\begin{array}{ll}
\Delta u_\varepsilon=0&\text{ in }\Omega_\varepsilon\,,\\
u_\varepsilon(x)=0&\text{ for all }x\in\partial\Omega\,,\\
u_\varepsilon(x)=U_\#[\varepsilon](x/\varepsilon)&\text{ for all }x\in\varepsilon\partial\omega\,.
\end{array}
\right.
\end{equation}
Clearly,
\[
\begin{split}
&V_{\varepsilon \omega, u}(x)= u_\varepsilon(x)\, , \qquad \forall x \in \Omega_\varepsilon\, ,\forall \varepsilon \in ]-\varepsilon_0,\varepsilon_0[\setminus \{0\}\, ,\\
&V_{\varepsilon \omega, u}(x)= u(x)\, , \qquad \forall x \in \varepsilon \overline{\omega}\, ,\forall \varepsilon \in ]-\varepsilon_0,\varepsilon_0[\setminus \{0\}\, {.}
\end{split}
\]
Accordingly, by the Divergence Theorem, we have
\begin{equation}\label{eq:div}
\begin{split}
\mathrm{Cap}_\Omega(\varepsilon\omega,u)&=\int_{\Omega_\varepsilon}|\nabla u_\varepsilon |^2 \, dx + \int_{\varepsilon \omega}|\nabla u|^2 \, dx\\
&= -\int_{\partial (\varepsilon \omega)}\frac{\partial u_\varepsilon}{\partial \nu_{\varepsilon \omega}}u_\varepsilon \, d\sigma+ \varepsilon^2 \int_{\omega}(\nabla u)(\varepsilon t) \cdot (\nabla u)(\varepsilon t) \, dt\\
&= -\int_{\partial \omega}\nu_{\omega}(t)\cdot \nabla \Big(u_\varepsilon(\varepsilon t)\Big) u(\varepsilon t) \, d\sigma_t+ \varepsilon^2 \int_{\omega}(\nabla u)(\varepsilon t) \cdot (\nabla u)(\varepsilon t) \, dt \, ,
\end{split}
\end{equation}
for all $\varepsilon \in ]-\varepsilon_0,\varepsilon_0[\setminus \{0\}$. Here above the symbols $\nu_\omega$ and $\nu_{\varepsilon \omega}$ denote the outward unit normal to $\partial \omega$ and to $\partial (\varepsilon \omega)$, respectively.

As we have mentioned, our goal is to provide a fully constructive and complete asymptotic expansion for $\mathrm{Cap}_\Omega(\varepsilon\omega,u)$ as $\varepsilon \to 0$. In order to do so, we follow the methods developed in \cite{DaMuRo15} for the solution of the Dirichlet problem in a planar perforated domain. However, in \cite{DaMuRo15} the Dirichlet datum on the boundary of the hole $\partial (\varepsilon \omega)$ is given by rescaling a fixed function $g$ defined on $\partial \omega$, {\it i.e.}, by considering the function $g(\cdot/\varepsilon)$. Here, instead, the boundary condition on $\partial (\varepsilon \omega)$ is given by the trace of $u$ on $\partial (\varepsilon \omega)$. Such a trace can be expressed as $u(\varepsilon(x/\varepsilon))$, that is  the rescaling of the $\varepsilon$-dependent function $U_\#[\varepsilon](\cdot)=u(\varepsilon \cdot)$. Thus we will need to take into account also such a dependence. By \eqref{eq:div},  the quantity $\mathrm{Cap}_\Omega(\varepsilon\omega,u)$ can be expressed  as the sum of
\begin{equation}\label{eq:enuom}
\varepsilon^2 \int_{\omega}(\nabla u)(\varepsilon t) \cdot (\nabla u)(\varepsilon t) \, dt 
\end{equation}
and of (the opposite of) the integral on $\partial  \omega$ of the function  
\begin{equation}\label{eq:fun}
t \mapsto \nu_{\omega}(t)\cdot \nabla \Big(u_\varepsilon(\varepsilon t)\Big) u(\varepsilon t)\, .
\end{equation}
By the analyticity of $u$ in a neighborhood of $0$, one can easily show that the term in \eqref{eq:enuom} is  a real analytic function of $\varepsilon$ around $0$ and accordingly it can be expanded in power series of $\varepsilon$. On the other hand, the integral on $\partial \omega$ of the function in \eqref{eq:fun} needs a more careful analysis.  Thus, as a preliminary step, we will need to provide an expansion for the function in \eqref{eq:fun}. Then, by integrating such an expansion, we will be able to obtain the corresponding result for $\mathrm{Cap}_\Omega(\varepsilon\omega,u)$. In particular, under vanishing assumption for $u$, we are able to prove the validity of the following result (cf. Theorem \ref{thm:cepsmseries} and Remark \ref{rem:cepsmseries} below).
\begin{theorem}\label{thm:introcap}
Let assumption \eqref{eq:assf} hold. Assume that there exists $\overline{k} \in \mathbb{N}\setminus \{0\}$ such that
\[
D^\gamma u(0)=0 \qquad \forall |\gamma| <\overline{k}\, ,
\]
{ and that there exists $\beta \in \mathbb{N}^2$ such that $|\beta|=\overline{k}$ and
\[
D^\beta u(0)\neq 0\, .
\]}
Then 
\[
\mathrm{Cap}_{\Omega}(\varepsilon \omega,u)=\varepsilon^{2 \overline{k}}\Bigg(\int_{\mathbb{R}^2 \setminus \overline{\omega}}|\nabla \mathsf{u}_{\overline{k}}|^2\, dt +\int_{\omega} |\nabla u_{\#,\overline{k}}|^2\, dt\Bigg)+o(\varepsilon^{2\overline{k}}) \qquad \text{as }\varepsilon \to 0\, ,
\]
where  $u_{\#,\overline{k}}$ is defined as in Proposition \ref{funk} and $\mathsf{u}_{\overline{k}}$ is the unique solution of problem \eqref{eq:bvp:mathsfu}.
\end{theorem}

As we shall see, the terms $\int_{\mathbb{R}^2 \setminus \overline{\omega}}|\nabla \mathsf{u}_{\overline{k}}|^2\, dt$ and $\int_{\omega} |\nabla u_{\#,\overline{k}}|^2\, dt$ depend  both on the {geometric} properties of the set $\omega$ and on the behavior at $0$ of the function $u$ (but not on $\Omega$). Here we note that {in the case of dimension higher than $3$} one could  expand  $\mathrm{Cap}_{\Omega}(\varepsilon \omega,u)$ as a power series in $\varepsilon$. Indeed, one may show that $\mathrm{Cap}_{\Omega}(\varepsilon \omega,u)$ depends real analytically on $\varepsilon$ (cf. \cite[Thm.~6.2]{La08}) in a neighborhood of $0$. The results of the present paper rely on the asymptotic expansions of \cite{DaMuRo15} for the solution of the Dirichlet problem for the Laplacian in a perforated planar domain. In such a case, one cannot hope to expand  $\mathrm{Cap}_{\Omega}(\varepsilon \omega,u)$ as a power series in $\varepsilon$ since, as is well known, a logarithmic term appears.

\subsection{Asymptotic expansions of the eigenvalues}

The asymptotic behavior of the eigenvalues of the Laplacian in domains with small holes has been long investigated by several authors. 

One of the first contributions is probably due to Samarski\u\i \   \cite{Sa48} that showed that the perturbation of an eigenvalue $\lambda_N$ for the Dirichlet-Laplacian when a small set $\omega_\varepsilon$ is removed from a subset $\Omega$ of $\mathbb{R}^3$ admits the following estimate
\begin{equation}\label{eq:Sam}
\Delta \lambda_N \leq 4\pi \kappa_N^2 \mathrm{Cap}_\Omega(\omega_\varepsilon)+O(\mathrm{Cap}_\Omega(\omega_\varepsilon)^2)\, ,
\end{equation}
where $\kappa_N$ is the maximum value of the $N$-th normalized eigenfunction on $\overline{\omega_\varepsilon}$ (cf. Maz'ya, Nazarov, and Plamenevski\u\i \   \cite{MaNaPl84}).

Later on, in the paper \cite{RaTa75}, Rauch and Taylor  studied the behavior of the eigenvalues and eigenfunctions of the Laplacian in a domain $\Omega$ where a ``thin'' set is removed. A consequence of their (more general) results is that if $\Omega$ and $\omega$ are sufficiently regular bounded open subsets of $\mathbb{R}^n$ containing the origin, and $\lambda_N(\cdot)$ is the $N$-th eigenvalue of Dirichlet-Laplacian then
\begin{equation}\label{eq:conveig}
\lambda_N(\Omega \setminus (\varepsilon \overline{\omega}))\to \lambda_N(\Omega) \qquad \text{as $\varepsilon\to 0^+$}\, .
\end{equation}
In view of the estimate \eqref{eq:Sam} of Samarski\u\i \   \cite{Sa48} and the convergence result \eqref{eq:conveig} of Rauch and Taylor \cite{RaTa75}, many authors have then started to compute asymptotic expansions for the eigenvalues of the Laplacian (under various boundary conditions) in domains with small holes.

For example, Ozawa has devoted a series of papers (cf., {\it e.g.}, \cite{Oz80, Oz81a, Oz81b, Oz82a, Oz82b}) to the computation of asymptotic expansions for the eigenvalues {of} the Laplacian, under many different boundary conditions, when we make a small perforation. In particular, Ozawa has shown in \cite{Oz81b} that
if $n=2$ and $\omega$ is the unit ball $\mathbb{B}_2(0,1)$ then 
\begin{equation}\label{eq:Ozawa}
\lambda_N(\Omega \setminus (\varepsilon \overline{\mathbb{B}_2(0,1)}))=\lambda_N(\Omega)-2\pi(\log \varepsilon)^{-1} (u_N(0))^2+O((\log \varepsilon)^{-2}) \qquad \text{as $\varepsilon \to 0^+$}\, ,
\end{equation}
where $\lambda_N(\Omega)$ is a simple eigenvalue for the Dirichlet-Laplacian in $\Omega$ and $u_N$ a corresponding $L^2(\Omega)$-normalized eigenfunction.

Moreover, Maz'ya, Nazarov, and Plamenevski\u\i \ (see, for example, \cite{MaNaPl84} and \cite[Chapter 9]{MaNaPl00i}) have produced asymptotic expansions of the eigenvalues of boundary value problems for the Laplace operator in domains with small holes. For example, in the three dimensional case, they have shown in \cite{MaNaPl84} that for the first eigenvalue of the Laplacian with Dirichlet condition we have
\begin{equation}\label{eq:Mazya}
\begin{split}
\lambda_1(\Omega \setminus (\varepsilon \overline{\omega}))=&\lambda_1(\Omega)+4\pi \mathrm{Cap}(\omega)(u_1(0))^2 \varepsilon +[4\pi u_1(0)\mathrm{Cap}(\omega)]^2 \\
& \times \Big\{ \Gamma(0)+\frac{u_1(0)}{4\pi}\int_\Omega u_1(x)|x|^{-1}\, dx\Big\}\varepsilon^2+O(\varepsilon^3)\qquad \text{as $\varepsilon \to 0^+$}\, ,
\end{split}
\end{equation}
where $u_1$ is a corresponding $L^2(\Omega)$-normalized eigenfunction in $\Omega$, $\mathrm{Cap}(\omega)$ the harmonic capacity of $\omega$ and $\Gamma$ is a function defined through an auxiliary boundary value problem. We note that since the first eigenfunction $u_1$ does not vanish inside $\Omega$, $u_1(0)\neq 0$ and thus the asymptotic expansion in \eqref{eq:Mazya} is sharp. However, this in general is not the case if we consider different eigenvalues $\lambda_N$. In particular, we note that if we consider the asymptotic expansion of \eqref{eq:Ozawa}, then if the origin belongs to a nodal line of the eigenfunction $u_N$, we have $u_N(0)=0$. Therefore \eqref{eq:Ozawa} reduces to
\begin{equation}\label{eq:Ozawabis}
\lambda_N(\Omega \setminus (\varepsilon \overline{\mathbb{B}_2(0,1)}))=\lambda_N(\Omega)+O((\log \varepsilon)^{-2}) \qquad \text{as $\varepsilon \to 0^+$}\, .
\end{equation}
As a consequence, in view of \eqref{eq:Ozawabis}, in case $u_N(0)=0$ one may need to compute further terms in the asymptotic expansion.

Subsequently, many authors have started to study the behavior of the spectrum of the Laplacian under removal of ``small'' sets in the Riemannian setting. As an example, we mention the works by Besson \cite{Be85}, Chavel \cite{Ch84}, Chavel and Feldman \cite{ChFe88}, Colbois and Courtois \cite{CoCo91}, Courtois \cite{Co95}. 

 As we have already mentioned, one can find in Courtois \cite[Proof of Theorem 1.2]{Co95} and in Abatangelo, Felli, Hillairet, and Lena \cite[Theorem 1.4]{AbFeHiLe} an asymptotic formula for the eigenvalues involving the notion of $u$-capacity (see equation \eqref{eq:asy:eig} of Theorem \ref{thm:asy:eig}). In particular, if $\Omega$ and $\omega$ are as in \eqref{e1} and $\lambda_N(\Omega)$ is simple, this reads as
\begin{equation}\label{eq:Lauraetal}
\lambda_N(\Omega \setminus(\varepsilon\overline{\omega}))=\lambda_N(\Omega)+\mathrm{Cap}_\Omega(\varepsilon\omega,u_N)+o(\mathrm{Cap}_\Omega(\varepsilon\omega,u_N)) \qquad \text{as $\varepsilon \to 0^+$}\, . 
\end{equation}
As a consequence of \eqref{eq:Lauraetal}, in order to find an asymptotic expansion in the parameter $\varepsilon$ we need to compute $\mathrm{Cap}_\Omega(\varepsilon \omega, u_N)$. Abatangelo, Felli, Hillairet, and Lena \cite{AbFeHiLe} have computed such quantity for specific sets as ball and segment and they have shown explicitly the dependence of $\mathrm{Cap}_\Omega(\varepsilon \omega, u_N)$ on the behavior of the eigenfunction $u_N$ around the origin. In particular if $n=2$ and $\omega=\mathbb{B}_2(0,1)$, they have proved that
\[
\mathrm{Cap}_\Omega(\varepsilon \mathbb{B}_2(0,1), u_N)=
\begin{cases}
\frac{2\pi}{|\log \varepsilon|}(u_N(0))^2(1+o(1))\, , & \text{if $\overline{k}=0$}\, ,\\
2k \pi \varepsilon^{2\overline{k}}b^2(1+o(1))\,  ,& \text{if $\overline{k}\geq 1$}\, ,
\end{cases}
\]
where $\overline{k} \in \mathbb{N}$ and $b \in \mathbb{R}\setminus \{0\}$ are such that
\[
r^{-\overline{k}}u(r(\cos t, \sin t))\to b \sin (a-\overline{k} t) \qquad \mbox{in } C^1([0,2\pi])\, ,
\]
as $r \to 0^+$ for some  $a \in [0,\pi[$ .

Here, instead, we wish to emphasize the interaction with the geometry of the hole and the structure of the eigenfunction near $0$. In order to do so, we  confine to the two-dimensional case and we exploit the power series expansion for $\mathrm{Cap}_\Omega(\varepsilon \omega, u)$ of Section \ref{sec5}, with $u=u_N$ and where $\omega$ is a quite general regular open set as in \eqref{e1}.

Under the assumption that the $N$-th eigenvalue $\lambda_N(\Omega)$ for the Dirichlet-Laplacian  is simple,  if  $u_N$ is a $L^2(\Omega)$-normalized eigenfunction related to $\lambda_N(\Omega)$ satisfying some vanishing assumption, we are able to prove the following (cf. Theorem \ref{thm:eig2} below).
\begin{theorem}\label{thm:our}
Let the N-th eigenvalue $\lambda_N(\Omega)$ for the Dirichlet-Laplacian  be simple and let $u_N$ be a $L^2(\Omega)$-normalized eigenfunction related to $\lambda_N(\Omega)$. Assume that there exists $\overline{k} \in \mathbb{N}\setminus \{0\}$ such that $D^\gamma u_N(0)=0$  for all $|\gamma| <\overline{k}$ { and that there exists $\beta \in \mathbb{N}^2$ such that $|\beta|=\overline{k}$ and $D^\beta u(0)\neq 0$.}  Then
\begin{equation}\label{eq:our}
\lambda_N(\Omega\setminus (\varepsilon \overline{\omega}))=\lambda_N(\Omega)+\varepsilon^{2 \overline{k}}C(u_N,\omega)+o(\varepsilon^{2\overline{k}}) \qquad \text{as }\varepsilon \to 0^+\, ,
\end{equation}
where $C(u_N,\omega)$ an explicitly defined positive constant depending on $u_N$ and on $\omega$ (cf.~Section \ref{sec6}).
\end{theorem}

One of the consequences of our asymptotic expansion \eqref{eq:our} of Theorem \ref{thm:our} is that it gives the order of the difference 
\[
\lambda_N(\Omega\setminus (\varepsilon \overline{\omega}))-\lambda_N(\Omega)
\]
for a wide family of holes $\omega$. A second important consequence is that the constant $C(u_N,\omega)$ in  \eqref{eq:our} is explicitly defined in terms of solutions to Dirichlet  problems in $\omega$ and $\mathbb{R}^2 \setminus \overline{\omega}$ for the Laplace equation. The Dirichlet data depend on the Taylor expansion of the normalized eigenfunction at $0$. In particular, it shows the dependence of $C(u_N,\omega)$ both on $u_N$ and $\omega$ and thus provides a starting point on the investigation of `optimal' inclusions $\omega$ for such constant (under different constraints).

We note that in the last years the investigation of this type of problems has been carried out in many different directions. Maz'ya, Movchan, and Nieves have \cite{MaMoNi17} have constructed the asymptotic approximation to the first eigenvalue and corresponding eigenfunctions of Laplace operator inside a domain containing a cloud of small rigid inclusions. Lanza de Cristoforis \cite{La12} has considered a Neumann eigenvalue problem and shown representation formulas in terms of analytic maps and $\log \varepsilon$ (depending on the dimension $n$). Sharp estimates when a ball is removed at the vertex of a sector are contained in Lamberti and Perin \cite{LaPe10}. Henrot \cite{He06} has considered  perforated domains in the frame of extremum problems for eigenvalues of elliptic operators. Finally, Ammari, Kang, and Lee \cite{AmKaLe09} have developed an asymptotic theory for eigenvalue problems under domain perturbations by a method  based on potential theory and on the Gohberg-Sigal theory of meromorphic operator-valued functions.

\subsection{Organization of the paper}

The paper is organized as follows. Section \ref{prel} is a section of preliminaries where we provide an integral equation formulation for the boundary value problem defining the $u$-capacity. In Sections \ref{sec3} and \ref{sec4} we compute power series expansions for some auxiliary functions. Section \ref{sec5} contain our main result on the power series expansion of the $u$-capacity of a small set. In Section \ref{sec6}, we compute the asymptotic expansion of an eigenvalue of the Dirichlet-Laplacian in the perforated domain $\Omega \setminus (\varepsilon \overline{\omega})$ as the size $\varepsilon$ of the hole $\varepsilon \omega$ tends to $0$ and in Section \ref{sec7} we discuss optimal locations of small holes.  Section \ref{sec8} is devoted to some numerical computations on the behavior of the eigenvalues of an ellipse with a small hole and Section \ref{sec9} to their analytical justification.

\section{Preliminaries}\label{prel}

\subsection{Classical notions of potential theory}\label{prel1}

In order to analyze the behavior of the solution to problem \eqref{eq:direps} as $\varepsilon \to 0$ we shall exploit an approach based on potential theory, which allows to convert a boundary value problem into a set of integral equations defined on the boundary of the domain. The method relies on the representation of the solution in terms of some specific integral operators, namely the single and the double layer potentials.

In order to define these operators, we denote by  $S$  the fundamental solution {of $\Delta\equiv \sum_{j=1}^2\partial_{j}^2$} in $\mathbb{R}^2$, that is the function from $\mathbb{R}^2\setminus\{0\}$ to $\mathbb{R}$ defined by
\[
S(x)\equiv\frac{1}{2\pi}\log\,|x|\qquad\forall x\in\mathbb{R}^2\setminus\{0\}\,.
\] 

Now let $\mathcal{O}$ be an open bounded subset of $\mathbb{R}^2$ of class $C^{1,\alpha}$. Let $\phi\in C^{0,\alpha}(\partial\mathcal{O})$.  Then $v[\partial\mathcal{O},\phi]$ denotes the single layer potential with density $\phi$, \textit{i.e.},
\[
v[\partial\mathcal{O},\phi](x)\equiv\int_{\partial\mathcal{O}}\phi(y)S(x-y)\,d\sigma_y\qquad\forall x\in\mathbb{R}^2\\,
\] 
where $d\sigma$ denotes the arc length  element on $\partial\mathcal{O}$. As is well known, $v[\partial\mathcal{O},\phi]$ is a continuous function from $\mathbb{R}^2$ to $\mathbb{R}$. The restriction $v^+[\partial\mathcal{O},\phi]\equiv v[\partial\mathcal{O},\phi]_{|\overline{\mathcal{O}}}$ belongs to $C^{1,\alpha}(\overline{\mathcal{O}})$. Moreover, if  we denote by $C^{1,\alpha}_{\mathrm{loc}}(\mathbb{R}^2\setminus\mathcal{O})$  the space of functions on $\mathbb{R}^2\setminus\mathcal{O}$ whose restrictions to $\overline{\mathcal{U}}$ belong to  $C^{1,\alpha}(\overline{\mathcal{U}})$ for all open bounded subsets $\mathcal{U}$ of $\mathbb{R}^2\setminus\mathcal{U}$, then $v^-[\partial\mathcal{O},\phi]\equiv v[\partial\mathcal{O},\phi]_{|\mathbb{R}^2\setminus\mathcal{O}}$  belongs to $C^{1,\alpha}_{\mathrm{loc}}(\mathbb{R}^2\setminus\mathcal{O})$. 

If $\psi\in C^{1,\alpha}(\partial\mathcal{O})$, then  the double layer potential is denoted by $w[\partial\mathcal{O},\psi]$:
\[
w[\partial\mathcal{O},\psi](x)\equiv-\int_{\partial\mathcal{O}}\psi(y)\;\nu_{\mathcal{O}}(y)\cdot\nabla S(x-y)\,d\sigma_y\qquad\forall x\in\mathbb{R}^2\,,
\] 
where $\nu_\mathcal{O}$ denotes the outer unit normal to $\partial\mathcal{O}$  and the symbol $\cdot$ denotes the scalar product in $\mathbb{R}^2$. Then the restriction $w[\partial\mathcal{O},\psi]_{|\mathcal{O}}$ extends to a function $w^+[\partial\mathcal{O},\psi]$ of $C^{1,\alpha}(\overline{\mathcal{O}})$ and  the restriction $w[\partial\mathcal{O},\psi]_{|\mathbb{R}^2\setminus\overline{\mathcal{O}}}$ extends to a function $w^-[\partial\mathcal{O},\psi]$ of $C^{1,\alpha}_{\mathrm{loc}}(\mathbb{R}^2\setminus\mathcal{O})$.

The single and the double layer potentials will be used to construct solutions to boundary value problems for the Laplace equation. To do so, we need to understand their boundary behavior. Accordingly, to describe the properties of the trace of the double layer potential on $\partial \mathcal{O}$ and of the normal derivative of the single layer potential, we introduce the boundary integral operators  $W_\mathcal{O}$ and $W^*_\mathcal{O}$:
\[
W_\mathcal{O}[\psi](x)\equiv -\int_{\partial\mathcal{O}}\psi(y)\;\nu_{\mathcal{O}}(y)\cdot\nabla S(x-y)\, d\sigma_y\qquad\forall x\in\partial\mathcal{O}\,,
\] for all $\psi\in C^{1,\alpha}(\partial\mathcal{O})$, and
\[
W^*_\mathcal{O}[\phi](x)\equiv \int_{\partial\mathcal{O}}\phi(y)\;\nu_{\mathcal{O}}(x)\cdot\nabla S(x-y)\, d\sigma_y\qquad\forall x\in\partial\mathcal{O}\,,
\] for all $\phi\in C^{0,\alpha}(\partial\mathcal{O})$. Then $W_\mathcal{O}$ is a compact operator  from   $C^{1,\alpha}(\partial\mathcal{O})$  to itself and $W^*_\mathcal{O}$   is a compact operator  from  $C^{0,\alpha}(\partial\mathcal{O})$ to itself (see Schauder \cite{Sc31} and \cite{Sc32}). The operators $W_\mathcal{O}$ and $W^*_\mathcal{O}$ are adjoint one to the other with respect to the duality on $C^{1,\alpha}(\partial\mathcal{O})\times C^{0,\alpha}(\partial\mathcal{O})$ induced by the inner product of the Lebesgue space $L^2(\partial\mathcal{O})$ (cf., {\it e.g.}, Kress \cite[Chap.~4]{Kr14}). For the theory of dual systems and the corresponding Fredholm Alternative Principle, we refer the reader to Kress \cite{Kr14} and Wendland \cite{We67,We70}. Moreover,
\begin{align*}
w^\pm[\partial\mathcal{O},\psi]_{|\partial\mathcal{O}}&=\pm\frac{1}{2}\psi+W_\mathcal{O}[\psi]&\forall\psi\in C^{1,\alpha}(\partial\mathcal{O})\,,\\
\nu_\mathcal{O}\cdot\nabla v^\pm[\partial\mathcal{O},\phi]_{|\partial\mathcal{O}}&=\mp\frac{1}{2}\phi+W^*_\mathcal{O}[\phi]&\forall\phi\in C^{0,\alpha}(\partial\mathcal{O})
\end{align*}
(see, {\it e.g.}, Folland \cite[Chap.~3]{Fo95}).

Finally, we shall need to consider subspaces of $C^{0,\alpha}(\partial \mathcal{O})$ and of $C^{1,\alpha}(\partial \mathcal{O})$, consisting of functions with zero integral on $\partial \mathcal{O}$. Accordingly, we set
\[
C^{k,\alpha}(\partial \mathcal{O})_0\equiv \Bigg\{f \in C^{k,\alpha}(\partial \mathcal{O})\colon \int_{\partial \mathcal{O}}f\, d\sigma=0\Bigg\} \qquad \text{for $k=0,1$}\, .
\]

\subsection{An integral formulation of the boundary value problem}

Our aim is now to convert problem \eqref{eq:direps} into a system of integral equations and we do so by following the strategy of Lanza de Cristoforis \cite{La08} and of \cite{DaMuRo15}. The first attempt to solve \eqref{eq:direps} would be to represent the solution in terms of a double layer potential. However, due to the presence of a hole in the domain, this in general is not possible for all boundary data and we may need to use, for example,  also single layer potentials (cf.~\textit{e.g.}, Folland \cite[Ch.~3]{Fo95}). Thus we need to split the problem in a part which can be solved in terms of the double layer potential and a part which will be represented by a single layer potential. This will be done via Fredholm Theory by characterizing the image of the trace of the double layer potential as the orthogonal to the kernel of the adjoint operator. The dimension of the kernel equals the number of holes in $\Omega_\varepsilon$, and therefore, in this specific case, is equal to one. A real analyticity result upon $\varepsilon$ for the generator of the kernel is provided by Proposition \ref{rhoeps} (see also Remark \ref{rem:rhoeps}).  Now we proceed as in \cite{DaMuRo15} and we introduce the map $M\equiv(M^o,M^i,M^c)$  from $]-\varepsilon_0,\varepsilon_0[\times C^{0,\alpha}(\partial\Omega)\times C^{0,\alpha}(\partial\omega)$ to $C^{0,\alpha}(\partial\Omega)\times C^{0,\alpha}(\partial\omega)_0\times\mathbb{R}$   by setting
\begin{align*}
&M^o[\varepsilon,\rho^o,\rho^i](x) \equiv \frac{1}{2}\rho^o(x)+W^*_{\Omega}[\rho^o](x)+\int_{\partial\omega}\rho^i(s)\;\nu_{\Omega}(x)\cdot\nabla S(x-\varepsilon s)\, d\sigma_s&\forall x\in\partial\Omega\,,\\
&M^i[\varepsilon,\rho^o,\rho^i](t)\equiv\frac{1}{2}\rho^i(t)-W^*_{\omega}[\rho^i](t)-\varepsilon \int_{\partial\Omega}\rho^o(y)\;\nu_{\omega}(t)\cdot\nabla S(\varepsilon t- y)\, d\sigma_y&\forall t\in\partial\omega\,,\\
&M^c[\varepsilon,\rho^o,\rho^i]\equiv \int_{\partial\omega}\rho^i\, d\sigma-1\,,
\end{align*} 
for all $(\varepsilon,\rho^o,\rho^i) \in ]-\varepsilon_0,\varepsilon_0[\times C^{0,\alpha}(\partial\Omega)\times C^{0,\alpha}(\partial\omega)$.
Then we can prove the following result of Lanza de Cristoforis \cite[\S3]{La08} (see also \cite[Prop.~4.1]{DaMuRo15}).

\begin{prop}\label{rhoeps}
The following statements hold.
\begin{itemize}
\item[(i)] The map $M$ is real analytic. 
\item[(ii)] If $\varepsilon\in]-\varepsilon_0,\varepsilon_0[$, then there exists a unique pair $(\rho^o[\varepsilon],\rho^i[\varepsilon])\in C^{0,\alpha}(\partial\Omega)\times C^{0,\alpha}(\partial\omega)$ such that $M[\varepsilon,\rho^o[\varepsilon],\rho^i[\varepsilon]]=0$. \item[(iii)] The map from $]-\varepsilon_0,\varepsilon_0[$ to $C^{0,\alpha}(\partial\Omega)\times C^{0,\alpha}(\partial\omega)$ which takes $\varepsilon$ to   $(\rho^o[\varepsilon],\rho^i[\varepsilon])$ is real analytic.
\end{itemize}
\end{prop}

\begin{rem}\label{rem:rhoeps}
For each $\varepsilon \in ]-\varepsilon_0,\varepsilon_0[\setminus \{0\}$, let $\tau_\varepsilon$ be defined by $\tau_\varepsilon(x)\equiv\rho^o[\varepsilon](x)$ for all $x \in \partial \Omega$ and $\tau_\varepsilon(x)\equiv |\varepsilon|^{-1}\rho^i[\varepsilon](x/\varepsilon)$ for all $x \in \partial (\varepsilon \omega)$. Then
\[
\frac{1}{2}\tau_\varepsilon+W^*_{\Omega_\varepsilon}[\tau_\varepsilon]=0 \, , \qquad \int_{\partial (\varepsilon \omega)}\tau_\varepsilon\, d\sigma=1\, ,
\]
for all $\varepsilon \in ]-\varepsilon_0,\varepsilon_0[\setminus \{0\}$.
\end{rem}

We now turn to consider the part which can be actually solved by the double layer potential. Indeed, by standard Fredholm theory and classical potential theory, one sees that for $\varepsilon \in ]-\varepsilon_0,\varepsilon_0[\setminus \{0\}$ the boundary datum $g_\varepsilon$ defined by
\[
g_\varepsilon(x)\equiv 0 \quad \forall x \in \partial \Omega\, , \qquad g_\varepsilon(x)=U_\#[\varepsilon](x/\varepsilon)-\int_{\partial (\varepsilon \omega)}U_\#[\varepsilon](x/\varepsilon)\tau_\varepsilon(x)\, d\sigma_x \quad \forall x \in \partial (\varepsilon \omega)\, ,
\]
belongs to the image of the trace of the double layer potential (for the definition of $U_\#$ see \eqref{eq:Usharp}). Then, as in \cite{DaMuRo15}, we define the map
$\Lambda\equiv(\Lambda^o,\Lambda^i)$  from $]-\varepsilon_0,\varepsilon_0[\times  C^{1,\alpha}(\partial\Omega)\times C^{1,\alpha}(\partial\omega)_0$ to $C^{1,\alpha}(\partial\Omega)\times C^{1,\alpha}(\partial\omega)$  by
\begin{align*}
&\Lambda^o[\varepsilon,\theta^o,\theta^i](x)\equiv\frac{1}{2}\theta^o(x)+W_{\Omega}[\theta^o](x) \\
&\qquad\qquad\qquad\qquad +\varepsilon\int_{\partial\omega}\theta^i(s)\;\nu_{\omega}(y)\cdot\nabla S(x-\varepsilon s)\, d\sigma_s&\forall x\in\partial\Omega\,,\\
&\Lambda^i[\varepsilon,\theta^o,\theta^i](t)\equiv\frac{1}{2}\theta^i(t)-W_{\omega}[\theta^i](t)+w[\partial\Omega,\theta^o](\varepsilon t)\\
\nonumber
&\qquad\qquad\qquad\qquad -U_\#[\varepsilon](t)+\int_{\partial\omega}U_\#[\varepsilon]\rho^i[\varepsilon]\,d\sigma&\forall t\in\partial\omega\, ,
\end{align*} 
for all $(\varepsilon,\theta^o,\theta^i) \in ]-\varepsilon_0,\varepsilon_0[\times  C^{1,\alpha}(\partial\Omega)\times C^{1,\alpha}(\partial\omega)_0$. Then we have the following result of Lanza de Cristoforis \cite[\S4]{La08} on the regularity of $\Lambda$ (cf. \cite[Prop.~4.3]{DaMuRo15}).

\begin{prop}\label{thetaeps} 
The following statements hold.
\begin{itemize}
\item[(i)] The map $\Lambda$ is real analytic.
\item[(ii)] If $\varepsilon\in ]-\varepsilon_0,\varepsilon_0[$, then there exists a unique pair $(\theta^o[\varepsilon],\theta^i[\varepsilon])\in C^{1,\alpha}(\partial\Omega)\times C^{1,\alpha}(\partial\omega)_0$ such that $\Lambda[\varepsilon,\theta^o[\varepsilon],\theta^i[\varepsilon]]=0$.
\item[(iii)]  The map from $]-\varepsilon_0,\varepsilon_0[$ to $C^{1,\alpha}(\partial\Omega)\times C^{1,\alpha}(\partial\omega)_0$ which takes $\varepsilon$ to $(\theta^o[\varepsilon], \theta^i[\varepsilon])$ is real analytic.
\end{itemize}
\end{prop}

\begin{rem}\label{rem:thetaeps}
For each $\varepsilon \in ]-\varepsilon_0,\varepsilon_0[\setminus \{0\}$, let $\mu_\varepsilon$ be defined by $\mu_\varepsilon(x)\equiv\theta^o[\varepsilon](x)$ for all $x \in \partial \Omega$ and $\mu_\varepsilon(x)\equiv \theta^i[\varepsilon](x/\varepsilon)$ for all $x \in \partial (\varepsilon \omega)$. Then
\[
\frac{1}{2}\mu_\varepsilon+W_{\Omega_\varepsilon}[\mu_\varepsilon]=g_\varepsilon \, ,
\]
for all $\varepsilon \in ]-\varepsilon_0,\varepsilon_0[\setminus \{0\}$.
\end{rem}

By summing the double layer potential with density $\mu_\varepsilon$ (cf.~Remark \ref{rem:thetaeps}) and a convenient multiple of the single layer potential with density $\tau_\varepsilon$ (cf.~Remark \ref{rem:rhoeps}), we can recover the solution $u_\varepsilon$. In particular, by arguing as in \cite[Prop.~4.5]{DaMuRo15}, the following Proposition \ref{solution}   shows how to represent  the rescaled function $u_\varepsilon(\varepsilon t)$ by means of the functions $\rho^o[\varepsilon]$, $\rho^i[\varepsilon]$, $\theta^o[\varepsilon]$, and $\theta^i[\varepsilon]$ introduced in Propositions \ref{rhoeps} and \ref{thetaeps} (see also Lanza de Cristoforis \cite[\S5]{La08} and \cite[\S 2.4]{DaMu13}).

\begin{prop}\label{solution}
Let $\varepsilon\in]-\varepsilon_0,\varepsilon_0[\setminus\{0\}$. Then
\[
\begin{split}
u_\varepsilon(\varepsilon t)& \equiv w^+[\partial\Omega,\theta^o[\varepsilon]](\varepsilon t)-w^-[\partial\omega,\theta^i[\varepsilon]](t)\\
&
+\int_{\partial\omega}U_\#[\varepsilon]\rho^i[\varepsilon]\,d\sigma \biggl(v^+[\partial\Omega,\rho^o[\varepsilon]](\varepsilon t) +v^-[\partial\omega,\rho^i[\varepsilon]](t)+\frac{\log |\varepsilon|}{2\pi}\biggr)\\
&
\times\biggl(\frac{1}{\int_{\partial\omega}d\sigma}\int_{\partial\omega}v[\partial\Omega,\rho^o[\varepsilon]](\varepsilon s)+v[\partial\omega,\rho^i[\varepsilon]](s)\,d\sigma_s+\frac{\log |\varepsilon|}{2\pi}\biggr)^{-1}
\end{split}
\] 
for all $t\in\overline{(\varepsilon^{-1}\Omega)}\setminus\omega$.
\end{prop}

\section{Power series expansions of the auxiliary functions  $(\rho^o[\varepsilon],\rho^i[\varepsilon])$ and $(\theta^o[\varepsilon],\theta^i[\varepsilon])$ around $\varepsilon=0$}\label{sec3}

As described in the Introduction, an intermediate goal is to provide a series expansion in $\varepsilon$ for the integral over $\partial \omega$ of the function
\[
t \mapsto \nu_{\omega}(t)\cdot \nabla \Big(u_\varepsilon(\varepsilon t)\Big) u(\varepsilon t)\, .
\]
Thus, the idea is first to construct an expansion for $\nu_{\omega}(t)\cdot \nabla \Big(u_\varepsilon(\varepsilon t)\Big) u(\varepsilon t)$ and then to integrate such an expansion on $\partial \omega$. Since $u_\varepsilon(\varepsilon t)$ is represented by means of the auxiliary {density} functions  $(\rho^o[\varepsilon],\rho^i[\varepsilon])$ and $(\theta^o[\varepsilon],\theta^i[\varepsilon])$, the plan is to obtain an expansion for those densities and then to get the one for $u_\varepsilon(\varepsilon t)$ by exploiting the representation formula of Proposition \ref{solution}.

In the following Proposition \ref{rhok} of \cite[Prop.~5.1]{DaMuRo15}, we present a power series expansion around $0$ of $(\rho^o[\varepsilon],\rho^i[\varepsilon])$. Throughout the paper, if $j\in\{1,2\}$, then $(\partial_j F)(y)$ denotes the partial derivative with respect to $x_j$ of the function $F(x)\equiv F(x_1,x_2)$ evaluated at $y\equiv(y_1,y_2)\in\mathbb{R}^2$.

\begin{prop}\label{rhok} Let $(\rho^o[\varepsilon],\rho^i[\varepsilon])$ be as in Proposition \ref{rhoeps} for all $\varepsilon\in]-\varepsilon_0,\varepsilon_0[$. Then there exist $\varepsilon_\rho\in]0,\varepsilon_0[$ and a sequence $\{(\rho^o_k,\rho^i_k)\}_{k\in\mathbb{N}}$ in $C^{0,\alpha}(\partial\Omega)\times C^{0,\alpha}(\partial\omega)$ such that 
\[
\rho^o[\varepsilon]=\sum_{k=0}^{+\infty}\frac{\rho^o_k}{k!}\varepsilon^k\quad\text{ and }\quad\rho^i[\varepsilon]=\sum_{k=0}^{+\infty}\frac{\rho^i_k}{k!}\varepsilon^k\qquad \forall \varepsilon\in]-\varepsilon_\rho,\varepsilon_\rho[\,,
\]
where the two series converge uniformly for $\varepsilon\in]-\varepsilon_\rho,\varepsilon_\rho[$ in $C^{0,\alpha}(\partial\Omega)$ and in $C^{0,\alpha}(\partial\omega)$, respectively. Moreover,  the pair of functions $(\rho^o_0,\rho^i_0)$ is the unique solution in $C^{0,\alpha}(\partial\Omega)
\times C^{0,\alpha}(\partial\omega)$ of the following system of integral equations
\begin{align}\nonumber
&\frac{1}{2}\rho^o_0(x)+W^*_{\Omega}[\rho^o_0](x)=-\nu_{\Omega}(x)\cdot\nabla S(x) & \forall x\in\partial\Omega\,,\\
\nonumber
&\frac{1}{2} \rho^i_0(t)-W^*_{\omega}[ \rho^i_0](t)=0 & \forall t\in\partial\omega\,,\\
\nonumber
&\int_{\partial\omega}\rho^i_0\, d\sigma=1\,, 
\end{align}
and for each $k\in\mathbb{N}\setminus\{0\}$ the pair $(\rho^o_k,\rho^i_k)$ is the unique solution in $C^{0,\alpha}(\partial\Omega)
\times C^{0,\alpha}(\partial\omega)$ of the following system of integral equations which involves $\{(\rho^o_j,\rho^i_j)\}_{j=0}^{k-1}$,
\begin{align}\nonumber
&\frac{1}{2}\rho^o_k(x)+W^*_{\Omega}[\rho^o_k](x)\\
\nonumber
&\quad =\sum_{j=0}^{k}\binom{k}{j}(-1)^{j+1}\sum_{h=0}^j\binom{j}{h}\nu_{\Omega}(x)\cdot(\nabla\partial_1^h\partial_2^{j-h} S)(x)\int_{\partial\omega}\rho^i_{k-j}(s)s_1^hs_2^{j-h}\, d\sigma_s& \forall x\in\partial\Omega\,,\\
\nonumber
&\frac{1}{2} \rho^i_k(t)-W^*_{\omega}[ \rho^i_k](t)\\
\nonumber
&\quad =k\sum_{j=0}^{k-1}\binom{k-1}{j}(-1)^{j+1}\sum_{h=0}^j\binom{j}{h}t_1^h t_2^{j-h}\nu_{\omega}(t)\cdot\int_{\partial\Omega}\rho^o_{k-1-j}(\nabla \partial_1^h\partial_2^{j-h} S)\, d\sigma&  \forall t\in\partial\omega\,,\\
\nonumber
&\int_{\partial\omega}\rho^i_k\, d\sigma=0\,.
\end{align}  
\end{prop}

In Proposition \ref{thetak}, instead we determine the coefficients in the power series expansion of $(\theta^o[\varepsilon],\theta^i[\varepsilon])$.

\begin{prop}\label{thetak} Let $(\theta^o[\varepsilon],\theta^i[\varepsilon])$ be as in Proposition \ref{thetaeps} for all $\varepsilon\in]-\varepsilon_0,\varepsilon_0[$. Then there exist $\varepsilon_\theta\in]0,\varepsilon_0[$ and a sequence $\{(\theta^o_k,\theta^i_k)\}_{k\in\mathbb{N}}$ in $C^{1,\alpha}(\partial\Omega)\times C^{1,\alpha}(\partial\omega)_0$ such that 
\begin{equation}\label{thetak0}
\theta^o[\varepsilon]=\sum_{k=0}^\infty\frac{\theta^o_k}{k!}\varepsilon^k\quad\text{and }\quad\theta^i[\varepsilon]=\sum_{k=0}^\infty\frac{\theta^i_k}{k!}\varepsilon^k\qquad\forall \varepsilon\in]-\varepsilon_\theta,\varepsilon_\theta[\,,
\end{equation} where the two series converge uniformly for $\varepsilon\in]-\varepsilon_\theta,\varepsilon_\theta[$ in $C^{1,\alpha}(\partial\Omega)$ and in $C^{1,\alpha}(\partial\omega)_0$, respectively. Moreover,  
\[
(\theta^o_0,\theta^i_0)=(0,0)\, , \qquad \theta^o_1=0\, ,
\]
and $\theta^i_1$ is the unique solution in $C^{1,\alpha}(\partial\omega)_0$ of
\begin{equation}\label{thetai1}
\begin{split}
\frac{1}{2}&\theta^i_1(t)-W_{\omega}[\theta^i_1](t)\\&=\sum_{h=0}^1t_1^h t_2^{1-h} (\partial_1^h\partial_2^{1-h}u)(0)-\sum_{l=0}^1 \sum_{h=0}^l \int_{\partial\omega}s_1^hs_2^{l-h}(\partial_1^h \partial_2^{l-h}u)(0)\rho^i_{1-l}(s)\,d\sigma_s\quad \forall t\in\partial\omega\,{,}
\end{split}
\end{equation}
and for each $k\in\mathbb{N}\setminus\{0,1\}$ the pair $(\theta^o_k,\theta^i_k)$ is the unique solution in $C^{1,\alpha}(\partial\Omega)\times C^{1,\alpha}(\partial\omega)_0$ of the following  system of integral equations which involves $\{(\theta^o_j,\theta^i_j)\}_{j=0}^{k-1}$,
\begin{align}\label{thetaok}
&\frac{1}{2}\theta^o_k(x)+W_{\Omega}[\theta^o_k](x)\\
\nonumber
&=k\sum_{j=0}^{k-2}\binom{k-1}{j}(-1)^{j+1}\sum_{h=0}^j\binom{j}{h}(\nabla \partial_1^h \partial_2^{j-h} S)(x)\cdot \int_{\partial\omega}\theta^i_{k-1-j}(s)\;\nu_{\omega}(s)s_1^h s_2^{j-h} d\sigma_s\\ 
\nonumber
&\qquad\qquad\qquad\qquad\qquad\qquad\qquad\qquad\qquad\qquad\qquad\qquad\qquad\qquad\qquad \forall x\in\partial\Omega\,,\\
\label{thetaik}
&\frac{1}{2}\theta^i_k(t)-W_{\omega}[\theta^i_k](t)=\sum_{j=0}^{k-1}\binom{k}{j}(-1)^{j+1}\sum_{h=0}^j\binom{j}{h}t_1^h t_2^{j-h}\int_{\partial\Omega}\theta^o_{k-j}\nu_{\Omega}\cdot \nabla \partial_1^h \partial_2^{j-h} S \,d\sigma\\ \nonumber
&+\sum_{h=0}^k\binom{k}{h}t_1^h t_2^{k-h} (\partial_1^h\partial_2^{k-h}u)(0)-\sum_{l=0}^k \sum_{h=0}^l \binom{k}{l}\binom{l}{h}\int_{\partial\omega}s_1^hs_2^{l-h}(\partial_1^h \partial_2^{l-h}u)(0)\rho^i_{k-l}(s)\,d\sigma_s\\ \nonumber
&\qquad\qquad\qquad\qquad\qquad\qquad\qquad\qquad\qquad\qquad\qquad\qquad\qquad\qquad\qquad  \forall t\in\partial\omega\,.\nonumber
\end{align}
\end{prop}
\proof  We follow the strategy of \cite[Prop.~5.2]{DaMuRo15}. We first note that the real analyticity of the map which takes $\varepsilon$ to $(\theta^o[\varepsilon],\theta^i[\varepsilon])$  (cf.~Proposition \ref{thetaeps} (iii)) imply the existence of $\varepsilon_\theta$ and $\{(\theta^o_k,\theta^i_k)\}_{k\in\mathbb{N}}$ such that \eqref{thetak0} holds.  Clearly, by Proposition \ref{thetaeps} (ii) we have 
\[
\Lambda[\varepsilon,\theta^o[\varepsilon],\theta^i[\varepsilon]]=0 \qquad \forall \varepsilon\in]-\varepsilon_0,\varepsilon_0[\, .
\] 
By computing the derivative with respect to $\varepsilon$ in the equality above, we deduce that 
\[
\partial_\varepsilon^k(\Lambda[\varepsilon,\theta^o[\varepsilon],\theta^i[\varepsilon]])=0 \qquad \forall \varepsilon\in]-\varepsilon_0,\varepsilon_0[\, , \forall k\in\mathbb{N}\, .
\]
Therefore,
\begin{align}\label{thetaok'}
&\partial_\varepsilon^k(\Lambda^o[\varepsilon,\theta^o[\varepsilon],\theta^i[\varepsilon]])(x)=\frac{1}{2}\partial^k_\varepsilon\theta^o[\varepsilon](x)+W_{\Omega}[\partial^k_\varepsilon\theta^o[\varepsilon]](x)\\
\nonumber
&+\varepsilon\sum_{j=0}^{k}\binom{k}{j}(-1)^j\sum_{h=0}^j\binom{j}{h}\int_{\partial\omega}\partial^{k-j}_\varepsilon\theta^i[\varepsilon](s)\;s_1^h s_2^{j-h}\nu_{\omega}(s) \cdot (\nabla\partial_1^h \partial_2^{j-h} S)(x-\varepsilon s)\, d\sigma_s\\
\nonumber
&+k\sum_{j=0}^{k-1}\binom{k-1}{j}(-1)^j\sum_{h=0}^j\binom{j}{h} \int_{\partial\omega}\partial^{k-1-j}_\varepsilon\theta^i[\varepsilon](s)\;s_1^h s_2^{j-h}\nu_{\omega}(s) \cdot (\nabla\partial_1^h \partial_2^{j-h} S)(x-\varepsilon s)\, d\sigma_s=0\\
\nonumber
 &\qquad\qquad\qquad\qquad\qquad\qquad\qquad\qquad\qquad\qquad\qquad\qquad \forall x\in\partial\Omega\,,
\end{align}
\begin{align}
\label{thetaik'}
&\partial_\varepsilon^k(\Lambda^i[\varepsilon,\theta^o[\varepsilon],\theta^i[\varepsilon]])(t)=\frac{1}{2}\partial^k_\varepsilon\theta^i[\varepsilon](t)-W_{\omega}[\partial^k_\varepsilon\theta^i[\varepsilon]](t)\\
\nonumber
&-\sum_{j=0}^k\binom{k}{j}\sum_{h=0}^j\binom{j}{h}t_1^h t_2^{j-h}\int_{\partial\Omega}\partial_\varepsilon^{k-j}\theta^o[\varepsilon](y)\,\nu_{\Omega}(y)\cdot(\nabla\partial_1^h \partial_2^{j-h} S)(\varepsilon t-y)\,d\sigma_y\\
\nonumber
&-\sum_{h=0}^k\binom{k}{h}t_1^h t_2^{k-h} (\partial_1^h\partial_2^{k-h}u)(\varepsilon t)\\ \nonumber
&+\sum_{l=0}^k \sum_{h=0}^l \binom{k}{l}\binom{l}{h}\int_{\partial\omega}t_1^ht_2^{l-h}(\partial_1^h \partial_2^{l-h}u)(\varepsilon t)\partial_\varepsilon^{k-l}\rho^i[\varepsilon](t)\,d\sigma_t=0\qquad\qquad\qquad\qquad\forall x\in\partial\omega\,, \nonumber
\end{align}
for all $\varepsilon\in]-\varepsilon_0,\varepsilon_0[$ and all  $k\in\mathbb{N}$, where we understand that $\sum_{j=0}^{k-1}$ is omitted for $k=0$.   By classical properties  of real analytic maps, we have $(\theta^o_k,\theta^i_k)=(\partial_\varepsilon^k\theta^o[0],\partial_\varepsilon^k\theta^i[0])$ for all  $k\in\mathbb{N}$. Therefore, by taking $\varepsilon=0$  in \eqref{thetaok'} and \eqref{thetaik'}, we deduce that  $(\theta^o_0,\theta^i_0)=(0,0)$, that $\theta^o_1=0$, that $\theta^i_1$ solves equation \eqref{thetai1},  and that $(\theta^o_k,\theta^i_k)$  is a solution of \eqref{thetaok} and \eqref{thetaik} for all $k\in\mathbb{N}\setminus\{0,1\}$. Then, to conclude, it suffices to note that the uniqueness in $C^{1,\alpha}(\partial\Omega)\times C^{1,\alpha}(\partial\omega)_0$ of the solutions of  \eqref{thetai1} and of \eqref{thetaok}, \eqref{thetaik} follows by classical potential theory (cf., {\it e.g.}, Folland \cite[Chap.~3]{Fo95}). 
 \qed

\section{Series expansion of $\nu_{\omega}(\cdot)\cdot \nabla \big(u_\varepsilon(\varepsilon \cdot)\big) u(\varepsilon \cdot)$ around $\varepsilon=0$}\label{sec4}

We now turn to construct a series expansion for $\nu_{\omega}(\cdot)\cdot \nabla \big(u_\varepsilon(\varepsilon \cdot)\big) u(\varepsilon \cdot)$ for $\varepsilon$ in a neighborhood of $0$, whose coefficients will be defined by means of the sequences $\{(\rho^o_k,\rho^i_k)\}_{k\in\mathbb{N}}$ and $\{(\theta^o_k,\theta^i_k)\}_{k\in\mathbb{N}}$  introduced in Section \ref{sec3}. The strategy is to compute the derivatives with respect to $\varepsilon$ in the representation formula of Proposition \ref{solution} and to exploit the power series expansions for the densities. As a consequence, as in \cite[Prop.~6.1]{DaMuRo15}, the first step is the following Proposition \ref{uk}, where we prove a representation formula which can be easily obtained by   Proposition \ref{solution}, Propositions \ref{rhok} and \ref{thetak}, and by   standard properties of real analytic maps (see also  Lanza de Cristoforis \cite[Theorem 5.3]{La08} and \cite[Theorem 3.1]{DaMu13}).

\begin{prop}\label{uk} Let $\{(\rho^o_k,\rho^i_k)\}_{k\in\mathbb{N}}$ and $\{(\theta^o_k,\theta^i_k)\}_{k\in\mathbb{N}}$ be as in Propositions \ref{rhok} and \ref{thetak}, respectively.  Let 
\[
\begin{split}
&u_{\mathrm{m},0}(t)\equiv0\qquad\qquad\qquad\qquad\qquad\qquad\qquad\qquad  \forall t\in\mathbb{R}^2\setminus\omega\,, \\
& u_{\mathrm{m},1}(t)\equiv- w^-[\partial\omega,\theta^i_1](t) \qquad\qquad\qquad\qquad\qquad  \forall t\in\mathbb{R}^2\setminus\omega\,, \\
&u_{\mathrm{m},k}(t)\equiv\frac{1}{k!}\sum_{j=0}^{k-1}\binom{k}{j}(-1)^{j}\sum_{h=0}^j\binom{j}{h} t_1^h t_2^{j-h} \int_{\partial\Omega}\theta^o_{k-j}\,\nu_{\Omega}\cdot(\nabla\partial_1^h \partial_2^{j-h} S)\,d\sigma\\
&\qquad\qquad - \frac{1}{k!}w^-[\partial\omega,\theta^i_k](t) \quad\qquad\qquad\qquad\qquad  \forall t\in\mathbb{R}^2\setminus\omega\,, \quad \forall k \geq 2 \\
\end{split}
\]
and
 \[
\begin{split}
&v_{\mathrm{m},k}(t)\equiv\frac{1}{k!}\sum_{j=0}^k \binom{k}{j}(-1)^j\sum_{h=0}^j\binom{j}{h} t_1^h t_2^{j-h}\int_{\partial\Omega}\rho^o_{k-j}\partial_1^h \partial_2^{j-h} S\,d\sigma
 +\frac{1}{k!}v^-[\partial\omega,\rho^i_k](t)\\
&\qquad\qquad\qquad\qquad\qquad\qquad\qquad\qquad\qquad\qquad\qquad\qquad\qquad\qquad \forall t\in\mathbb{R}^2\setminus\omega\,, \\
&g_k\equiv\frac{1}{k!}\sum_{l=0}^k \sum_{h=0}^l \binom{k}{l}\binom{l}{h}\int_{\partial\omega}s_1^hs_2^{l-h}(\partial_1^h \partial_2^{l-h}u)(0)\rho^i_{k-l}(s)\,d\sigma_s\,,\\
&r_k\equiv\frac{1}{k!\int_{\partial\omega}d\sigma}\sum_{j=0}^k\binom{k}{j}(-1)^j\sum_{h=0}^j\binom{j}{h} \int_{\partial\omega}s_1^h s_2^{j-h}\,d\sigma_s\int_{\partial\Omega}\rho^o_{k-j}\partial_1^h \partial_2^{j-h}S\,d\sigma\\
&\qquad+\frac{1}{k!\int_{\partial\omega}d\sigma}\int_{\partial\omega}v[\partial\omega,\rho^i_k]\,d\sigma\,,
\end{split}
\] for all $k\in\mathbb{N}$. Then  the following statements hold.
\begin{enumerate}
\item[(i)] There exists $\varepsilon^*\in]0,\varepsilon_0]$ such that the series $\sum_{k=0}^\infty g_{k}\varepsilon^k$ and $\sum_{k=0}^\infty r_{k}\varepsilon^k$ converge absolutely in $]-\varepsilon^*,\varepsilon^*[$. Moreover,
\[
g_0=u(0)\, .
\]
\item[(ii)] If $\Omega_\mathrm{m}\subseteq\mathbb{R}^2\setminus\overline{\omega}$ is open and bounded,   then there exists $\varepsilon_\mathrm{m}\in]0,\varepsilon^*]\cap]0,1[$  such that $\varepsilon\overline{\Omega}_\mathrm{m}\subseteq\Omega$  for all $\varepsilon\in]-\varepsilon_\mathrm{m},\varepsilon_\mathrm{m}[$ and such that
\begin{equation}\label{uepsm}
u_{\varepsilon}(\varepsilon\cdot)_{|\overline{\Omega}_\mathrm{m}}=\sum_{k=1}^\infty u_{\mathrm{m},k|\overline{\Omega}_\mathrm{m}}\varepsilon^k+(\sum_{k=0}^\infty g_{k}\varepsilon^k)\frac{\sum_{k=0}^\infty v_{\mathrm{m},k|\overline{\Omega}_\mathrm{m}}\varepsilon^k+(2\pi)^{-1}\log|\varepsilon|}{\sum_{k=0}^\infty r_{k}\varepsilon^k+(2\pi)^{-1}\log|\varepsilon|}
\end{equation} for all $\varepsilon\in]-\varepsilon_\mathrm{m},\varepsilon_\mathrm{m}[\setminus\{0\}$. Moreover, the series $\sum_{k=1}^\infty u_{\mathrm{m},k|\overline{\Omega}_\mathrm{m}}\varepsilon^k$ and $\sum_{k=0}^\infty v_{\mathrm{m},k|\overline{\Omega}_\mathrm{m}}\varepsilon^k$ converge   in $C^{1,\alpha}(\overline{\Omega}_\mathrm{m})$ uniformly for $\varepsilon\in]-\varepsilon_\mathrm{m},\varepsilon_\mathrm{m}[$.
\end{enumerate}
\end{prop}

By Proposition \ref{uk}, we can then prove an expansion for the map in \eqref{eq:fun}.

\begin{prop}\label{funk} 
With the notation introduced in Proposition \ref{uk}, let 
\[
\begin{split}
&u_{\#,k}(t)\equiv\sum_{\substack{(h,j)\in \mathbb{N}^2\\ h+j=k}}\frac{\partial_1^h \partial_2^{j}u(0)}{h! j!}t^h_1 t_2^j\qquad\qquad\quad  \forall t\in  \mathbb{R}^2\,, \\
& \tilde{u}_{k}(t)\equiv\sum_{l=0}^k \nu_{\omega}(t)\cdot\nabla u_{\mathrm{m},l |\partial \omega}(t) u_{\#,k-l}(t)
\qquad  \forall t\in\partial\omega\,, \\
& \tilde{v}_k(t)\equiv \nu_{\omega}(t)\cdot\nabla v_{\mathrm{m},k|\partial \omega}(t) \qquad  \forall t\in\partial\omega\,, \\
&\tilde{g}_{k}(t)\equiv\sum_{l=0}^k g_l u_{\#,k-l}(t) \quad\qquad\qquad\qquad\qquad  \forall t\in\partial \omega\,, \end{split}
\]
for all $k\in\mathbb{N}$. Then  there exists $\tilde{\varepsilon}\in]0,\varepsilon^{*}]\cap]0,1[$  such that
\begin{equation}\label{funepsm}
\nu_{\omega}(\cdot) \cdot \nabla \big(u_{\varepsilon}(\varepsilon\cdot)\big)_{|\partial \omega}u(\varepsilon \cdot)_{|\partial \omega}=\sum_{k=1}^\infty \tilde{u}_{k}(\cdot)\varepsilon^k+\Bigg(\sum_{k=0}^\infty \tilde{g}_{k}(\cdot)\varepsilon^k\Bigg)\frac{\sum_{k=0}^\infty \tilde{v}_{k}(\cdot)\varepsilon^k}{\sum_{k=0}^\infty r_{k}\varepsilon^k+(2\pi)^{-1}\log|\varepsilon|}
\end{equation} for all $\varepsilon\in]-\tilde{\varepsilon},\tilde{\varepsilon}[\setminus\{0\}$. Moreover, the series $\sum_{k=0}^\infty \tilde{g}_{k}\varepsilon^k$, $\sum_{k=0}^\infty \tilde{u}_{k}\varepsilon^k$, and $\sum_{k=0}^\infty \tilde{v}_{k}\varepsilon^k$ converge   in $C^{0,\alpha}(\partial \omega)$ uniformly for $\varepsilon\in]-\tilde{\varepsilon},\tilde{\varepsilon}[$. 
\end{prop}
\proof
We first note that if we take $\tilde{\varepsilon}\in ]0,\varepsilon^{\ast}[$ small enough, then for $\varepsilon \in ]-\tilde{\varepsilon},\tilde{\varepsilon}[$ we have that
\[
\begin{split}
u(\varepsilon t)&=\sum_{(i,j)\in \mathbb{N}^2}\varepsilon^{i+j}\frac{\partial_1^i \partial_2^{j}u(0)}{i! j!}t^i_1 t_2^j\\
&=\sum_{h=0}^\infty \Bigg(\sum_{\substack{(i,j)\in \mathbb{N}^2\\ i+j=h}}\frac{\partial_1^i \partial_2^{j}u(0)}{i! j!}t^i_1 t_2^j\Bigg) \varepsilon^{h}=\sum_{h=0}^\infty u_{\#,h}(t) \varepsilon^{h} \qquad \forall t \in \partial \omega\, ,
\end{split}
\]
and that the power series $\sum_{h=0}^\infty u_{\#,h|\partial \omega} \varepsilon^{h} $ converges in $C^{0,\alpha}(\partial \omega)$ uniformly for $\varepsilon \in ]-\tilde{\varepsilon},\tilde{\varepsilon}[$.
Possibly taking a smaller $\tilde{\varepsilon}$,  we observe that for $\varepsilon \in ]-\tilde{\varepsilon},\tilde{\varepsilon}[$ we have
\[
\Big(\sum_{k=1}^\infty \nu_{\omega} \cdot \nabla u_{\mathrm{m},k|\partial \omega}\varepsilon^k\Big)\Big(\sum_{h=0}^\infty u_{\#,h|\partial \omega} \varepsilon^{h}\Big)=\sum_{k=0}^\infty \tilde{u}_{k}\varepsilon^k\, ,  \qquad \Big(\sum_{k=0}^\infty g_{k}\varepsilon^k\Big)\Big(\sum_{h=0}^\infty u_{\#,h|\partial \omega} \varepsilon^{h}\Big)=\sum_{k=0}^\infty \tilde{g}_{k}\varepsilon^k
\]
where the series converge in $C^{0,\alpha}(\partial \omega)$ uniformly for $\varepsilon \in ]-\tilde{\varepsilon},\tilde{\varepsilon}[$ and we have set
\[
\tilde{u}_{k}\equiv\sum_{l=0}^k \nu_{\omega} \cdot \nabla u_{\mathrm{m},l |\partial \omega} u_{\#,k-l|\partial \omega}\, , \qquad \tilde{g}_{k}\equiv\sum_{l=0}^k g_l u_{\#,k-l|\partial \omega}\, . 
\]
Then the validity of \eqref{funepsm} follows by Proposition \ref{uk} (see formula \eqref{uepsm}). \qed

Now we would like to obtain an expression for $\nu_{\omega}(\cdot) \cdot \nabla \big(u_{\varepsilon}(\varepsilon\cdot)\big)_{|\partial \omega}u(\varepsilon \cdot)_{|\partial \omega}$ in the form of a convergent series of the type
\[
\sum_{n=0}^\infty  \varphi_\varepsilon(\cdot) \varepsilon^n \, .
\]
On the other hand, because  of the quotient in \eqref{funepsm}, we don't have yet an expression as above. However, by exploiting exactly the same argument of \cite[Thm.~6.3]{DaMuRo15}, we can prove  Theorem  \ref{umk} below where we exhibit  a series expansion for the map which takes $\varepsilon$ to  $\nu_{\omega}(\cdot) \cdot \nabla \big(u_{\varepsilon}(\varepsilon\cdot)\big)_{|\partial \omega}u(\varepsilon \cdot)_{|\partial \omega}$.

\begin{theorem}\label{umk}
With the notation introduced in Proposition \ref{uk}, let $\{\tilde{a}_{n}\}_{n\in\mathbb{N}}$ be the sequence of functions from $\partial \omega$ to $\mathbb{R}$ defined by 
\[
\tilde{a}_{n}\equiv\sum_{k=0}^n \tilde{g}_{n-k}\tilde{v}_{k}\qquad\forall n\in\mathbb{N}\,.
\] Let $\{\tilde{\lambda}_{(n,l)}\}_{(n,l)\in\mathbb{N}^2\,,\;l\le n+1}$ be the family of functions from $\partial \omega$ to $\mathbb{R}$ defined by 
\[
\tilde{\lambda}_{(n,0)}\equiv \tilde{u}_{n}\,,\quad\tilde{\lambda}_{(n,1)}\equiv \tilde{a}_{n}\,,
\] for all $n\in\mathbb{N}$, and
\[
\tilde{\lambda}_{(n,l)}\equiv (-1)^{l-1} \sum_{k=l-1}^n \tilde{a}_{n-k} \sum_{\beta\in(\mathbb{N}\setminus\{0\})^{l-1}\,,\;|\beta|=k}\ \prod_{h=1}^{l-1} r_{\beta_h}
\] for all $n,l\in\mathbb{N}$ with $2\le l\le n+1$. 
  Then there exists $\tilde{\varepsilon}'\in]0,\varepsilon_0]\cap]0,1[$ such that 
  \begin{equation}\label{fuepsmseries}
\nu_{\omega}(\cdot) \cdot \nabla \big(u_{\varepsilon}(\varepsilon\cdot)\big)_{|\partial \omega}u(\varepsilon \cdot)_{|\partial \omega}=\sum_{n=0}^\infty\varepsilon^n\sum_{l=0}^{n+1} \frac{\tilde{\lambda}_{(n,l)}(\cdot)}{(r_0+(2\pi)^{-1}\log|\varepsilon|)^{l}}
\end{equation} for all $\varepsilon\in]-\tilde{\varepsilon}',\tilde{\varepsilon}'[\setminus\{0\}$. Moreover, the series 
\[
\sum_{n=0}^\infty\varepsilon^n\sum_{l=0}^{n+1} \frac{\tilde{\lambda}_{(n,l)}\eta^l}{(r_0\eta+(2\pi)^{-1})^{l}}\] converges in $C^{1,\alpha}(\partial \omega)$ uniformly for $(\varepsilon,\eta)\in]-\tilde{\varepsilon}',\tilde{\varepsilon}'[\times]1/\log\tilde{\varepsilon}',-1/\log\tilde{\varepsilon}'[$.
 \end{theorem}

\medskip

\begin{rem}\label{rem:uepsm}
With the notation of Theorem \ref{umk}, a straightforward computation shows that
\[
\begin{split}
 \tilde{\lambda}_{(0,0)}=&\tilde{u}_{0}=0\, ,  \\
  \tilde{\lambda}_{(0,1)}=&\tilde{a}_{0}=\big(u(0)\big)^2 \frac{\partial}{\partial \nu_{\omega}}v^-[\partial \omega, \rho^i_0]\, .
\end{split}
\]
\end{rem}

\medskip

\section{Series expansion of $\mathrm{Cap}_{\Omega}(\varepsilon \omega,u)$}\label{sec5}

Our aim is now to deduce a full expansion for the $u$-capacity $\mathrm{Cap}_\Omega(\varepsilon\omega,u)$, which is given  as the sum of $\int_{\Omega_\varepsilon}|\nabla u_\varepsilon |^2 \, dx$ and of $\int_{\varepsilon \omega}|\nabla u|^2 \, dx$. As a first step, we provide an expansion for $\int_{\varepsilon \omega}|\nabla u|^2 \, dx$ around $\varepsilon=0$. As we shall see, the term $\int_{\varepsilon \omega}|\nabla u|^2 \, dx$ depends analytically on $\varepsilon$ and thus can be expanded in a power series. Therefore, we compute such a power series in the following lemma.

\begin{lem}\label{lem:nrguom}
Let $\{\xi_{n}\}_{n\in\mathbb{N}}$ be the sequence of real numbers defined by 
\[
\begin{split}
& \xi_0\equiv0\, , \qquad \xi_1\equiv0\, , \qquad\xi_{n}\equiv\sum_{j=1}^2\sum_{l=0}^{n-2} \int_{\omega} \partial_j u_{\#,l+1}(t)\partial_j u_{\#,n-l-1}(t)\, dt\qquad \forall n \geq 2\ .
\end{split}
\] 
Then there exists $\varepsilon_\xi\in]0,\varepsilon_0]$ such that 
\[
\int_{\varepsilon \omega}|\nabla u|^2 \, dx=\sum_{n=2}^\infty\xi_n \varepsilon^n
\]
for all $\varepsilon\in]-\varepsilon_\xi,\varepsilon_\xi[\setminus\{0\}$. Moreover, 
\[
\xi_2 = |\nabla u(0)|^2 m_2(\omega) \, ,
\]
and the series 
\[
\sum_{n=2}^\infty\xi_n \varepsilon^n
\] 
converges uniformly for $\varepsilon \in]-\varepsilon_\xi,\varepsilon_\xi[$. (The symbol $m_2(\dots)$ denotes the two-dimensional Lebesgue measure of a set).
\end{lem}
\proof 
If $\varepsilon \in ]-\varepsilon_0,\varepsilon_0[\setminus\{0\}$, by the Theorem of change of variable in integrals, we have
\[
\int_{\varepsilon \omega}|\nabla u|^2 \, dx=\varepsilon^2\int_{ \omega}|(\nabla u)(\varepsilon t)|^2 \, dt\, .
\]
Then we note that by assumption \eqref{eq:assf}  on the analyticity of $u$, by analyticity results for the composition operator (cf.~B\"{o}hme and Tomi~\cite[p.~10]{BoTo73}, 
Henry~\cite[p.~29]{He82}, Valent~\cite[Thm.~5.2, p.~44]{Va88}), there exists $\varepsilon_\xi \in ]0,\varepsilon_0]$ such that the map from $]-\varepsilon_\xi,\varepsilon_\xi[$ to $C^{0,\alpha}(\overline{\omega})$ which takes $\varepsilon$ to $(\partial_j u)(\varepsilon \cdot)_{|\overline{\omega}}$ is real analytic. Possibly shrinking $\varepsilon_\xi$, one verifies that  for $\varepsilon \in ]-\varepsilon_\xi,\varepsilon_\xi[\setminus \{0\}$, 
\[
\begin{split}
(\partial_j u )(\varepsilon t) &= \frac{1}{\varepsilon} \partial_j(u(\varepsilon t))\\
&=\frac{1}{\varepsilon}\sum_{h=0}^\infty \partial_j u_{\#,h}(t) \varepsilon^{h}\\
&=\frac{1}{\varepsilon}\sum_{h=1}^\infty \partial_j u_{\#,h}(t) \varepsilon^{h-1}\varepsilon=\sum_{h=0}^\infty \partial_j u_{\#,h+1}(t) \varepsilon^{h} \qquad \forall t \in \overline{\omega}\, ,
\end{split}
\]
where the series $\sum_{h=0}^\infty \partial_j u_{\#, h+1|\overline{\omega}} \varepsilon^h$ converges in $C^{0,\alpha}(\overline{\omega})$ uniformly for $\varepsilon \in ]-\varepsilon_\xi,\varepsilon_\xi[$. As a consequence,
\[
(\partial_j u )^2(\varepsilon t) =  \sum_{n=0}^\infty \Bigg(\sum_{l=0}^n \partial_j u_{\#,l+1}(t)\partial_j u_{\#,n-l+1}(t)\Bigg) \varepsilon^n \qquad \forall t \in \overline{\omega}\, , \forall \varepsilon \in ]-\varepsilon_\xi,\varepsilon_\xi[\setminus \{0\}\, .
\]
By the continuity of the linear operator from $C^{0,\alpha}(\overline{\omega})$ to $\mathbb{R}$ which takes a function $h$ to its integral $\int_{\omega}h\, dt$, by summing on $j \in \{1,2\}$, one deduces that {possibly taking a smaller $\varepsilon_\xi$}
\begin{equation}\label{eq:nrguom1}
\int_{ \omega}|(\nabla u)(\varepsilon t)|^2 \, dt=\sum_{n=0}^\infty \bigg(\sum_{j=1}^2\sum_{l=0}^n \int_{\omega}\partial_j u_{\#,l+1}(t)\partial_j u_{\#,n-l+1}(t)\, dt\bigg) \varepsilon^n\, ,
\end{equation}
for all $\varepsilon \in ]-\varepsilon_\xi,\varepsilon_\xi[\setminus \{0\}$. In particular,
\[
\begin{split}
\sum_{j=1}^2 \int_{\omega}\partial_j u_{\#,1}(t)\partial_j u_{\#,1}(t)\, dt&= \int_{\omega}\Big( (\partial_1 u(0))^2+(\partial_2 u(0))^2\Big)\, dt\\
&= |\nabla u(0)|^2\int_{\omega}\, dt=|\nabla u(0)|^2m_2(\omega)\, .
\end{split}
\]
Then, by multiplying equation \eqref{eq:nrguom1} by $\varepsilon^2$, we deduce the validity of the lemma.\qed

\medskip

By integrating over $\partial \omega$  formula \eqref{fuepsmseries} and adding the coefficients of Lemma \ref{lem:nrguom}, by Theorem \ref{umk} we can immediately deduce the validity of our main result on the asymptotic behavior of $\mathrm{Cap}_{\Omega}(\varepsilon \omega,u)$.

\begin{theorem}\label{capk}
With the notation introduced in Proposition \ref{uk}, Theorem \ref{umk} and Lemma \ref{lem:nrguom}, let $\{c_{(n,l)}\}_{\substack{(n,l)\in\mathbb{N}^2\\\;l\le n+1}}$ be the family of real numbers defined by
\[
c_{(n,l)}\equiv-\int_{\partial \omega}\tilde{\lambda}_{(n,l)}\, d\sigma+\delta_{0,l}\xi_n\, ,
\] 
for all $n,l\in\mathbb{N}$ with $l\le n+1$ (where $\delta_{0,l}=1$ if $l=0$ and $\delta_{0,l}=0$ if $l\neq0$).   Then there exists $\varepsilon_\mathrm{c}\in]0,\varepsilon_0]\cap]0,1[$ such that 
\[
\mathrm{Cap}_{\Omega}(\varepsilon \omega,u)=\sum_{n=0}^\infty\varepsilon^n\sum_{l=0}^{n+1} \frac{c_{(n,l)}}{(r_0+(2\pi)^{-1}\log|\varepsilon|)^{l}}
\]
for all $\varepsilon\in]-\varepsilon_{\mathrm{c}},\varepsilon_\mathrm{c}[\setminus\{0\}$. Moreover, the series 
\[
\sum_{n=0}^\infty\varepsilon^n\sum_{l=0}^{n+1} \frac{c_{(n,l)}\eta^l}{(r_0\eta+(2\pi)^{-1})^{l}}\] converges  uniformly for $(\varepsilon,\eta)\in]-\varepsilon_\mathrm{c},\varepsilon_\mathrm{c}[\times]1/\log\varepsilon_\mathrm{c},-1/\log\varepsilon_\mathrm{c}[$.
 \end{theorem}
 
 \medskip
 
 \begin{rem}\label{rem:cepsm}
With the notation of Theorem \ref{capk}, we observe that Remark \ref{rem:uepsm} and a straightforward computation   based on Folland \cite[Lem.~3.30]{Fo95} imply that
\[
\begin{split}
 c_{(0,0)}=&0\, ,  \\
c_{(0,1)}=& -\int_{\partial \omega} \big(u(0)\big)^2 \frac{\partial}{\partial \nu_{\omega}}v^-[\partial \omega, \rho^i_0]\, d\sigma=-\big(u(0)\big)^2\int_{\partial \omega}\rho^i_0\, d\sigma=-\big(u(0)\big)^2\, .
\end{split}
\]
Moreover, if we denote by  $H^o_0$   the unique solution in $C^{1,\alpha}(\overline{\Omega})$ of
\[
\left\{
\begin{array}{ll}
\Delta H^o_0=0&\text{in }\Omega\,,\\
H^o_{0}(x)=S(x)&\text{for all }x\in\partial\Omega\,,
\end{array}
\right.
\] 
and by $H^i_0$  the unique solution  in $C^{1,\alpha}_{\mathrm{loc}}(\mathbb{R}^2\setminus\omega)$ of
\[
\left\{
\begin{array}{ll}
\Delta H^i_0=0&\text{in }\mathbb{R}^2\setminus\overline{\omega}\,,\\
H^i_{0}(t)=S(t)&\text{for all }t\in\partial\omega\,,\\
\sup_{t\in\mathbb{R}^2\setminus\omega}|H^i_0(t)|<+\infty\,,
\end{array} 
\right.
\] 
then by \cite[Prop.~7.3]{DaMuRo15} we have
\[
r_0= \lim_{t\to\infty}H^i_0(t)-H^o_0(0)\,.
\]
Accordingly,
\begin{equation}\label{eq:capepsfirst}
\begin{split}
\mathrm{Cap}_{\Omega}(\varepsilon \omega,u)=& -\frac{\big(u(0)\big)^2}{\lim_{t\to\infty}H^i_0(t)-H^o_0(0)+(2\pi)^{-1}\log|\varepsilon|}\\
&+\varepsilon\bigg(\sum_{n=1}^\infty\varepsilon^{n-1}\sum_{l=0}^{n+1} \frac{c_{(n,l)}}{(\lim_{t\to\infty}H^i_0(t)-H^o_0(0)+(2\pi)^{-1}\log|\varepsilon|)^{l}}\bigg)
\end{split}
\end{equation}
for all $\varepsilon\in]-\varepsilon_{\mathrm{c}},\varepsilon_\mathrm{c}[\setminus\{0\}$. Moreover, in case $\omega$ is a Jordan domain, we deduce by \cite[\S 4]{MuMi15} that $e^{2\pi \lim_{t\to\infty}H^i_0(t)}$ is the logarithmic capacity (or  outer conformal radius)  of $\omega$. $H^o_0(0)$ is  the value at $0$ of the unique harmonic function in $\Omega$ which agrees with $S$ on $\partial \Omega$. In other words, 
\[
H_{(0,0)}=-H^o_0(0) \, , \qquad N=\lim_{t\to\infty}H^i_0(t)\, ,
\]
where $H_{(0,0)}$ and $N$ are as in formula \eqref{eq:MNPexp}. Finally, we note that the if we look at the first summand in the right hand side of equality \eqref{eq:capepsfirst}, then the information on the function $u$ is in the numerator, whereas the geometry of $\Omega$ and $\omega$ is taken into account in the denominator. 
\end{rem}

\subsection{Asymptotic behavior of $\mathrm{Cap}_{\Omega}(\varepsilon \omega,u)$ under vanishing assumption for $u$}

We now assume that there exists $\overline{k} \in \mathbb{N}\setminus \{0\}$ such that
\begin{equation}\label{eq:vanu}
D^\gamma u(0)=0 \quad \forall |\gamma| <\overline{k}\, , \qquad D^\beta u(0)\neq 0 \quad \mbox{for some $\beta \in \mathbb{N}^2$ with $|\beta|=\overline{k}$}\, .
\end{equation}
Then condition \eqref{eq:vanu} and Proposition \ref{thetak} imply that
\begin{equation}\label{eq:vantheta}
(\theta^o_k,\theta^i_k)=(0,0) \quad \forall k < \overline{k}\, , \qquad \theta^o_{\overline{k}}=0\, ,
\end{equation}
and that $\theta^i_{\overline{k}}$ is the unique solution in $C^{1,\alpha}(\partial\omega)_0$ of  
\[
\begin{split}
\frac{1}{2}\theta^i_{\overline{k}}(t)-W_{\omega}[\theta^i_{\overline{k}}](t)&=\sum_{h=0}^{\overline{k}}\binom{{\overline{k}}}{h}t_1^h t_2^{{\overline{k}}-h} (\partial_1^h\partial_2^{\overline{k}-h}u)(0)\\
&-\sum_{h=0}^{\overline{k}}  \binom{{\overline{k}}}{h}\int_{\partial\omega}s_1^hs_2^{{\overline{k}}-h}(\partial_1^h \partial_2^{{\overline{k}}-h}u)(0)\rho^i_{0}(s)\,d\sigma_s \qquad \forall t\in\partial\omega\,,
\end{split}
\]
\textit{i.e.}, 
\begin{equation}\label{eq:thetaooverk}
\begin{split}
\frac{1}{2}\theta^i_{\overline{k}}(t)-W_{\omega}[\theta^i_{\overline{k}}](t)&=\overline{k}! \Bigg (u_{\#,\overline{k}}(t)-\int_{\partial\omega}u_{\#,\overline{k}}\rho^i_{0}\,d\sigma\Bigg) \qquad \forall t\in\partial\omega\,.
\end{split}
\end{equation}
Then equations \eqref{eq:vantheta}, \eqref{eq:thetaooverk}, and Proposition \ref{uk} imply that
\begin{equation}\label{eq:vanum}
u_{\mathrm{m},k}=0 \qquad \forall k <\overline{k}\, ,\qquad u_{\mathrm{m},\overline{k}}=-\frac{1}{\overline{k}!}w^-[\partial \omega, \theta^i_{\overline{k}}]\, .
\end{equation}
As a consequence, by classical potential theory, $u_{\mathrm{m},\overline{k}}$ is the unique solution in $C^{1,\alpha}_{\mathrm{loc}}(\mathbb{R}^2 \setminus \omega)$ of the following problem
\begin{equation}\label{eq:bvpumk}
\left\{
\begin{array}{ll}
\Delta u_{\mathrm{m},\overline{k}}=0&\text{in }\mathbb{R}^2\setminus\overline{\omega}\,,\\
u_{\mathrm{m},\overline{k}}(t)=u_{\#,\overline{k}}(t)-\int_{\partial\omega}u_{\#,\overline{k}}\rho^i_{0}\,d\sigma&\text{for all }t\in\partial\omega\,,\\
\sup_{t\in\mathbb{R}^2\setminus\omega}|u_{\mathrm{m},\overline{k}}(t)|<+\infty\,.
\end{array} 
\right.
\end{equation}
Moreover, by assumption \eqref{eq:vanu} and Proposition \ref{uk} we have
\begin{equation}\label{eq:vang}
g_{k}=0 \quad \forall k <\overline{k}\, ,\qquad g_{\overline{k}}=\frac{1}{\overline{k}!}\sum_{h=0}^{\overline{k}}  \binom{{\overline{k}}}{h}\int_{\partial\omega}s_1^hs_2^{{\overline{k}}-h}(\partial_1^h \partial_2^{{\overline{k}}-h}u)(0)\rho^i_{0}(s)\,d\sigma_s=\int_{\partial\omega}u_{\#,\overline{k}}\rho^i_{0}\,d\sigma\, .
\end{equation}
Then by \eqref{eq:vanu} and by Propostion \ref{funk} we verify that
\begin{equation}\label{eq:vanusharp}
u_{\#,k}=0 \qquad \forall k <\overline{k}\, ,
\end{equation}
and accordingly Proposition \ref{funk} and equations \eqref{eq:vanum}, \eqref{eq:vanusharp} imply
\begin{equation}\label{eq:vantildeu}
\tilde{u}_{k}=0 \qquad \forall k <2\overline{k}\, ,\qquad \tilde{u}_{2\overline{k}}=u_{\#,\overline{k}|\partial \omega}\Bigg( \frac{\partial u_{\mathrm{m},\overline{k}}}{\partial \nu_\omega}\Bigg)\, .
\end{equation}
Furthermore, by \eqref{eq:vang} and \eqref{eq:vanusharp} we have
\begin{equation}\label{eq:vantildeg}
\tilde{g}_{k}=0 \qquad \forall k <2\overline{k}\, ,\qquad \tilde{g}_{2\overline{k}}=g_{\overline{k}}u_{\#,\overline{k}|\partial \omega}=u_{\#,\overline{k}|\partial \omega}
 \int_{\partial\omega}u_{\#,\overline{k}}\rho^i_{0}\,d\sigma \, .
\end{equation}
Then, as an intermediate step for computing the coefficients of the expansion of the $u$-capacity $\mathrm{Cap}_{\Omega}(\varepsilon \omega,u)$, we consider the quantities $\tilde{a}_{n}, \tilde{\lambda}_{(n,l)}$ introduced in Theorem \ref{umk} for representing the behavior of $\nu_{\omega}(\cdot) \cdot \nabla \big(u_{\varepsilon}(\varepsilon\cdot)\big)_{|\partial \omega}u(\varepsilon \cdot)$. A straightforward computation based on \eqref{eq:vantildeu}, \eqref{eq:vantildeg} implies that
\[
\tilde{a}_{n}=0 \qquad \forall n <2\overline{k}\, ,\qquad \tilde{a}_{2\overline{k}}=\tilde{g}_{2\overline{k}}\tilde{v}_0=  \tilde{v}_0
u_{\#,\overline{k}|\partial \omega}\int_{\partial\omega}u_{\#,\overline{k}}\rho^i_{0}\,d\sigma\, ,
 \]
and accordingly
\begin{equation}\label{eq:vantilde0}
\tilde{\lambda}_{(n,0)}=0 \qquad \forall n <2\overline{k}\, ,\qquad \tilde{\lambda}_{2\overline{k},0}=\tilde{u}_{2\overline{k}}=u_{\#,\overline{k}|\partial \omega} \Bigg( \frac{\partial u_{\mathrm{m},\overline{k}}}{\partial \nu_\omega}\Bigg)\, ,
\end{equation}
\begin{equation}\label{eq:vantilde1}
\tilde{\lambda}_{(n,1)}=0 \qquad \forall n <2\overline{k}\, ,\qquad \tilde{\lambda}_{2\overline{k},1}=\tilde{a}_{2\overline{k}}=  \tilde{v}_0
 u_{\#,\overline{k}|\partial \omega} \int_{\partial\omega}u_{\#,\overline{k}}\rho^i_{0}\,d\sigma\, ,
\end{equation}
and
\begin{equation}\label{eq:vantildel}
\tilde{\lambda}_{(n,l)}=0 \qquad \text{$\forall (n,l)$ such that  $n-l+1<2\overline{k}$ and that $2\leq l \leq n+1$}\, .
\end{equation}
In particular, $\tilde{\lambda}_{(n,l)}=0$ for all $(n,l)$ such that  $n<2\overline{k}+1$ and that $2\leq l \leq n+1$. Moreover, a simple computation shows that
\[
\xi_{n}=0 \qquad \forall n < 2 \overline{k}\, , \qquad \xi_{2 \overline{k}}= \int_{\omega} |\nabla u_{\#,\overline{k}}|^2\, dt\, .
\]
Finally, by Theorem \ref{capk} and by integrating equalities \eqref{eq:vantilde0}-\eqref{eq:vantildel}, we obtain
\[
\begin{split}
&c_{(n,0)}=0 \qquad \forall n <2\overline{k}\, ,\\ 
&c_{2\overline{k},0}=-\int_{\partial \omega}\tilde{u}_{2\overline{k}}\, d\sigma+\int_{\omega} |\nabla u_{\#,\overline{k}}|^2\, dt\\
&\qquad=-\int_{\partial \omega} u_{\#,\overline{k}|\partial \omega} \Bigg( \frac{\partial u_{\mathrm{m},\overline{k}}}{\partial \nu_\omega}\Bigg)\, d\sigma+\int_{\omega} |\nabla u_{\#,\overline{k}}|^2\, dt\,, 
\end{split}
\]
\[
c_{(n,1)}=0 \qquad \forall n <2\overline{k}\, ,\qquad c_{2\overline{k},1}=-\int_{\partial \omega}\tilde{a}_{2\overline{k}}\,d\sigma=-\int_{\partial \omega} \tilde{v}_0
 u_{\#,\overline{k}|\partial \omega}\, d\sigma \int_{\partial\omega}u_{\#,\overline{k}}\rho^i_{0}\,d\sigma\, , 
\]
and
\[
c_{(n,l)}=0 \qquad \text{$\forall (n,l)$ such that  $n-l+1<2\overline{k}$ and that $2\leq l \leq n+1$}\, .
\]
In particular, $c_{(n,l)}=0$ for all $(n,l)$ such that  $n<2\overline{k}+1$ and that $2\leq l \leq n+1$.  Since {$u_{\mathrm{m},\overline{k}}=-\frac{1}{\overline{k}!}w^-[\partial \omega, \theta^i_{\overline{k}}]$, then $u_{\mathrm{m},\overline{k}}$ is harmonic at infinity (cf. \eqref{eq:vanum}). As a consequence,} the decay properties of its radial derivative (cf. Folland \cite[Prop.~2.75]{Fo95}) and the Divergence Theorem imply that
\[
\int_{\partial \omega}\frac{\partial u_{\mathrm{m},\overline{k}}}{\partial \nu_\omega}\, d\sigma=0\, .
\]
Accordingly, 
\[
-\int_{\partial \omega} u_{\#,\overline{k}|\partial \omega} \Bigg( \frac{\partial u_{\mathrm{m},\overline{k}}}{\partial \nu_\omega}\Bigg)\, d\sigma=-\int_{\partial \omega}\Big( u_{\#,\overline{k}|\partial \omega}-\int_{\partial \omega} u_{\#,\overline{k}|\partial \omega}\, d\sigma \Big) \Bigg( \frac{\partial u_{\mathrm{m},\overline{k}}}{\partial \nu_\omega}\Bigg)\, d\sigma\, .
\]
{Since, $u_{\mathrm{m},\overline{k}}$ solves problem \eqref{eq:bvpumk},  we have $u_{\mathrm{m},\overline{k}}=u_{\#,\overline{k}}-\int_{\partial \omega} u_{\#,\overline{k}}\,d\sigma$ on $\partial \omega$, and thus}
\[
-\int_{\partial \omega}\Big( u_{\#,\overline{k}|\partial \omega}-\int_{\partial \omega} u_{\#,\overline{k}|\partial \omega}\, d\sigma \Big) \Bigg( \frac{\partial u_{\mathrm{m},\overline{k}}}{\partial \nu_\omega}\Bigg)\, d\sigma=-\int_{\partial \omega}u_{\mathrm{m},\overline{k}} \Bigg( \frac{\partial u_{\mathrm{m},\overline{k}}}{\partial \nu_\omega}\Bigg)\, d\sigma\, .
\] 
On the other hand,  the harmonicity at infinity of $u_{\mathrm{m},\overline{k}}$ and the Divergence Theorem imply that
\[
0<\int_{\mathbb{R}^2 \setminus \overline{\omega}}|\nabla u_{\mathrm{m},\overline{k}}|^2\, dt=-\int_{\partial \omega}u_{\mathrm{m},\overline{k}} \Bigg( \frac{\partial u_{\mathrm{m},\overline{k}}}{\partial \nu_\omega}\Bigg)\, d\sigma\, 
\]
(cf.~Folland \cite[p.~118]{Fo95}). As a consequence,
\begin{equation}\label{eq:ibp}
-\int_{\partial \omega} u_{\#,\overline{k}|\partial \omega} \Bigg( \frac{\partial u_{\mathrm{m},\overline{k}}}{\partial \nu_\omega}\Bigg)\, d\sigma=\int_{\mathbb{R}^2 \setminus \overline{\omega}}|\nabla u_{\mathrm{m},\overline{k}}|^2\, dt>0\, .
\end{equation}
Moreover, if we denote by $\mathsf{u}_{\overline{k}}$ the unique solution in $C^{1,\alpha}_{\mathrm{loc}}(\mathbb{R}^2\setminus \omega)$ of
\begin{equation}\label{eq:bvp:mathsfu}
\left\{
\begin{array}{ll}
\Delta \mathsf{u}_{\overline{k}}=0&\text{in }\mathbb{R}^2\setminus\overline{\omega}\,,\\
\mathsf{u}_{\overline{k}}(t)=u_{\#,\overline{k}}(t)&\text{for all }t\in\partial\omega\,,\\
\sup_{t\in\mathbb{R}^2\setminus\omega}|\mathsf{u}_{\overline{k}}(t)|<+\infty\,,
\end{array} 
\right.
\end{equation}
then clearly
\[
\mathsf{u}_{\overline{k}}=u_{\mathrm{m},\overline{k}}+\int_{\partial\omega}u_{\#,\overline{k}}\rho^i_{0}\,d\sigma\, ,
\]
and thus
\[
\int_{\mathbb{R}^2 \setminus \overline{\omega}}|\nabla \mathsf{u}_{\overline{k}}|^2\, dt=\int_{\mathbb{R}^2 \setminus \overline{\omega}}|\nabla u_{\mathrm{m},\overline{k}}|^2\, dt\, .
\]
We now turn to consider the product
\[
-\int_{\partial \omega} \tilde{v}_0
 u_{\#,\overline{k}|\partial \omega}\, d\sigma \int_{\partial\omega}u_{\#,\overline{k}}\rho^i_{0}\,d\sigma\, .
\]
We first note that
\[
\tilde{v}_0= \nu_{\omega}\cdot\nabla v_{\mathrm{m},0|\partial \omega}=\nu_{\omega}\cdot\nabla v^-[\partial \omega,\rho^i_0]_{|\partial \omega}\, .
\]
On the other hand, by Proposition \ref{rhok} and the jump formula for the normal derivative of the single layer potential,
\[
\nu_{\omega}\cdot\nabla v^-[\partial \omega,\rho^i_0]_{|\partial \omega}=\frac{1}{2}\rho^i_0+W^*_{\omega}[ \rho^i_0]=\frac{1}{2}\rho^i_0+\frac{1}{2}\rho^i_0=\rho^i_0\, .
\]
Accordingly,
\[
\int_{\partial \omega} \tilde{v}_0
 u_{\#,\overline{k}|\partial \omega}\, d\sigma=\int_{\partial \omega}u_{\#,\overline{k}}\rho^i_{0}\,d\sigma\, .
\]
By \cite[Proof of Lem.~7.2]{DaMuRo15}, we have
\[
\int_{\partial \omega}u_{\#,\overline{k}}\rho^i_{0}\,d\sigma=\lim_{t\to \infty}\mathsf{u}_{\overline{k}}(t)\, , 
\]
which implies
\[
-\int_{\partial \omega} \tilde{v}_0
 u_{\#,\overline{k}|\partial \omega}\, d\sigma \int_{\partial\omega}u_{\#,\overline{k}}\rho^i_{0}\,d\sigma=-\bigg(\lim_{t\to \infty}\mathsf{u}_{\overline{k}}(t)\bigg)^2\, .
\]
As a consequence, under assumption \eqref{eq:vanu}, by Remark \ref{rem:cepsm} and formula \eqref{eq:capepsfirst}, we can deduce the validity of the following (cf. Theorem \ref{thm:introcap}).
 
 \begin{theorem}\label{thm:cepsmseries}
Let assumption \eqref{eq:vanu} hold. Then
 \begin{equation}\label{cepsmseries}
\begin{split}
\mathrm{Cap}_{\Omega}(\varepsilon \omega,u)=&\varepsilon^{2 \overline{k}}\Bigg(\int_{\mathbb{R}^2 \setminus \overline{\omega}}|\nabla  \mathsf{u}_{\overline{k}}|^2\, dt+\int_{\omega} |\nabla u_{\#,\overline{k}}|^2\, dt-\frac{\Big(\lim_{t\to \infty}\mathsf{u}_{\overline{k}}(t)\Big)^2}{(\lim_{t\to\infty}H^i_0(t)-H^o_0(0)+(2\pi)^{-1}\log|\varepsilon|)}\Bigg) \\
 &+\sum_{n=2\overline{k}+1}^\infty\varepsilon^n\sum_{l=n-2\overline{k}+1}^{n+1} \frac{c_{(n,l)}}{(\lim_{t\to\infty}H^i_0(t)-H^o_0(0)+(2\pi)^{-1}\log|\varepsilon|)^{l}}\, ,
\end{split}
\end{equation} 
for all $\varepsilon\in]-\varepsilon_{\mathrm{c}},\varepsilon_\mathrm{c}[\setminus\{0\}$. 
\end{theorem}

\begin{rem}\label{rem:cepsmseries}
Therefore, by \eqref{cepsmseries} we have
 \begin{equation}\label{cepsmseriesbis}
\mathrm{Cap}_{\Omega}(\varepsilon \omega,u)=\varepsilon^{2 \overline{k}}\Bigg(\int_{\mathbb{R}^2 \setminus \overline{\omega}}|\nabla \mathsf{u}_{\overline{k}}|^2\, dt +\int_{\omega} |\nabla u_{\#,\overline{k}}|^2\, dt\Bigg)+o(\varepsilon^{2\overline{k}}) \qquad \text{as }\varepsilon \to 0\, .
\end{equation}
Moreover, we note that the terms $\int_{\mathbb{R}^2 \setminus \overline{\omega}}|\nabla \mathsf{u}_{\overline{k}}|^2\, dt$ and $\int_{\omega} |\nabla u_{\#,\overline{k}}|^2\, dt$ depend  both on the geometrical properties of the set $\omega$ and on the behavior at $0$ of the function $u$ (but not on $\Omega$).
\end{rem}

\section{Asymptotic expansion of $\lambda_N(\Omega\setminus (\varepsilon \overline{\omega}))$}\label{sec6}

The aim of this section is to obtain an asymptotic expansion of $\lambda_N(\Omega\setminus (\varepsilon \overline{\omega}))$ by combining the results on $\mathrm{Cap}_{\Omega}(\varepsilon \omega,u)$ of Section \ref{sec5} and the approximation fomula  \eqref{eq:asy:eig} for the eigenvalues (see Courtois \cite[Proof of Theorem 1.2]{Co95} and Abatangelo, Felli, Hillairet, and Lena \cite[Theorem 1.4]{AbFeHiLe}).  

To do so, we take $\alpha \in ]0,1[$, $\Omega$ and $\omega$ as in \eqref{e1} and we assume that
\begin{equation}\label{ass:eig}
\begin{split}
&\text{the N-th eigenvalue $\lambda_N(\Omega)$ for the Dirichlet-Laplacian  is simple}\\
&\text{and $u_N$ is a $L^2(\Omega)$-normalized eigenfunction related to $\lambda_N(\Omega)$.}
\end{split}
\end{equation} 

In order to study  $\lambda_N(\Omega\setminus (\varepsilon \overline{\omega}))$ as $\varepsilon \to 0$, by \eqref{eq:asy:eig} we need to consider the behavior of $\mathrm{Cap}_{\Omega}(\varepsilon \omega,u_N)$. By elliptic regularity theory  (see for instance Theorem 1.2, page 205, in \cite{Fr69}), $u_N$ is analytic in a neighborhood of $0$. Next we note that  by \eqref{eq:capepsfirst} we have
\[
\begin{split}
\mathrm{Cap}_{\Omega}(\varepsilon \omega,u_N)=& -\frac{\big(u_N(0)\big)^2}{\lim_{t\to\infty}H^i_0(t)-H^o_0(0)+(2\pi)^{-1}\log|\varepsilon|}\\
&+\varepsilon\bigg(\sum_{n=1}^\infty\varepsilon^{n-1}\sum_{l=0}^{n+1} \frac{c_{(n,l)}}{(\lim_{t\to\infty}H^i_0(t)-H^o_0(0)+(2\pi)^{-1}\log|\varepsilon|)^{l}}\bigg) \quad \text{as $\varepsilon \to 0$}\, ,
\end{split}
\]
where  $\{c_{(n,l)}\}_{\substack{(n,l)\in\mathbb{N}^2\\\;l\le n+1}}$ as in Theorem \ref{capk} and $H^i_0$ and $H^o_0$ are as in Remark \ref{rem:cepsm}. 

Then by formula  \eqref{eq:asy:eig} we immediately deduce the validity of the following well-known result.

\begin{theorem}\label{thm:eig1}
Let assumption \eqref{ass:eig} hold.Then
\begin{equation}\label{eq:eign1:1}
\begin{split}
\lambda_N(\Omega&\setminus (\varepsilon \overline{\omega}))\\&=\lambda_N(\Omega)-\frac{\big(u_N(0)\big)^2}{(2\pi)^{-1}\log\varepsilon }+o\Big(\frac{1}{\log \varepsilon}\Big)\qquad \text{as $\varepsilon \to 0^+$}\, {.}
\end{split}
\end{equation}
\end{theorem}

Clearly, formula \eqref{eq:eign1:1} of Theorem \ref{thm:eig1} in case 
\begin{equation}\label{eq:van:uN}
u_N(0)=0
\end{equation}
reduces to
\[
\begin{split}
\lambda_N(\Omega\setminus (\varepsilon \overline{\omega}))=\lambda_N(\Omega)+o\Big(\frac{1}{\log \varepsilon}\Big)\qquad \text{as $\varepsilon \to 0^+$}\, .
\end{split}
\]
Therefore, if \eqref{eq:van:uN} holds, we would like to obtain a more accurate  asymptotic expansion of  $\lambda_N(\Omega\setminus (\varepsilon \overline{\omega}))$. We now assume that 
\begin{equation}\label{eq:vanuN}
\begin{split}
&\text{there exists $\overline{k} \in \mathbb{N}\setminus \{0\}$ such that}\  D^\gamma u_N(0)=0 \qquad \forall |\gamma| <\overline{k}\\
&\text{and that $D^\beta u_N(0)\neq 0$ for some $\beta \in \mathbb{N}^2$ with $|\beta|=\overline{k}$}\, .
\end{split}
\end{equation}
and we set
\begin{equation}\label{eq:uNsharp}
u_{N,\#,\overline{k}}(t)\equiv\sum_{\substack{(h,j)\in \mathbb{N}^2\\ h+j=\overline{k}}}\frac{\partial_1^h \partial_2^{j}u_{N}(0)}{h! j!}t^h_1 t_2^j\qquad\qquad\quad  \forall t\in  \mathbb{R}^2\,, 
\end{equation}

Moreover, we denote by $\mathsf{u}_{N,\overline{k}}$ the unique solution in $C^{1,\alpha}_{\mathrm{loc}}(\mathbb{R}^2\setminus \omega)$ of
\begin{equation}\label{eq:uNbf}
\left\{
\begin{array}{ll}
\Delta \mathsf{u}_{N,\overline{k}}=0&\text{in }\mathbb{R}^2\setminus\overline{\omega}\,,\\
\mathsf{u}_{N,\overline{k}}(t)=u_{N,\#,\overline{k}}(t)&\text{for all }t\in\partial\omega\,,\\
\sup_{t\in\mathbb{R}^2\setminus\omega}|\mathsf{u}_{N,\overline{k}}(t)|<+\infty\,.
\end{array} 
\right.
\end{equation}

Then by \eqref{cepsmseriesbis} we have
\[
\mathrm{Cap}_{\Omega}(\varepsilon \omega,u_N)=\varepsilon^{2 \overline{k}}\Bigg(\int_{\mathbb{R}^2 \setminus \overline{\omega}}|\nabla \mathsf{u}_{N,\overline{k}}|^2\, dt +\int_{\omega} |\nabla u_{N,\#, \overline{k}}|^2\, dt\Bigg)+o(\varepsilon^{2\overline{k}}) \qquad \text{as }\varepsilon \to 0\, .
\]

Then, again, by formula  \eqref{eq:asy:eig} of Theorem \ref{thm:asy:eig} we  deduce the validity of the following  result (from which we deduce Theorem \ref{thm:our}).

\begin{theorem}\label{thm:eig2}
Let assumptions \eqref{ass:eig}, \eqref{eq:vanuN} hold. Let $u_{N,\#,\overline{k}}$ be as in \eqref{eq:uNsharp}. Let  $\mathsf{u}_{N,\overline{k}}$ be the unique solution in $C^{1,\alpha}_{\mathrm{loc}}(\mathbb{R}^2\setminus \omega)$ of \eqref{eq:uNbf}. Then
\begin{equation}\label{eq:eign2:1}
\begin{split}
\lambda_N(\Omega&\setminus (\varepsilon \overline{\omega}))\\&=\lambda_N(\Omega)+\varepsilon^{2 \overline{k}}\Bigg(\int_{\mathbb{R}^2 \setminus \overline{\omega}}|\nabla \mathsf{u}_{N,\overline{k}}|^2\, dt +\int_{\omega} |\nabla u_{N,\#, \overline{k}}|^2\, dt\Bigg)+o(\varepsilon^{2\overline{k}}) \qquad \text{as }\varepsilon \to 0^+\, .
\end{split}
\end{equation}
\end{theorem}

\begin{rem}We note that in formula \eqref{eq:eign2:1} the term
\begin{equation}\label{eq:eign2:2}
\Bigg(\int_{\mathbb{R}^2 \setminus \overline{\omega}}|\nabla \mathsf{u}_{N,\overline{k}}|^2\, dt +\int_{\omega} |\nabla u_{N,\#, \overline{k}}|^2\, dt\Bigg)
\end{equation}
depends both on the behavior near $0$ of the eigenfunction $u_N$ and on the geometry $\omega$ of the perforation. We emphasize that the way the term in \eqref{eq:eign2:2} depends on $\Omega$ is only through the eigenfunction $u_N$.
\end{rem}

\section{Optimal location of small holes}\label{sec7}

\cbk
Let us now use the above results to discuss how to position a hole in a domain in order to maximize or minimize an eigenvalue. Let $\Omega$ and $\omega$ satisfy the hypotheses \eqref{e1} for a given $\alpha\in ]0,1[$. Moreover, let us assume that the integer $N\ge1$ is such that  $\lambda_N(\Omega)$ is simple.  The small holes we are considering are sets of the form $p+\varepsilon\omega$, where $p\in\Omega$ and $\varepsilon>0$ are such that $p+\varepsilon\overline\omega\subset\overline\Omega$. Let us point out that we do not \emph{a priori} exclude that $\overline\omega$ touches the boundary of $\Omega$.

For a fixed $\varepsilon\in]0.\varepsilon_0[$ (see the comments following Condition \eqref{e1}), we may define a maximum and a minimum problem. More precisely, we can look for two points $p_\varepsilon^M$ and $p_\varepsilon^m$ in $\Omega$, if they exist, such that
\begin{itemize}
\item[(M1)] For each $p \in \Omega$ such that $p+\varepsilon\overline\omega\subset\overline\Omega$ , $\lambda_N(\Omega \setminus (p^M_\varepsilon+\varepsilon \overline{\omega}))\geq \lambda_N(\Omega \setminus (p+\varepsilon \overline{\omega})) $ ;
\item[(m1)] For each $p \in \Omega$ such that $p+\varepsilon\overline\omega\subset\overline\Omega$,  $\lambda_N(\Omega \setminus (p^m_\varepsilon+\varepsilon \overline{\omega}))\leq \lambda_N(\Omega \setminus (p+\varepsilon \overline{\omega})) $ .
\end{itemize} 
These problems are studied in much more detail, when $N=1$, in Section 3.5 of \cite{He06}.
 
 In the rest of this section, we are discussing slightly different problems, which can be understood as an asymptotic version of (M1) and (m1).
More specifically, we would like to find, under the same assumptions, two points $p^M$ and $p^m$ in $\Omega$, if they exist, such that:
\begin{itemize}
\item[(M2)] For each $p \in \Omega$, there exists $\varepsilon_p^M>0$ such that 
\[
\lambda_N(\Omega \setminus (p^M+\varepsilon \overline{\omega}))\geq \lambda_N(\Omega \setminus (p+\varepsilon \overline{\omega})) \qquad \forall \varepsilon \in ]0,\varepsilon_p^M[\, ;
\]
\item[(m2)] For each $p \in \Omega$, there exists $\varepsilon_p^m>0$ such that 
\[
\lambda_N(\Omega \setminus (p^m+\varepsilon \overline{\omega}))\leq \lambda_N(\Omega \setminus (p+\varepsilon \overline{\omega})) \qquad \forall \varepsilon \in ]0,\varepsilon_p^m[\, .
\]
\end{itemize}  
Also we do not have a complete solution of these problems, we wish to present some remarks.

Let us first consider Problem (M2). As before, we denote by $u_N$ a normalized eigenfunction associated  with $\lambda_N(\Omega)$. If the function $|u_N|^2$ has a unique point of maximum $p^*$ in $\Omega$, then $p^*$ is the unique solution of Problem (M2). This follows directly from Theorem \ref{thm:eig1}. If $|u_N|^2$ has more that one point of maximum, a solution of Problem (M2), if it exists, must be one of them. In order to be more precise, we would have to look at higher order terms in the expansions. 

Problem (m2) seems more difficult. Indeed, we can sometime prove that it has no solution. In the case $N=1$, there exists a (unique) positive and normalized eigenfunction associated with $\lambda_1(\Omega)$, which we denote by $u_1$. Since $u_1$ is continuous on $\overline\Omega$ and vanishes on $\partial\Omega$, for any $p\in \Omega$, there exists $q\in\Omega$ such that $0<u_1(q)<u_1(p)$. Using again Theorem \ref{thm:eig1}, it follows that $p$ is not a solution of (m2), showing that the problem has no solution. In the case $N\ge2$, any eigenfunction associated with $\lambda_N(\Omega)$ is orthogonal to $u_1$ and therefore has a non-empty nodal set. Nevertheless, it is still possible that (m2) has no solution. For instance, let us consider the case where the nodal set  of $u_N$ consists of a single simple curve $\gamma$ connecting two points $p_1$ and $p_2$ of $\partial\Omega$. If $p$ belongs to $\gamma$,  by Theorem \ref{thm:eig2},
\begin{equation}
\label{eq:line}
\begin{split}
\lambda_N(\Omega&\setminus (p+\varepsilon \overline{\omega}))\\&=\lambda_N(\Omega)+\varepsilon^{2}\Bigg(\int_{\mathbb{R}^2 \setminus \overline{\omega}}|\nabla \mathsf{u}_{N,1}^{p}|^2\, dt +\int_{\omega} |\nabla u_{N,\#, 1}^{p}|^2\, dt\Bigg)+o(\varepsilon^{2}) \qquad \text{as }\varepsilon \to 0^+\,,
\end{split}
\end{equation}
where $u_{N,\#, 1}^{p}$ and $\mathsf{u}_{N,1}^{p}$ are defined by \eqref{eq:uNsharp} and \eqref{eq:uNbf}, after a translation sending $p$ to $0$. {In other words,
\begin{equation}\label{eq:uN1}
u_{N,\#,{1}}^{p}(t)\equiv \partial_1 u(p) t_1 + \partial_1 u(p) t_2 \qquad\qquad\quad  \forall t\in  \mathbb{R}^2\,, 
\end{equation}
and $\mathsf{u}_{N,1}^{p}$ is the unique solution in $C^{1,\alpha}_{\mathrm{loc}}(\mathbb{R}^2\setminus \omega)$ of
\begin{equation}\label{eq:uN1bis}
\left\{
\begin{array}{ll}
\Delta \mathsf{u}_{N,1}^{p}=0&\text{in }\mathbb{R}^2\setminus\overline{\omega}\,,\\
\mathsf{u}_{N,1}^{p}(t)=u_{N,\#,1}^{p}(t)&\text{for all }t\in\partial\omega\,,\\
\sup_{t\in\mathbb{R}^2\setminus\omega}|\mathsf{u}_{N,1}^{p}(t)|<+\infty\,.
\end{array} 
\right.
\end{equation}}
From this and Theorem \ref{thm:eig1}, it follows that whenever $q\in \Omega\setminus\gamma$ and $\varepsilon>0$ small enough,
\[\lambda_N(\Omega\setminus (p+\varepsilon \overline{\omega}))<\lambda_N(\Omega\setminus (q+\varepsilon \overline{\omega})).\]
On the other hand, since $\nabla u_N$ vanishes at $p_1$ and $p_2$, { we have that $u_{N,\#, 1}^{p}$ and $\mathsf{u}_{N,1}^{p}$ converge to $0$ as $p$ moves on $\gamma$ towards $p_1$ or $p_2$. Accordingly,} the coefficient following $\varepsilon^2$ in Formula \eqref{eq:line} goes to $0$ as $p$ moves on $\gamma$ towards $p_1$ or $p_2$ {(cf. \eqref{eq:uN1} and \eqref{eq:uN1bis})}. It follows that for any fixed $p\in\gamma\cap\Omega$ such that the coefficient is non-zero, we can find $p'\in\gamma\cap\Omega$ such that 
\[\lambda_N(\Omega\setminus (p'+\varepsilon \overline{\omega}))<\lambda_N(\Omega\setminus (p+\varepsilon \overline{\omega}))\]
for all $\varepsilon>0$ small enough. As a result, Problem (m2) has no solution in this case, assuming that $\omega$ is such that the coefficient never vanishes. We will see in Section \ref{sec9} that this last condition is satisfied when $\omega$ is the interior of an ellipse, in particular when $\omega$ is a disk.

If instead $u_N$ has an order of   vanishing greater than one at some points inside $\Omega$, that is to say if at least two nodal lines meet at some points, Problem (m2) may have a solution. As a first step to find it, we need to look for the set $\mathcal{N}$ of those points $p\in \Omega$ where the largest number of nodal lines intersect. We denote this number by $\overline{k}$. Again by Theorem \ref{thm:eig2}, we know that for each $p\in \mathcal{N}$ we have
\[
\begin{split}
\lambda_N(\Omega&\setminus (p+\varepsilon \overline{\omega}))\\&=\lambda_N(\Omega)+\varepsilon^{2 \overline{k}}\Bigg(\int_{\mathbb{R}^2 \setminus \overline{\omega}}|\nabla \mathsf{u}_{N,\overline{k}}^{p}|^2\, dt +\int_{\omega} |\nabla u_{N,\#, \overline{k}}^{p}|^2\, dt\Bigg)+o(\varepsilon^{2\overline{k}}) \qquad \text{as }\varepsilon \to 0^+\,.
\end{split}
\]
{Here above, 
\[
u_{N,\#,\overline{k}}^{p}(t)\equiv\sum_{\substack{(h,j)\in \mathbb{N}^2\\ h+j=\overline{k}}}\frac{\partial_1^h \partial_2^{j}u(p)}{h! j!}t^h_1 t_2^j\qquad\qquad\quad  \forall t\in  \mathbb{R}^2\,, 
\]
and $\mathsf{u}_{N,\overline{k}}^{p}$ is the unique solution in $C^{1,\alpha}_{\mathrm{loc}}(\mathbb{R}^2\setminus \omega)$ of
\[
\left\{
\begin{array}{ll}
\Delta \mathsf{u}_{N,\overline{k}}^{p}=0&\text{in }\mathbb{R}^2\setminus\overline{\omega}\,,\\
\mathsf{u}_{N,\overline{k}}^{p}(t)=u_{N,\#,\overline{k}}^{p}(t)&\text{for all }t\in\partial\omega\,,\\
\sup_{t\in\mathbb{R}^2\setminus\omega}|\mathsf{u}_{N,\overline{k}}^{p}(t)|<+\infty\,.
\end{array} 
\right.
\]}
If $\mathcal N$ contains a single point, this point is the unique solution of (m2). If not, we have to move to a second step: we need to minimize the coefficient in front of $\varepsilon^{2\overline{k}}$. If there exists a point $p^*\in \mathcal{N}$ such that
\[
\Bigg(\int_{\mathbb{R}^2 \setminus \overline{\omega}}|\nabla \mathsf{u}_{N,\overline{k}}^{p^*}|^2\, dt +\int_{\omega} |\nabla u_{N,\#, \overline{k}}^{p^*}|^2\, dt\Bigg)< \Bigg(\int_{\mathbb{R}^2 \setminus \overline{\omega}}|\nabla \mathsf{u}_{N,\overline{k}}^{p}|^2\, dt +\int_{\omega} |\nabla u_{N,\#, \overline{k}}^{p}|^2\, dt\Bigg) \qquad \forall p\in \mathcal{N}\setminus\{p^*\}\, ,
\]
then this point is the unique solution of Problem (m2). If not, we cannot conclude that a solution exists without looking at higher order terms in the expansions.

We have not explored the relation between Problems (M1) and (M2), respectively (m1) and (m2), or even the existence of solutions of Problems (M1) and (m1). Assuming existence, we could for instance ask the following question: if (M2) has a unique solution $p^M$, do solutions $p_\varepsilon^M$ of (M1) converge to $p^M$ when $\varepsilon\to0^+$, and similarly for (m2) and (m1)? Such questions are discussed on page 60 of \cite{He06}, but answering them would require a deeper analysis. In particular, we probably need to understand eigenvalue variation for a hole close $\partial \Omega$ in order to find the connection between (m1) and (m2).
\cbk

\section{Numerical simulations}\label{sec8}

\noindent In this section, we present some numerical simulations on the asymptotic behavior of the eigenvalues in a domain with a small hole. Both the domains $\Omega$ and $\omega$ will have elliptic shapes, but we will consider different rotations of the small hole in order to show the dependence of the asymptotic behavior of the eigenvalues on the geometry of the hole $\varepsilon \omega$ and on the relation of its orientation with respect to the nodal lines of a suitably normalized eigenfunction in the unperturbed domain $\Omega$.

We take $a,b>0$ and we consider the ellipse $\cE_{0}(a,b)$ parametrized by 
\begin{equation}\label{ellisseE0ab}
	\cE_{0}(a,b) = \left\{(x,y)\in\mathbb R^2,  \frac{x^2}{a^2}+\frac{y^2}{b^2}<1 \right\}.
\end{equation}
We denote by $\mathcal R_{\theta}$ the rotation of angle $\theta\in[0,\pi/2]$.
For $\varepsilon>0$ small enough, we define the perforated domain $\cE_{\varepsilon,\theta}(a,b)$ by setting
$$\cE_{\varepsilon,\theta}(a,b) = \cE_{0}(a,b) \setminus \varepsilon\mathcal R_{\theta} \cE_{0}(\frac a4,\frac b4).$$
In other words, we set
\[
\Omega \equiv \cE_{0}(a,b)\, , \qquad \omega\equiv \mathcal R_{\theta} \cE_{0}(\frac a4,\frac b4)\, ,
\]
so that
\[
\Omega \setminus (\varepsilon \overline{\omega})= \cE_{\varepsilon,\theta}(a,b)\, .
\]

In the sequel, we fix $a=3, b=2$ and we omit these parameters in the notation. Namely, 
\[
\cE \equiv \cE_{0}(3,2)\, , \qquad \cE(\varepsilon,\theta)\equiv  \cE_{\varepsilon,\theta}(3,2)\, .
\]
As discretization of the parameters, we choose
$$\varepsilon\in\{1.5^{-k}, 0\leq k\leq 20\}\quad\mbox{ and }\quad\theta\in\left\{\frac{j}{10}\frac\pi2, 0\leq j\leq 10\right\}.$$
We denote by $\lambda_N$ and $\lambda_N(\varepsilon,\theta)$  the $N$-th eigenvalue of the Laplacian  with Dirichlet boundary condition in  the (unperturbed) ellipse $\cE$ and in the perforated domain $\cE(\varepsilon,\theta)$, respectively. 

We first note that the first 16 eigenvalues of the Dirichlet Laplacian in the ellipse  $\cE$ are the following:\\

$$\begin{array}{llllllll}
\lambda_1 =  1.04 & \lambda_2 =  2.13 & \lambda_3 =  3.14 & \lambda_4 =  3.69 &
\lambda_5 =  4.71  & \lambda_6 =  5.74  & \lambda_7 =  6.52 & \lambda_8 =  6.69 \\
\lambda_9 =  8.26 & \lambda_{10} =  8.65 & \lambda_{11} =  9.09 & \lambda_{12} =  11.12 
&\lambda_{13} =  11.21 & \lambda_{14} =  11.25 & \lambda_{15} =  11.92 & \lambda_{16} =  13.82
\end{array}$$

\begin{figure}[h!]\begin{center}
\includegraphics[]{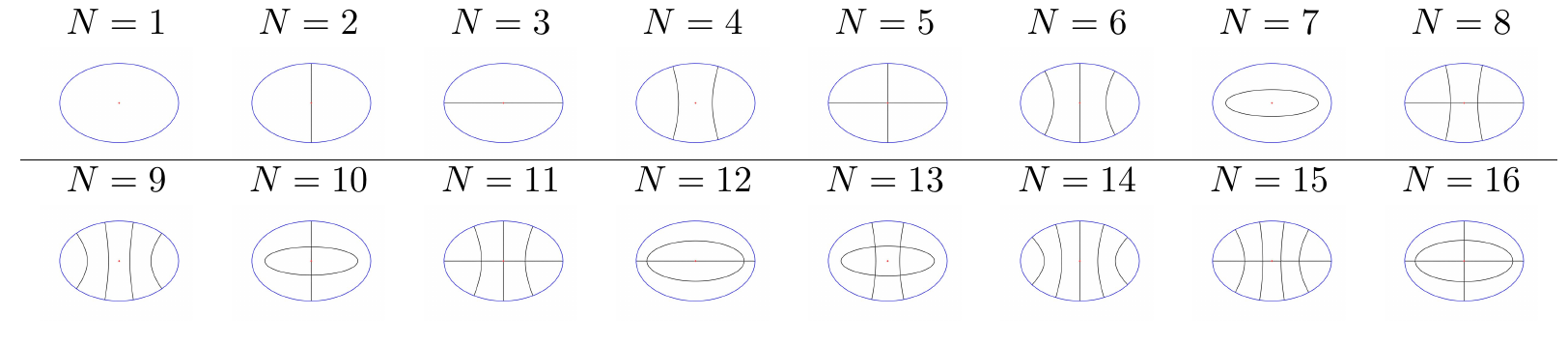}
\caption{Nodal lines of the first 16 eigenfunctions in the ellipse}\label{fig:nodE}
\end{center}\end{figure}

As we can see, all the eigenvalues $\lambda_1, \dots, \lambda_{16}$ are  simple. Then in Figure \ref{fig:nodE} we trace the nodal line of the corresponding eigenfunctions.

Figure \ref{fig:nodE} shows that the origin $0$ (which is the point where the hole collapses when $\varepsilon =0$) belongs to a nodal line of an  eigenfunction associated to the eigenvalue $\lambda_N$ for $N\in \{2, 3, 5, 6, 8, 10, 11, 12, 14, 15, 16\}$.


As a preliminary step in our analysis, we look at the behavior of the nodal lines in the perforated domain as $\varepsilon$ approaches $0$ and for different values of the angle $\theta$. We note that for $N\in \{2,3,6,8,10,12,14,15\}$ there is only one nodal line passing through $0$, whereas for $N\in \{5,11,16\}$ there are two nodal lines. As a consequence, in view of the results of Section \ref{sec6}, we expect that
\[
\begin{aligned}
&\lambda_N(\varepsilon,\theta)-\lambda_N \sim -\frac{c_N (\theta)}{\log \varepsilon} \qquad \mbox{ if } N  \in  \{1,4,7,9,13\}\, ,\\
&\lambda_N(\varepsilon,\theta)-\lambda_N \sim  c_N (\theta)\varepsilon^2 \qquad \mbox{ if } N  \in \{2,3,6,8,10,12,14,15\}\, ,\\
&\lambda_N(\varepsilon,\theta)-\lambda_N \sim c_N (\theta) \varepsilon^4  \qquad \mbox{ if } N  \in \{5,11,16\}\, ,
\end{aligned}
\]
as $\varepsilon \to 0^+$, for some  constant $c_N (\theta)>0$ which depends on $N $ and $\theta$. 

\begin{figure}[h!]\begin{center}
\includegraphics[]{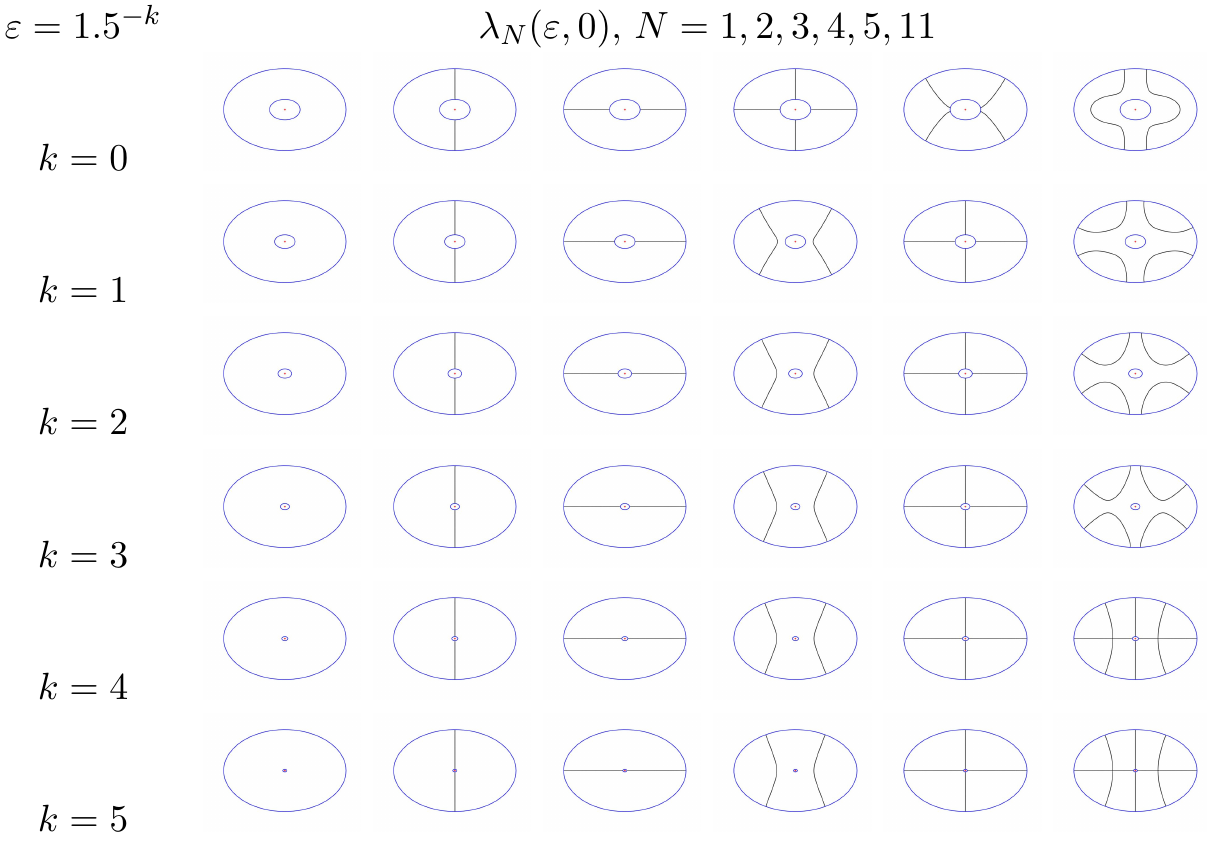}
\caption{Nodal lines of the eigenfunctions associated to  $\lambda_N(\varepsilon,\theta)$, $N =1,2,3, 4, 5, 11$, $\theta=0$}\label{fig:nodE0}
\end{center}\end{figure}

\begin{figure}[h!]\begin{center}
\includegraphics[]{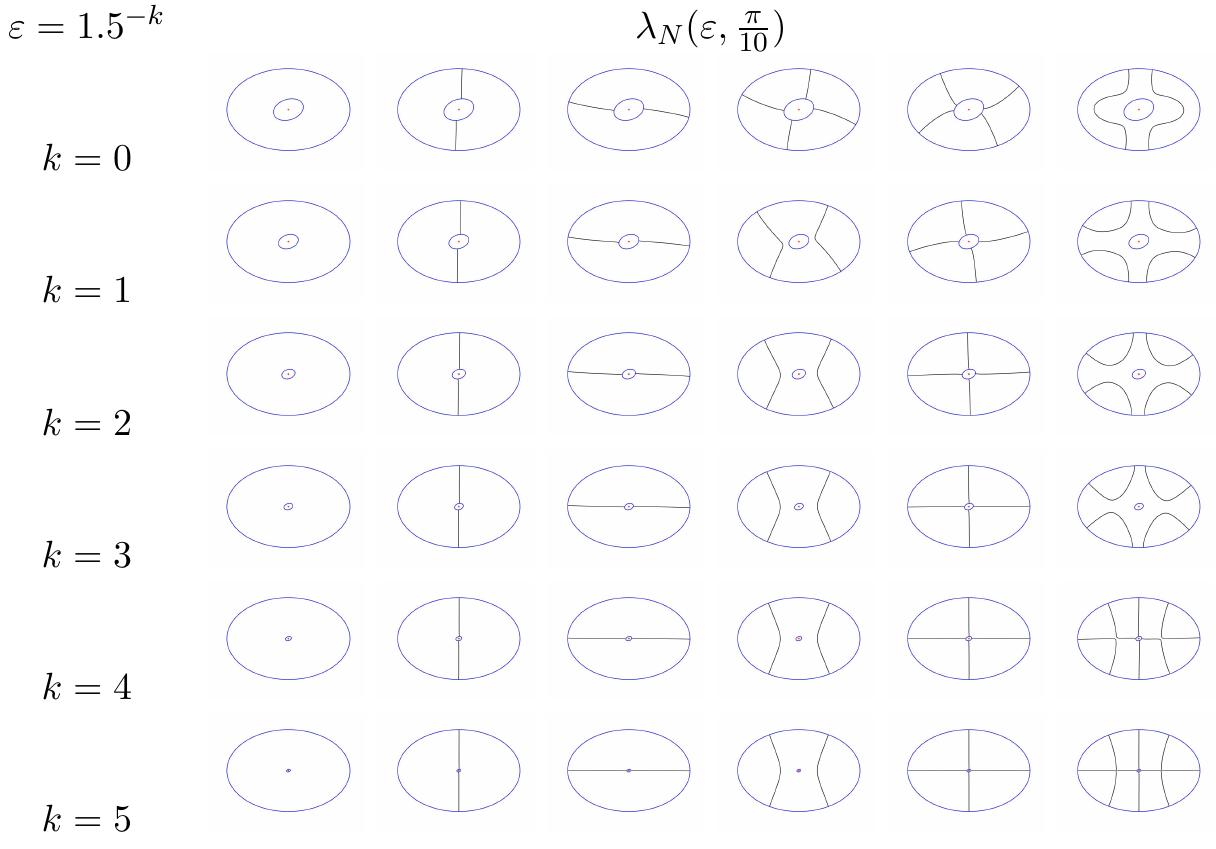}
\caption{Nodal lines of the eigenfunctions associated to  $\lambda_N(\varepsilon,\theta)$, $N =1,2,3, 4, 5, 11$, $\theta=\frac\pi{10}$}\label{fig:nodE10}
\end{center}\end{figure}

\begin{figure}[h!]\begin{center}
\includegraphics[]{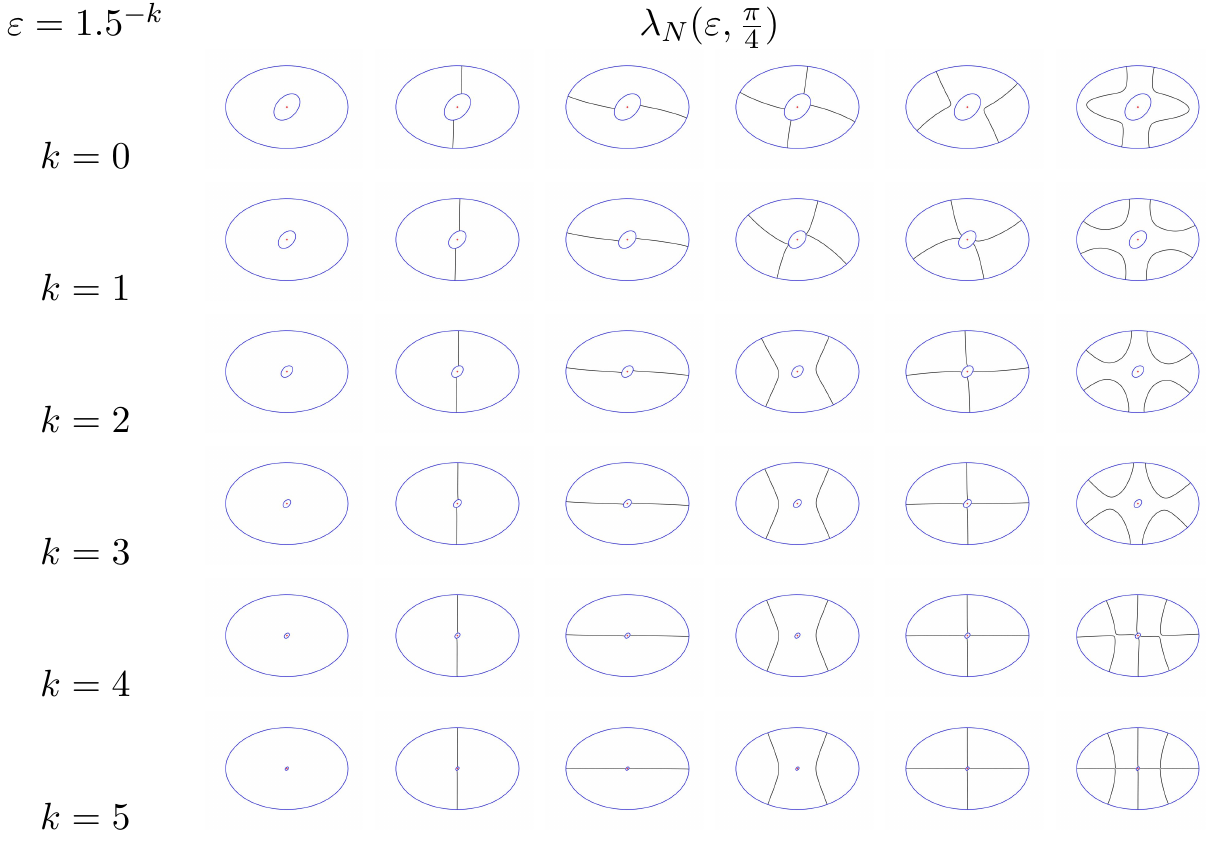}
\caption{Nodal lines of the eigenfunctions associated to  $\lambda_N(\varepsilon,\theta)$, $N =1,2,3, 4, 5, 11$, $\theta=\frac\pi4$}\label{fig:nodE4}
\end{center}\end{figure}

\begin{figure}[h!]\begin{center}
\includegraphics[]{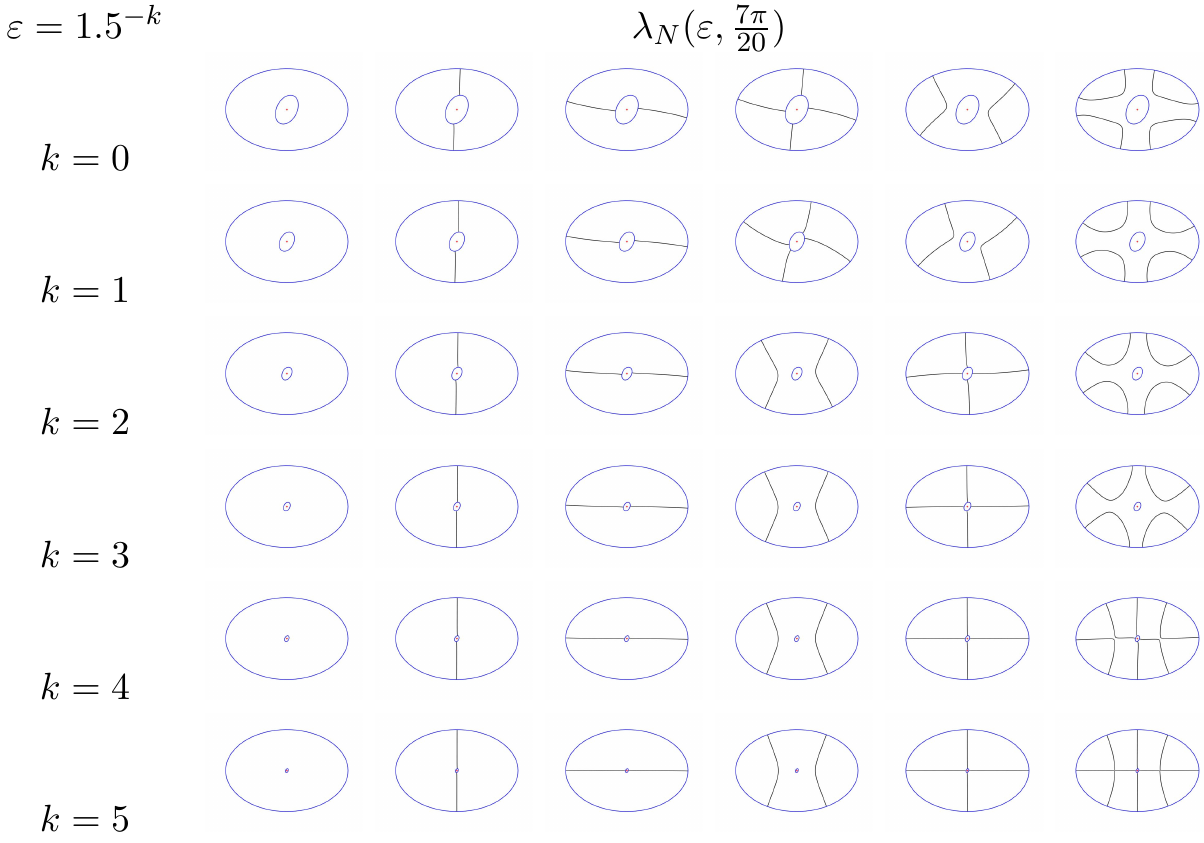}
\caption{Nodal lines of the eigenfunctions associated to  $\lambda_N(\varepsilon,\theta)$, $N =1,2,3, 4, 5, 11$, $\theta=\frac{7\pi}{20}$}\label{fig:nodE720}
\end{center}\end{figure}

\begin{figure}[h!]\begin{center}
\includegraphics[]{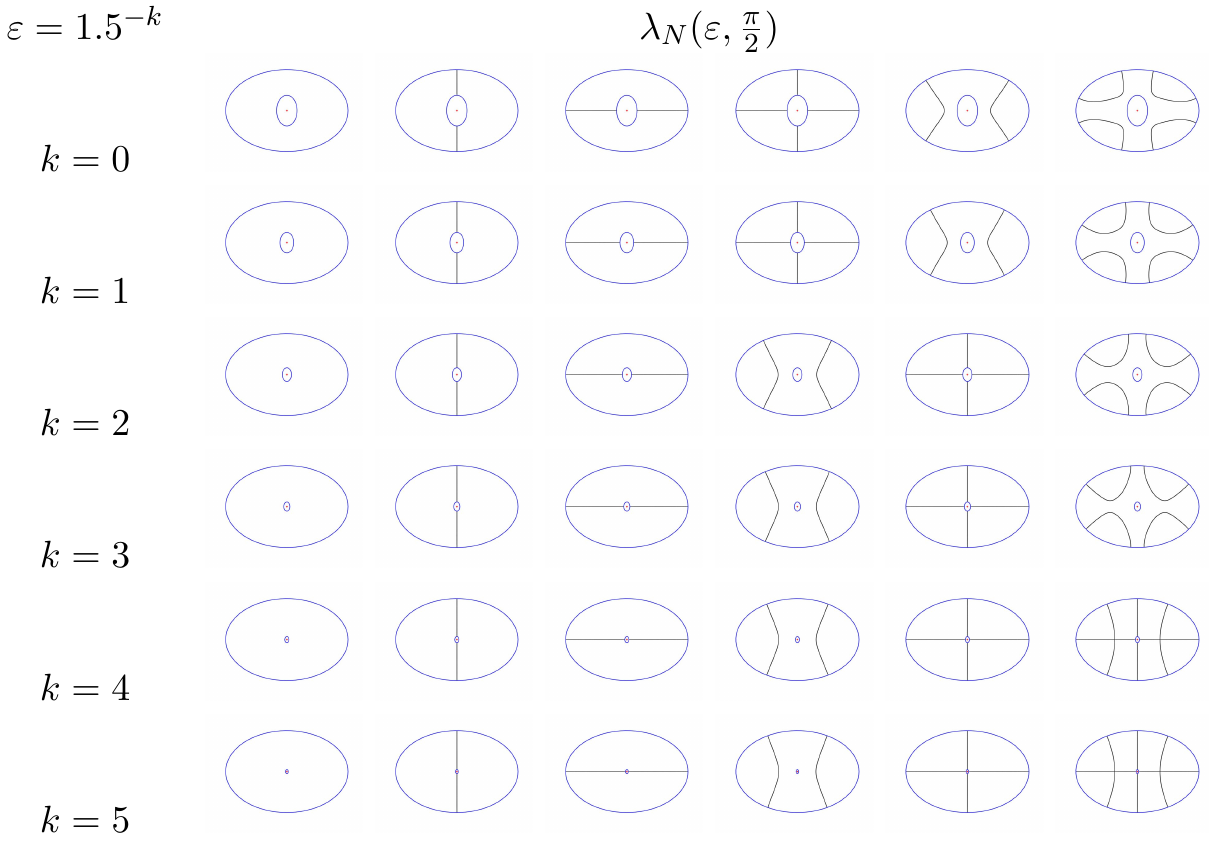}
\caption{Nodal lines of the eigenfunctions associated to  $\lambda_N(\varepsilon,\theta)$, $N =1,2,3, 4, 5, 11$, $\theta=\frac\pi2$}\label{fig:nodE2}
\end{center}\end{figure}

Figures \ref{fig:nodE0}--\ref{fig:nodE2} show how the angle $\theta$ affects the convergence of the nodal lines as $\varepsilon$ tends  to $0$.

\newpage

\newpage

\newpage

Then we turn to consider the behavior of the functions $\lambda_N(\varepsilon,\theta)-\lambda_N$ as $\varepsilon$ approaches $0$, for different values of the angle $\theta$ (namely for $\theta \in \{0,\frac{\pi}{10}, \frac{\pi}{4}, \frac{7\pi}{20}, \frac{\pi}{2}\}$) and for $N  \in \{2,3,5,11\}$. As already mentioned, in view of Section \ref{sec6}, we expect
\[
\begin{aligned}
&\lambda_N(\varepsilon,\theta)-\lambda_N \sim  c_N (\theta)\varepsilon^2 \qquad \mbox{ for } N \in\{2,3\}\, ,\\
&\lambda_N(\varepsilon,\theta)-\lambda_N \sim c_N (\theta) \varepsilon^4 \qquad \mbox{ for } N \in\{5,11\}\, ,
\end{aligned}
\]
as $\varepsilon \to 0^+$, for some  constant $c_N (\theta)>0$ which depends on $N $ and $\theta$. 

Figures \ref{fig:eigE0} -- \ref{fig:eigE2} below show in a $\log$--$\log$ plot a good fitting with the expected behavior.
To compute the eigenmodes, we use a finite element method of degree $\mathbb P_{4}$ with at least 2800 triangular elements. We work in a simple precision, so computations are relevant when the gap $\log(\lambda_N(\varepsilon,0)-\lambda_N)$ is larger than $10^{-8}$. It is the reason the computations for $N=5,11$ and $\varepsilon \leq 10^{-2}$ are irrelevant.
Moreover, for $\varepsilon$ fixed, they show  a decreasing behavior in $\theta \in [0,\pi/2]$ for $N =2$, in contrast with an increasing behavior in $\theta \in [0,\pi/2]$ for $N =3$. Instead, for $N =5$ and $\varepsilon$ is fixed and small, the quantity $\lambda_N(\varepsilon,\theta)-\lambda_N$ is first increasing and then decreasing. 

\begin{figure}[h!]
\begin{center}
\includegraphics[]{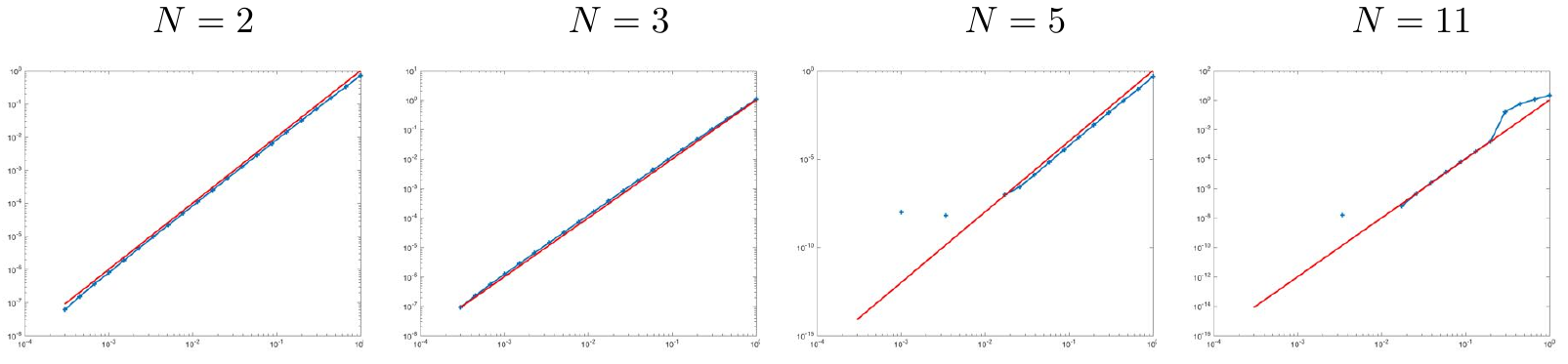}
\caption{In blue, the plot of $\log\varepsilon\mapsto \log(\lambda_N(\varepsilon,0)-\lambda_N)$; in red, for $N =2,3$ the plot of $\log \varepsilon\mapsto \log \varepsilon^2$ and for $N =5,11$ the plot of $\log \varepsilon\mapsto \log  \varepsilon^4$}\label{fig:eigE0}
\end{center}
\end{figure}


\begin{figure}[h!]
\begin{center}
\includegraphics[]{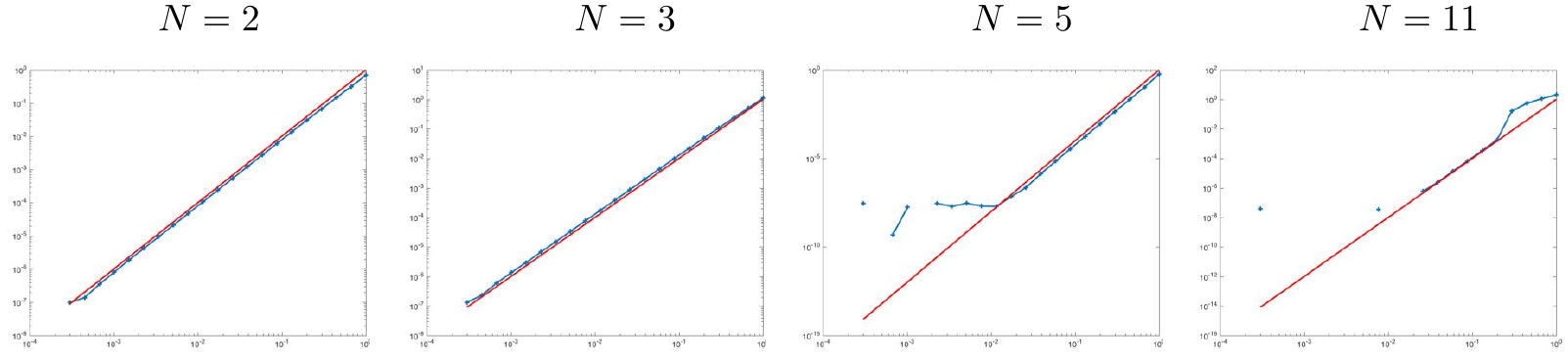}
\caption{In blue, the plot of  $\log\varepsilon\mapsto \log(\lambda_N(\varepsilon,\theta)-\lambda_N)$, $\theta=\frac\pi{10}$;  in red, for  $N =2,3$ the plot of $\log \varepsilon\mapsto \log \varepsilon^2$ and for $N =5,11$ the plot of $\log \varepsilon\mapsto \log \varepsilon^4$.}\label{fig:eigE10}
\end{center}
\end{figure}

\begin{figure}[h!]
\begin{center}
\includegraphics[]{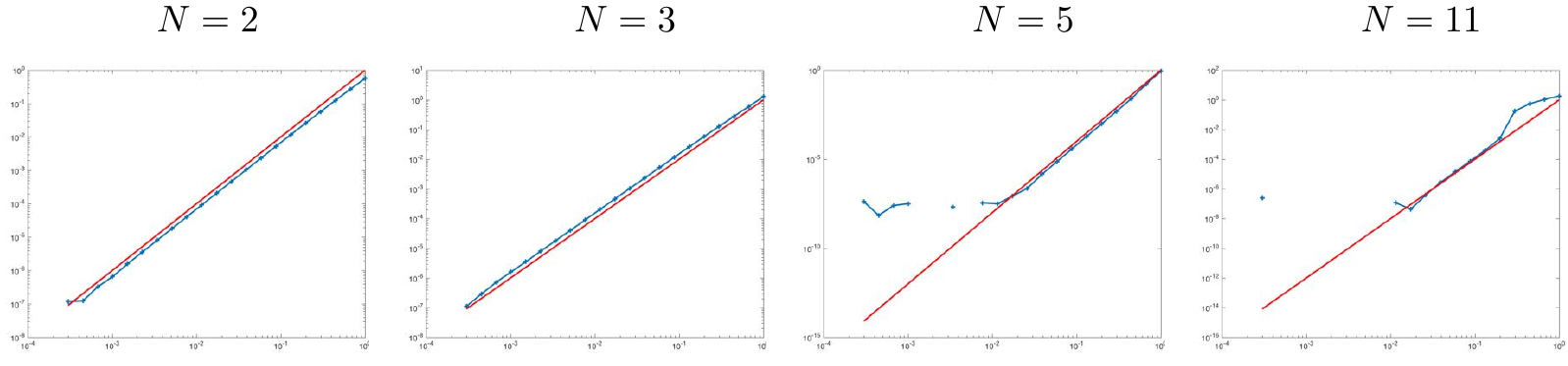}
\caption{In blue, the plot of $\log\varepsilon\mapsto \log(\lambda_N(\varepsilon,\theta)-\lambda_N)$, $\theta=\frac\pi{4}$; in red, for $N =2,3$  the plot of $\log \varepsilon\mapsto \log \varepsilon^2$ and for $N =5,11$ the plot of $\log \varepsilon\mapsto \log \varepsilon^4$.}\label{fig:eigE4}
\end{center}
\end{figure}


\begin{figure}[h!]
\begin{center}
\includegraphics[]{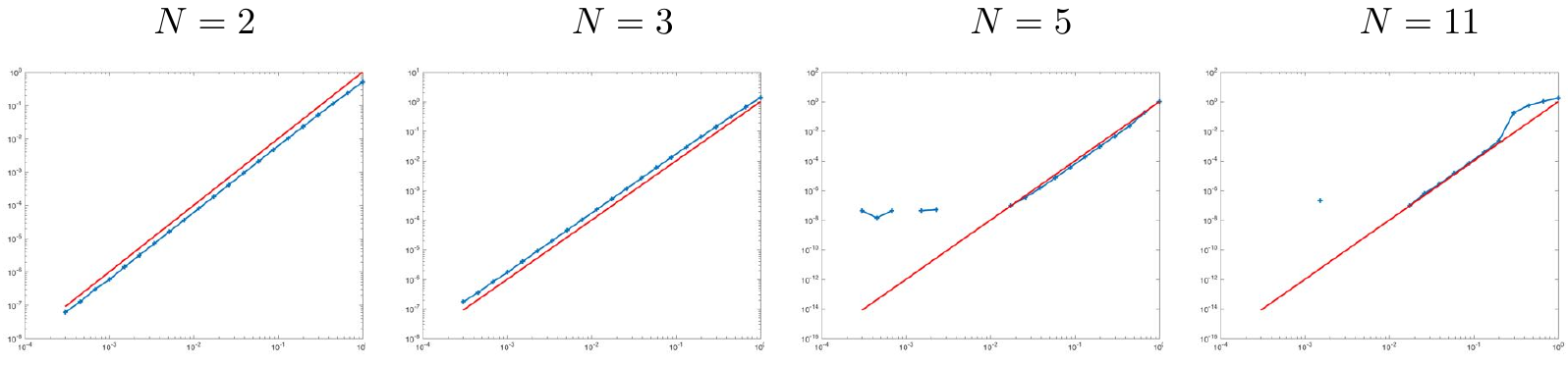}
\caption{In blue, the plot of $\log\varepsilon\mapsto \log(\lambda_N(\varepsilon,\theta)-\lambda_N)$, $\theta=\frac{7\pi}{20}$; in red, for $N =2,3$ the plot of $\log \varepsilon\mapsto \log \varepsilon^2$ and for $N =5,11$ the plot of $\log \varepsilon\mapsto \log \varepsilon^4$.}\label{fig:eigE720}
\end{center}
\end{figure}


\begin{figure}[h!]
\begin{center}
\includegraphics[]{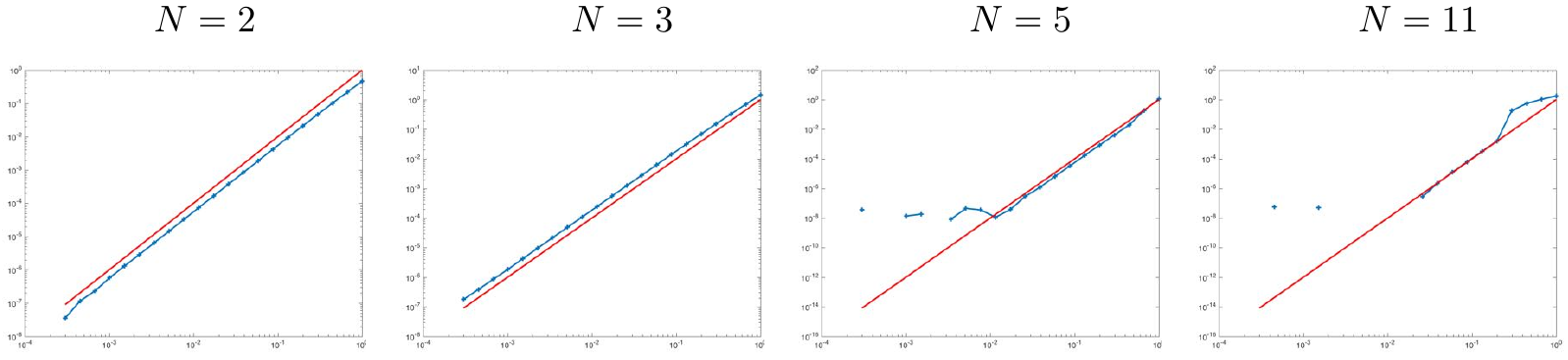}
\caption{In blue the plot of $\log\varepsilon\mapsto \log(\lambda_N(\varepsilon,\theta)-\lambda_N)$, $\theta=\frac\pi{2}$; in red, for $N =2,3$ the plot of $\log \varepsilon\mapsto \log \varepsilon^2$ and for $N =5,11$ the plot of $\log \varepsilon\mapsto \log \varepsilon^4$.}\label{fig:eigE2}
\end{center}
\end{figure}

Finally, we study the limiting behavior as a function of the angle $\theta$. To do so, we set
$$\mu_N (\varepsilon,\theta) = \frac{1}{\varepsilon^{\alpha_N }} \left(\lambda_N(\varepsilon,\theta)-\lambda_N\right),$$
where $\alpha_N $ is the order of the second term in the asymptotic expansion, {\it i.e.}:
$$\alpha_{2}=\alpha_{3}=2, \qquad \alpha_{5}=\alpha_{11}=4.$$
We plot the curve 
$$\theta\mapsto \mu_N (\varepsilon,\theta),$$
with $\varepsilon=1.5^{-k}$ and $k=14$ or $16$ when $N =2,3$ (in this case, we have a convergence at order $2$)
and $k=7$ or $8$ when $N =5, 11$.
This choice is done to ensure that   $\alpha_N k$ is constant in both cases.

\begin{figure}[h!]
\begin{center}
\includegraphics[]{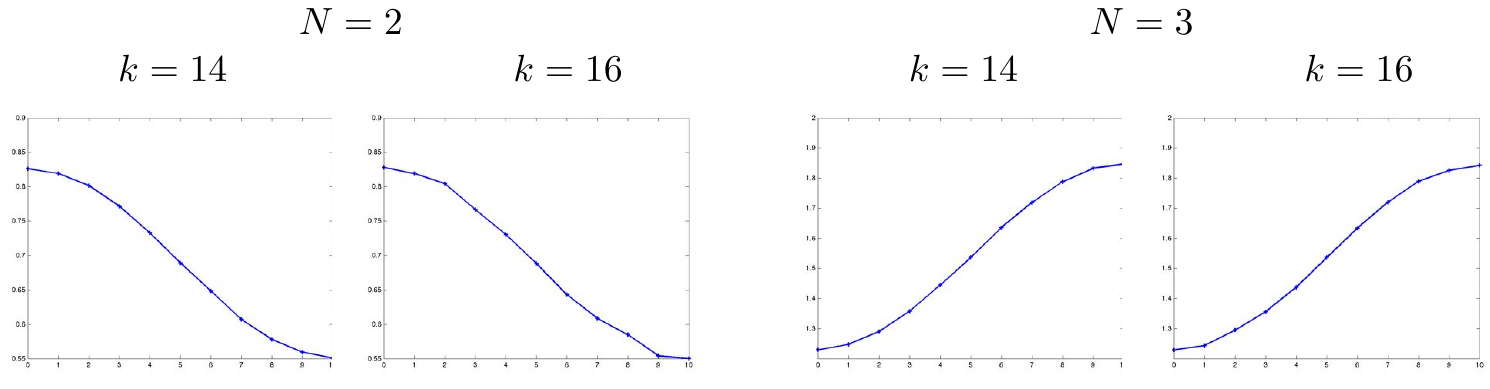}
\caption{$\theta\mapsto \mu_N (1.5^{-k},\theta)$, $N=2,3$, $k=14, 16$}\label{fig:plot1}
\end{center}
\end{figure}

\begin{figure}[h!]
\begin{center}
\includegraphics[]{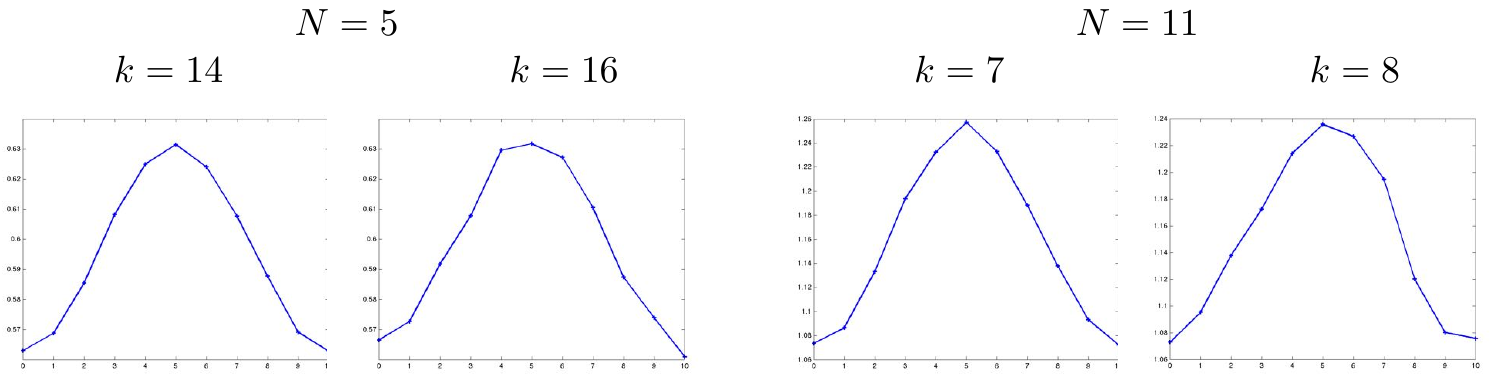}
\caption{$\theta\mapsto \mu_N (1.5^{-k},\theta)$, $N=5,11$, $k=7, 8$}\label{fig:plot2}
\end{center}
\end{figure}

Figure \ref{fig:plot1} confirms that for $\varepsilon$ fixed  the function $[0,\frac{\pi}{2}]\ni \theta \mapsto \mu_N (\varepsilon,\theta)$ is decreasing for $N =2$ and increasing for $N =3$. Instead, for $N =5, 11$, Figure \ref{fig:plot2} shows that the function $[0,\frac{\pi}{2}]\ni \theta \mapsto \mu_N (\varepsilon,\theta)$ is first increasing and then decreasing. 

\section{Theoretical analysis of the simulations}\label{sec9}

In this {last} section we are going to prove analytically what we have seen in Section \ref{sec8}, that is the dependence of simple eigenvalues' behavior on the angle between $x_1$-axis and the small ellipse's major axis.

If we consider the ellipse $\mathcal E_0$ in \eqref{ellisseE0ab} with $a>b>0$, it can be written as
\[
 \mathcal E_0(a,b) = \Big\{ (x_1,x_2)\in \R^2:\ \frac{x_1^2}{b^2+c^2} + \frac{x_2^2}{b^2} < 1 \Big\},
\]
where $c$ is  the distance between the two foci, which satisfies $c^2=a^2-b^2$.
Up to replacing $\eps/4$ with $\eps$, we can think 
\[
 \mathcal E(\eps,\theta) = \mathcal E_0(a,b) \setminus \eps\mathcal R_\theta\mathcal E_0(a,b),
\]
being $\mathcal R_\theta$ the rotation of angle $\theta \in [0,\pi/2]$, as in Section 8. 

In view of Theorem \ref{thm:eig2}, we aim at computing the quantity 
\begin{equation}\label{eq:aim}
\mbox{E}\left(\omega,u_{N,\#,\overline{k}}\right):= \int_{\mathbb{R}^2 \setminus \overline{\omega}}|\nabla \mathsf{u}_{N,\overline{k}}|^2\, dt +\int_{\omega} |\nabla u_{N,\#, \overline{k}}|^2\, dt
\end{equation}
which is indeed the coefficient of the leading term of the eigenvalues' difference expansion.  
We remark that in this case 
\begin{equation}\label{eq:omegasim}
 \omega = \mathcal R_\theta\mathcal E_0(a,b),
\end{equation}
depending on $\theta$, so that even \eqref{eq:aim} is in fact depending on $\theta$. An explicit computation of it will show how it depends on this angle.  

As seen in Section \ref{sec6}, the eigenfunction $u_N$ is analytic in a neighborhood of $0$. In this Section, we assume $u_N(0)=0$ (cf. assumption \eqref{eq:vanuN}).  Accordingly,  there exist $\overline k\in \N\setminus \{0\}$ and $P_{\overline k}$, a homogeneous polynomial of degree $\overline k$ in  two variables, such that 
\begin{equation*}
	u_N(x)=P_{\overline k}(x)+O\left(\vert x\vert^{\overline k+1}\right).
\end{equation*}
It follows, from differentiating the series expansion of $u_N$ at $0$, that
\begin{equation*}
	\Delta u_N(x)=\Delta P_{\overline k}(x)+O\left(\vert x\vert^{\overline k-1}\right),
\end{equation*}
and since $(\Delta+\lambda_N) u_N=0$, we obtain $\Delta P_{\overline k}=0$, that is to say the polynomial $P_{\overline k}$ is harmonic. 
Therefore, there exists 
$\beta\in\R\setminus\{0\}$ and $\alpha\in\cbk]-\frac\pi2,\frac\pi2]\cbk$ such that, 
\begin{align}\label{eq:orderk}
r^{-\overline k} u_N(r\cos t,r\sin t ) \to \beta
\sin(\overline k t+\alpha) \text{ as }r\to 0 \text{ in }C^{1,\tau}([0,2\pi]) 
\end{align}
for any $\tau\in\cbk]0,1[\cbk$. 

Moreover, as noted in \cite{AFL2019}, $\beta$ is directly linked to the norm of the ${\overline k}$-th differential of $u_N$ at $0$. More precisely, if we consider 
\[
 \| d^{j}u (x) \|^2 := \sum_{i_1,\ldots,i_j=1}^2 \left| \frac{\partial ^j u}{\partial x_{i_1}\ldots \partial x_{i_j}}(x) \right|^2 ,
\]
then 
\[
\beta^2 = \frac{\| d^{{\overline k}}u_N (0) \|^2}{(\overline k!)^2 \,2^{{\overline k}-1}}.
\]
For the sake of simplicity and without loss of generality, we perform a change of variables by rotating  the domain,  in such a way that
\begin{enumerate} 
	\item[(i)] in the new domain, the major axis of the  small elliptic hole is lying along the $x_1$-axis, so that Equation \eqref{eq:omegasim} reads
		\begin{equation} 
			\omega = \mathcal E_0(a,b) \label{eq:omegasimnew};
		\end{equation}
	\item[(ii)] Equation \eqref{eq:orderk} now reads 
		\begin{equation}
			r^{-\overline k} u_N(r\cos t,r\sin t ) \to \beta\sin(\overline k t + \overline k\varphi) \qquad \text{ as }r\to 0 \text{ in }C^{1,\tau}([0,2\pi]), 				\label{eq:orderknew}
		\end{equation}
	with $\varphi\in \cbk]-\pi/2\overline k,\pi/2\overline k]\cbk$.
\end{enumerate}
\begin{rem}\label{rem:phi}
Given the above condition, $\varphi$ is unique and $-\varphi$ is, in absolute value, the smallest angle at the origin between the major axis of $\mathcal R_\theta\mathcal E_0(a,b)$ and  a nodal line of the eigenfunction $u_N$. We denote this unique angle $\varphi$ by $\varphi(u_N,\theta)$.
\end{rem}

In order to compute explicitly the quantity in \eqref{eq:aim} under assumptions \eqref{eq:omegasimnew} and \eqref{eq:orderknew}, we define the elliptic coordinates $(\xi,\eta)$ (see for instance \cite{Sokolov2011} or \cite{AbFeHiLe,AFL2019}) by
\begin{equation}\label{eq:ellcoord}
\begin{cases}
x_1=c\cosh(\xi)\cos(\eta),\\
x_2=c\sinh(\xi)\sin(\eta),
\end{cases}
\quad \cbk \xi\in [0, +\infty[,\ \eta\in[0,2\pi[.\cbk
\end{equation}
The boundary $\partial\omega=\partial\mathcal E_0(a,b)$ has equation $\xi = \bar\xi$, 
where $\bar \xi$ is defined by the relation 
$c\sinh(\bar\xi)=b$, that is $\bar\xi = \log\left( \frac{b}{c} + \sqrt{1+\frac{b^2}{c^2}}\right)$.

More precisely, we are considering the function $F: (\xi,\eta)\mapsto (x_1,x_2)$ defined by \eqref{eq:ellcoord}. It is a $C^\infty$ diffeomorphism from $\cbk D:= [0,+\infty[\times [0,2\pi[$ onto $\R^2$. $F$ is actually a conformal map, as noted in \cite[Subsection 3.2]{AFL2019}. Let us denote $\mathcal D^{1,2}(\R^2)$ the functions space which is the closure of $C^\infty_c(\R^2)$ with respect to the $L^2$ norm of the gradient.  
For any function $u\in \mathcal D^{1,2}(\R^2)$, let us define $U:= u\circ F$. Since $F$ is conformal, $|\nabla U|\in L^2(D)$ with $\int_{D} |\nabla U|^2\,d\xi d\eta = \int_{\R^2} |\nabla u|^2\,dx_1\, dx_2$ and $U$ is harmonic in $\tilde D\subseteq D$ if and only if $u$ is harmonic in $F(\tilde D)$.   

Let us now denote 
\[
 \psi_{\overline k}^\varphi (r\cos t,r\sin t) := \beta\,r^{\overline k} \sin(\overline k t+\overline k\varphi) \quad \text{for }r>0,\ t\in\cbk [0,2\pi[ 
\]
and define the complex variables $z:= x_1+ i x_2$ and $\zeta:= \xi+ i \eta$. Then we have
\[
 \psi_{\overline k}^\varphi (x_1,x_2) = \textrm{Im}(\beta e^{i{\overline k}\varphi}z^{\overline k});
\]
since $z = F(\xi,\eta)= c\cosh(\zeta)$ and taking into account the Binomial Theorem we obtain  
\begin{align*}
 \Psi_{\overline k}^\varphi (\xi,\eta) := \big(\psi_{\overline k}^\varphi \circ F\big) (\xi,\eta) &= \textrm{Im}(\beta e^{i{\overline k}\varphi}(c \cosh \zeta)^{\overline k}) = \textrm{Im} \left( \frac{\beta c^{\overline k} e^{i{\overline k}\varphi}}{2^{\overline k}} \sum_{j=0}^{\overline k} \left(\begin{array}{c}{\overline k}\\ j\end{array}\right) e^{({\overline k}-2j)\zeta}\right) \\
 &= \frac{\beta c^{\overline k}}{2^{\overline k}} \sum_{j=0}^{\overline k} \left(\begin{array}{c}{\overline k}\\ j\end{array}\right) e^{({\overline k}-2j)\xi} \sin\big(({\overline k}-2j)\eta+{\overline k}\varphi\big)
\end{align*}

In this way, the first contribution in \eqref{eq:aim} is precisely
\begin{equation}\label{eq:energyW}
 \int_{\mathbb{R}^2 \setminus \overline{\omega}}|\nabla \mathsf{u}_{N,\overline{k}}|^2\, dt =
 \int_{\cbk]\bar\xi,+\infty[\times ]0,2\pi[} |\nabla W_{N,\overline k}|^2\,d\xi d\eta
\end{equation}
where $W_{N,\overline k}$ is the unique solution in $C^{1,\alpha}_{\mathrm{loc}}\big(\cbk] \bar\xi,+\infty[\times [0,2\pi[ \cbk \big)$ to the problem 
\begin{equation}\label{eq:W}
 \begin{cases}
  -\Delta W_{N,\overline k} =0, &\text{in }\cbk]\bar\xi,+\infty[\times [0,2\pi[\cbk,\\
  W_{N,\overline k}(\xi,\eta)=\Psi_{\overline k}^\varphi (\bar \xi,\eta) , &\text{on }\xi=\bar\xi,\\
  \displaystyle\sup_{(\xi,\eta)\in\cbk ]\bar\xi,+\infty[\times [0,2\pi[\cbk} |W_{N,\overline k}(\xi,\eta)| < +\infty \\
  W(\xi,0)=W(\xi,2\pi), &\text{for all } \xi\in\cbk]\bar\xi,+\infty[\cbk.
 \end{cases}
\end{equation}
which is the analogous of problem \eqref{eq:uNbf} in elliptic coordinates, that is $W_{N,\overline k}= \mathsf{u}_{N,\overline{k}} \circ F $. 

As well, the second contribution in \eqref{eq:aim} is 
\[
 \int_{\omega} |\nabla u_{N,\#, \overline{k}}|^2\, dt =
 \int_{\cbk ]0,\bar\xi[\times ]0,2\pi[} |\nabla \Psi_{\overline k}^\varphi|^2 \,d\xi d\eta,
\]
that is $ \Psi_{\overline k}^\varphi = u_{N,\#, \overline{k}} \circ F $, since $u_{N,\#, \overline{k}}=\psi_{\overline k}^\varphi$ in view of \eqref{eq:uNsharp}.

\subsection{Computation of the first contribution}

In order to compute the first contribution, we need to compute explicitely the potential $W_{N,\overline k}$ solution to \eqref{eq:W}.  
Let us consider the Fourier expansion of $W$ in elliptic coordinates:
\[
 W_{N,\overline k}(\xi,\eta)= \dfrac{a_0(\xi)}{2} + \sum_{j\ge1} \left(a_j(\xi)\cos(j\eta) + b_j(\xi)\sin(j\eta)\right)
\]
where
\begin{align*}
 &a_j(\xi)= \frac1\pi \int_0^{2\pi} W_{N,\overline k}(\xi,\eta) \cos(j\eta) \,d\eta \quad \text{ for }j\in\N,\\ 
 &b_j(\xi)= \frac1\pi \int_0^{2\pi} W_{N,\overline k}(\xi,\eta) \sin(j\eta) \,d\eta \quad \text{ for }j\in\N\setminus\{0\}.
\end{align*}
Therefore we have
\[
 0=-\Delta_{(\xi,\eta)}W_{N,\overline k} = \dfrac{a_0''(\xi)}{2} 
 +  \sum_{j\ge1} \left((a_j''(\xi)-j^2a_j(\xi))\cos(j\eta) + (b_j''(\xi)-j^2b_j(\xi))\sin(j\eta)\right).
\]
Imposing the boundary conditions for $\xi\in\cbk]\bar\xi,+\infty[$, the latter equation implies 
\begin{equation}
 \begin{cases}
  a_j''(\xi)-j^2 a_j(\xi)=0 \quad \text{for }\xi\geq\bar\xi\\
  a_j(\bar\xi) = \frac1\pi\int_0^{2\pi} \Psi_{\overline k}^\varphi(\bar\xi,\eta)\cos(j\eta)\,d\eta\\
  \sup_{\xi\geq\bar\xi} |a_j(\xi)| < +\infty 
 \end{cases}
 \qquad \text{and}\qquad
 \begin{cases}
  b_j''(\xi)-j^2 b_j(\xi)=0 \quad \text{for }\xi\geq\bar\xi\\
  b_j(\bar\xi) = \frac1\pi\int_0^{2\pi} \Psi_{\overline k}^\varphi(\bar\xi,\eta)\sin(j\eta)\,d\eta\\
  \sup_{\xi\geq\bar\xi} |b_j(\xi)| < +\infty 
 \end{cases}
\end{equation}
for any $j\in\N$ and any $j\in\N\setminus \{0\}$, respectively.  
We solve the latter problems by 
\begin{align*}
 &a_0(\xi) \equiv a_0(\bar\xi) \quad \text{for }\xi\geq\bar\xi;\\
 &a_j(\xi) = a_j(\bar\xi) \,e^{-j(\xi-\bar\xi)} \quad \text{for }\xi\geq\bar\xi, \quad \text{for }j\geq1;\\
 &b_j(\xi) = b_j(\bar\xi) \,e^{-j(\xi-\bar\xi)} \quad \text{for }\xi\geq\bar\xi, \quad \text{for }j\geq1.
\end{align*}
By rewriting Formula \eqref{eq:ibp} in the elliptic coordinates $(\xi,\eta)$, we obtain
\begin{align}\label{eq:pos}
 &\int_{\cbk]\bar\xi,+\infty[\times ]0,2\pi[} |\nabla W_{N,\overline k}|^2\,d\xi d\eta 
 = \int_0^{2\pi} -\frac{\partial W_{N,\overline k}}{\partial\xi}(\bar\xi,\eta)\,W_{N,\overline k}(\bar\xi,\eta) \,d\eta \notag\\
 &\qquad= \int_0^{2\pi} \sum_{j,l\geq 1} j\left( a_j(\bar\xi)\cos(j\eta) + b_j(\bar\xi)\sin(j\eta) \right) \left( a_l(\bar\xi)\cos(l\eta) + b_l(\bar\xi)\sin(l\eta) \right)\,d\eta \notag\\
 &\qquad = \pi \sum_{j\geq 1} j\left( a_j^2(\bar\xi) + b_j^2(\bar\xi) \right). 
\end{align}
In order to conclude the analysis on this first contribution, let us compute the quantities $a_j^2(\bar\xi)$ and $b_j^2(\bar\xi)$. 
By definition, for any $j\geq 1$
\begin{align}\label{eq:aj}
 a_j(\bar\xi) &= \frac1\pi \int_0^{2\pi} W_{N,\overline k}(\bar\xi,\eta)\,\cos(j\eta)\,d\eta \notag \\
 &= \frac{\beta c^{\overline k}}{2^{\overline k}\pi} \sum_{l=0}^{\overline k} \left(\begin{array}{c}{\overline k}\\ l\end{array}\right)
 e^{({\overline k}-2l)\bar\xi} \int_0^{2\pi} \sin\big(({\overline k}-2l)\eta + {\overline k}\varphi\big)\,\cos(j\eta) \,d\eta \notag\\
 &= \frac{ \beta c^{\overline k}}{2^{\overline k}\pi}\sin({\overline k}\varphi)\sum_{l=0}^{\overline k} \left(\begin{array}{c}{\overline k}\\ l\end{array}\right)
 e^{({\overline k}-2l)\bar\xi} \int_0^{2\pi} \cos\big(({\overline k}-2l)\eta\big)\,\cos(j\eta) \,d\eta,
\end{align}
where the last equality follows  the addition formula for the sine and the mutual orthogonality of trigonometric functions.
As well,  
\begin{equation}\label{eq:bj}
 b_j(\bar\xi) =\frac{ \beta c^{\overline k}}{2^{\overline k}\pi}\cos({\overline k}\varphi)\sum_{l=0}^{\overline k} \left(\begin{array}{c}{\overline k}\\ l\end{array}\right)
 e^{({\overline k}-2l)\bar\xi} \int_0^{2\pi} \sin\big(({\overline k}-2l)\eta\big)\,\sin(j\eta) \,d\eta,
\end{equation}
the computation being  similar to the previous one. We note that the terms in the sums in the right-hand side of \eqref{eq:aj} and \eqref{eq:bj} are nontrivial only if  ${\overline k}-2l=\pm j$, and obtain the values of the coefficients:
\begin{equation*}
	a_j(\bar \xi)=\begin{cases}0&\mbox{ if $\overline k+j$ odd;}\\
	\frac{ \beta c^{\overline k}}{2^{\overline k-1}}\sin({\overline k}\varphi)\left(\begin{array}{c}\overline k\\ \frac{\overline k+j}2\end{array}\right)\cosh j\bar \xi&\mbox{ if $\overline k+j$ even,}\end{cases}
\end{equation*}
and 
\begin{equation*}
	b_j(\bar \xi)=\begin{cases}0&\mbox{ if $\overline k+j$ odd;}\\
	\frac{ \beta c^{\overline k}}{2^{\overline k-1}}\cos({\overline k}\varphi)\left(\begin{array}{c}\overline k\\ \frac{\overline k+j}2\end{array}\right)\sinh \overline k\bar\xi&\mbox{ if $\overline k+j$ even.}\end{cases}
\end{equation*}
Finally,
\begin{align*}
\int_{\R^2\setminus\overline\omega} |\nabla \mathsf{u}_{N,\overline{k}}|^2\, dt&=\sum_{\begin{array}{c}1\le j\le \overline k\\ \cbk \overline k+j\mbox{ even} \cbk\end{array}}\frac{\pi \beta^2c^{2\overline k}}{4^{\overline k-1}}j\left(\begin{array}{c}\overline k\\ \frac{\overline k+j}2\end{array}\right)^2\left(\sin^2\overline k\varphi\cosh^2j\bar\xi+\cos^2\overline k\varphi\sinh^2j\bar\xi\right)\\
&=\sum_{\begin{array}{c}1\le j\le \overline k\\ \cbk \overline k+j\mbox{ even} \cbk \end{array}}\frac{\pi \beta^2c^{2\overline k}}{2^{2\overline k-1}}j\left(\begin{array}{c}\overline k\\ \frac{\overline k+j}2\end{array}\right)^2\left(\cosh 2j\bar\xi-\cos 2\overline k\varphi\right).
\end{align*}
The latter sum can be rewritten to give
\begin{equation*}
	\int_{\R^2\setminus\overline\omega} |\nabla \mathsf{u}_{N,\overline{k}}|^2\, dt=\frac{\pi \beta^2 c^{2{\overline k}}}{2^{2{\overline k}}} \sum_{j=0}^{\overline k}  \left\vert{\overline k}-2j\right\vert\,\left(\begin{array}{c}{\overline k}\\ j\end{array}\right)^2\left(e^{2({\overline k}-2j)\bar\xi}-\cos 2\overline k\varphi\right).
\end{equation*}
In accordance with \cite{AbFeHiLe, AFL2019}, we use the notation
\begin{equation*}
	C_{\overline k}:= \frac1{2^{2\overline k-1}}\sum_{j=0}^{\overline k}  \left\vert{\overline k}-2j\right\vert\,\left(\begin{array}{c}{\overline k}\\ j\end{array}\right)^2= \frac1{4^{\overline k-1}}\sum_{j=0}^{\left\lfloor\frac{\overline k-1}2\right\rfloor}  \left({\overline k}-2j\right)\,\left(\begin{array}{c}{\overline k}\\ j\end{array}\right)^2.
\end{equation*}
Furthermore, we define
\begin{equation*}
	D_{\overline k}(\bar\xi):=\frac1{2^{2{\overline k}}} \sum_{j=0}^{\overline k}  \left\vert{\overline k}-2j\right\vert\,\left(\begin{array}{c}{\overline k}\\ j\end{array}\right)^2\,e^{2({\overline k}-2j)\bar\xi}.
\end{equation*}
We summarize the analysis of this subsection in the following statement.
\begin{prop}\label{p:firstcontr}
 Let $\mathsf{u}_{N,\overline{k}}$ be the unique $C^{1,\alpha}_{\mathrm{loc}}(\R^2\setminus \omega)$ solution to Problem \eqref{eq:uNbf}. Then 
 \[  \int_{\R^2\setminus\overline\omega} |\nabla \mathsf{u}_{N,\overline{k}}|^2\, dt
 = -\frac{\pi\beta^2 c^{2\overline k}}2C_{\overline k}\cos 2\overline k\varphi(u_N,\theta)+\pi \beta^2 c^{2{\overline k}}D_{\overline k}(\bar\xi)
 \]
for any $\theta \in [0,\pi/2]$, with $\varphi(u_N,\theta)$ defined in Remark \ref{rem:phi}.
\end{prop}


\subsection{Computation of the second contribution}

We recall that $u_{N,\#, \overline{k}}$ is a harmonic homogeneous polynomial. We perform an integration by parts, pass to elliptic coordinates, apply the addition formula for sines and thanks to the mutual orthogonality of trigonometric functions we obtain
\begin{align}\label{eq:secondcontr}
 &\int_{\omega} |\nabla u_{N,\#, \overline{k}}|^2\, dt = 
 \int_{\partial\omega} \dfrac{\partial u_{N,\#, \overline{k}}}{\partial \nu_\omega} u_{N,\#, \overline{k}} \, dt 
 =
 \int_{0}^{2\pi}\frac{\partial \Psi_{\overline k}^\varphi}{\partial\xi}(\bar\xi,\eta)\,  \Psi_{\overline k}^\varphi(\bar\xi,\eta)\,d\eta \notag \\
  &= \frac{\beta^2 c^{2{\overline k}}}{2^{2{\overline k}}} \sum_{j,l=0}^{\overline k} \left(\begin{array}{c}{\overline k}\\ j\end{array}\right) \left(\begin{array}{c}{\overline k}\\ l\end{array}\right) ({\overline k}-2j)e^{({\overline k}-2j)\bar\xi+({\overline k}-2l)\bar\xi} \int_0^{2\pi}\sin\big(({\overline k}-2j)\eta+{\overline k}\varphi\big)\sin\big(({\overline k}-2l)\eta+{\overline k}\varphi\big)\,d\eta \notag\\
   &= \frac{ \beta^2 c^{2{\overline k}}}{2^{2{\overline k}}} \sum_{j,l=0}^{\overline k} \left(\begin{array}{c}{\overline k}\\ j\end{array}\right) \left(\begin{array}{c}{\overline k}\\ l\end{array}\right) ({\overline k}-2j)e^{({\overline k}-2j)\bar\xi+({\overline k}-2l)\bar\xi} 
   \left\{ \cos^2({\overline k}\varphi) \int_0^{2\pi} \sin\big(({\overline k}-2j)\eta\big)\,\sin\big(({\overline k}-2l)\eta\big)\,d\eta \right.\notag\\ &\indent\left.+ \sin^2({\overline k}\varphi) \int_0^{2\pi} \cos\big(({\overline k}-2j)\eta\big)\,\cos\big(({\overline k}-2l)\eta\big)\,d\eta \right\}\notag\\
   &= \frac{\pi  \beta^2 c^{2{\overline k}}}{2^{2{\overline k}}} \sum_{j=0}^{\overline k} \left(\begin{array}{c}{\overline k}\\ j\end{array}\right)^2 ({\overline k}-2j)\,e^{2({\overline k}-2j)\bar\xi} 
   \ -\  \frac{\pi}{2^{2{\overline k}}} \sum_{j=0}^{\overline k} \left(\begin{array}{c}{\overline k}\\ j\end{array}\right)
   \left(\begin{array}{c}{\overline k}\\ {\overline k}-j\end{array}\right)
   ({\overline k}-2j) \cos(2{\overline k}\varphi)\notag\\
   &=  \frac{\pi  \beta^2 c^{2{\overline k}}}{2^{2{\overline k}}} \sum_{j=0}^{\overline k} \left(\begin{array}{c}{\overline k}\\ j\end{array}\right)^2 ({\overline k}-2j)\,e^{2({\overline k}-2j)\bar\xi}
\end{align}
where the second to last equality follows from the fact that every term of the sum in the third line is zero except when $l=j$ or $l={\overline k}-j$.
Moreover, the last equality follows easily recalling that $\left(\begin{array}{c}{\overline k}\\ j\end{array}\right) = \left(\begin{array}{c}{\overline k}\\ {\overline k}-j\end{array}\right)$.

\subsection{Comparison with the numerical simulations}

According to Theorem \ref{thm:eig2}, we have 
\begin{align*}
 \lambda_N(\eps,\theta) - \lambda_N &\sim \eps^{2\overline k}\mbox{E}\left(\mathcal R_\theta\mathcal E(a,b),u_{N,\#,\overline{k}}\right) \qquad \text{as }\eps\to0,
\end{align*}
where $\mbox{E}\left(\mathcal R_\theta\mathcal E(a,b),u_{N,\#,\overline{k}}\right)$ is the quantity defined in Equation \eqref{eq:aim}. Summing up the contributions in Proposition \ref{p:firstcontr} and Equation \eqref{eq:secondcontr}, we find
\begin{equation}	
\label{eq:capLim}
	\mbox{E}\left(\mathcal R_\theta\mathcal E(a,b),u_{N,\#,\overline{k}}\right)=-\frac{\pi\beta^2 c^{2\overline k}}2C_{\overline k}\cos 2\overline k\varphi(u_N,\theta)+\pi\beta^2 c^{2\overline k}E_{\overline k}(\bar\xi), 
\end{equation}
where
\begin{align*}
	E_{\overline k}(\bar\xi)&=\frac1{2^{2{\overline k}}} \sum_{j=0}^{\overline k}  \left(\left\vert{\overline k}-2j\right\vert+({\overline k}-2j)\right)\,\left(\begin{array}{c}{\overline k}\\ j\end{array}\right)^2\,e^{2({\overline k}-2j)\bar\xi}\\
	&=\frac1{2^{2{\overline k}-1}}\sum_{j=0}^{\left\lfloor\frac{\overline k-1}2\right\rfloor} ({\overline k}-2j)\,\left(\begin{array}{c}{\overline k}\\ j\end{array}\right)^2\,e^{2({\overline k}-2j)\bar\xi}.
\end{align*}
and $\varphi(u_N,\theta)$ is defined in Remark \ref{rem:phi}. Let us note that the second term in the right-hand side of Equation \eqref{eq:capLim} can be written as a polynomial in $a$ and $b$. Indeed, we have
\begin{equation*}
	\bar \xi = \log\left( \frac{b}{c} + \sqrt{1+\frac{b^2}{c^2}}\right),
\end{equation*}
so that, for any non-negative integer $m$,
\begin{equation*}
 	e^{m\bar\xi}=\left(\frac{b}{c}+\sqrt{1+\frac{b^2}{c^2}}\right)^m=\left(\frac{b}{c}+\sqrt{\frac{c^2+b^2}{c^2}}\right)^m=\left(\frac{a+b}{c}\right)^m.
\end{equation*}
Using the above identity, we get
\begin{equation}	
\label{eq:capLim2}
	\mbox{E}\left(\mathcal R_\theta\mathcal E(a,b),u_{N,\#,\overline{k}}\right)=-\frac{\pi\beta^2 c^{2\overline k}}2C_{\overline k}\cos 2\overline k\varphi(u_N,\theta)+\pi\beta^2 Q_{\overline k}(a,b), 
\end{equation}
with
\begin{align}
\label{eq:poly}
	Q_{\overline k}(a,b)=c^{2\overline k}E_{\overline k}(\bar\xi)&= \frac1{2^{2{\overline k}-1}}\sum_{j=0}^{\left\lfloor\frac{\overline k-1}2\right\rfloor} ({\overline k}-2j)\,\left(\begin{array}{c}{\overline k}\\ j\end{array}\right)^2\,c^{4j}\,(a+b)^{2({\overline k}-2j)}\notag\\
	&=\frac1{2^{2{\overline k}-1}}\sum_{j=0}^{\left\lfloor\frac{\overline k-1}2\right\rfloor} ({\overline k}-2j)\,\left(\begin{array}{c}{\overline k}\\ j\end{array}\right)^2\,(a^2-b^2)^{2j}\,(a+b)^{2({\overline k}-2j)}.
\end{align}

Formula \eqref{eq:capLim2} confirms the simulations on Figure \ref{fig:plot1}, which correspond to a vanishing order ${\overline k}=1$, where $\varphi(u_2,\theta)=\theta-\pi/2$ for $\theta\in\cbk]0,\pi/2]\cbk$,  $\varphi(u_2,0)=\pi/2$ and  $\varphi(u_3,\theta)=\theta$. It also confirms the simulations on Figure \ref{fig:plot2},  corresponding to a vanishing order ${\overline k}=2$; there, for $N=11$ and $N=15$,  $\varphi(u_N,\theta)=\theta$ when $\theta\in[0,\pi/4]$ and $\varphi(u_N,\theta)=\theta-\pi/2$ when $\theta\in\cbk]\pi/4,\pi/2]$. We have thus explained the variations of the functions $\theta\mapsto  \lambda_N(\eps,\theta) - \lambda_N$.

\cbk Finally, starting from Formula \eqref{eq:capLim2}, we can recover the $u$-capacity of a disk and that of a segment, given respectively in Theorems 1.13 and 1.9 of \cite{AbFeHiLe}. We achieve this by letting either $b$ go to $a$ or $b$ go to $0$ and by a suitable scaling.\cbk

{\section*{Acknowledgement}
L.~Abatangelo and P.~Musolino are members of the `Gruppo Nazionale per l'Analisi Matematica, la Probabilit\`a e le loro Applicazioni' (GNAMPA) of the `Istituto Nazionale di Alta Matematica' (INdAM). }

\end{document}